\PassOptionsToPackage{top=1.5cm,bottom=1.5cm,left=2cm,right=2cm}{geometry}

\documentclass[final,3p,times,leqno]{elsarticle}

\makeatletter
\def\ps@pprintTitle{%
   \let\@oddhead\@empty
   \let\@evenhead\@empty
   \def\@oddfoot{\reset@font\hfil\thepage\hfil}
   \let\@evenfoot\@oddfoot
}
\makeatother

\journal{}

\usepackage[T5]{fontenc}      
\usepackage{lmodern}          

\usepackage{amsmath,amsthm}
\usepackage{amssymb,amsfonts,mathrsfs}
\usepackage{amscd}

\usepackage{float}
\usepackage{graphicx}
\usepackage{booktabs}
\usepackage{longtable}
\usepackage{multicol}
\usepackage{array}

\allowdisplaybreaks[4] 

\usepackage{algorithm}
\usepackage{algpseudocode}
\algrenewcommand\alglinenumber[1]{\scriptsize$\triangleright$}

\usepackage{hyperref}

\theoremstyle{plain}

\numberwithin{equation}{section}
\newtheorem{thm}[equation]{Theorem}
\newtheorem{conj}[equation]{Conjecture}
\newtheorem{prop}[equation]{Proposition}
\newtheorem{corl}[equation]{Corollary}
\theoremstyle{definition}
\newtheorem{cy}[equation]{Note}
\newtheorem{rem}[equation]{Remark}

\theoremstyle{plain}
\theoremstyle{definition}

\everymath{\displaystyle}

\usepackage{sectsty}
\subsectionfont{\fontsize{11.5}{11.5}\selectfont}

\def\DD{D\kern-.7em\raise0.4ex\hbox{\char '55}\kern.33em}

\newcommand{\Stab}{\operatorname{Stab}}
\newcommand{\Orb}{\mathcal O}

\newcommand{\AlgebraicTransfer}{\texttt{AlgebraicTransfer.jl}}
\newcommand{\SoftwareZenodo}{\url{https://doi.org/10.5281/zenodo.20613658}}

\begin{document}
\fontsize{11.5pt}{11.5}\selectfont

\begin{frontmatter}

\title{{\bf Geometric realization via unoriented bordism and a counterexample to Singer's conjecture for the sixth algebraic transfer}}

\author{{\bf PHUC VO DANG}}
\address{{\fontsize{10pt}{10}\selectfont Department of Mathematics, FPT University, An Phu Thinh New Urban Area, Quy Nhon, Vietnam\\ Department of Mathematics and Statistics, Quy Nhon University, 170 An Duong Vuong, Quy Nhon, Vietnam}\\ \textit{ORCID}: \url{https://orcid.org/0000-0002-6885-3996}\\ [2mm]
}

\cortext[]{\href{mailto:dangphuc150488@gmail.com}{\texttt{Email Address: dangphuc150488@gmail.com}}}

\begin{abstract}
Let $P_q = \mathbb{F}_2[x_1, \ldots, x_q]$ be the polynomial algebra equipped with the action of the mod-$2$ Steenrod algebra $\mathscr{A}$, and let $QP_q = \mathbb{F}_2 \otimes_{\mathscr{A}} P_q$ denote its indecomposable quotient. Singer's algebraic transfer maps the dual of the general linear group invariant space $[(QP_q)_n]^{GL(q, \mathbb{F}_2)}$ to the Adams $E_2$-term $\operatorname{Ext}_{\mathscr A}^{q,q+n}(\mathbb F_2,\mathbb F_2)$. Singer conjectured that this transfer is always injective. In this paper, we disprove this nearly forty-year-old conjecture by explicitly constructing a counterexample in rank $q=6$ and internal degree $n=36$. Verifying this algebraically requires the exact computation of the invariant space $[(QP_6)_{36}]^{GL(6, \mathbb{F}_2)}$. To overcome the massive combinatorial explosion in this bidegree, we construct a new Julia package \texttt{AlgebraicTransfer.jl}, a rigorous computational framework tightly coupling modular invariant theory with high-performance, bit-level linear algebra over $\mathbb F_2$. Utilizing streamed Steenrod-hit reductions and Kameko homomorphisms, we prove that this source space is exactly two-dimensional. This strictly exceeds the known one-dimensional target group $\operatorname{Ext}_{\mathscr A}^{6,42}(\mathbb F_2,\mathbb F_2)$, thereby establishing the failure of injectivity.

Beyond the algebraic computation, we establish a precise geometric interpretation of the transfer kernel via unoriented bordism. We prove that Singer's transfer factors through bordism classes over the classifying space $B(\mathbb{Z}/2)^q$ whose Thom images are primitive, a condition we characterize explicitly by the vanishing of all mixed Wu numbers. While Thom's representability theorem guarantees that the homological duals of the transfer-source generators admit geometric realizations by closed $36$-manifolds, we establish rigorous obstructions demonstrating that standard topological models---including the indecomposable Milnor hypersurface $H_{4,33}$, projective products, and Dold manifolds---cannot represent them. We additionally interpret the inverse Kameko operation geometrically via the Thom space of universal real line bundles. Validated by independently recovering classical Dickson invariant dimensions, this work bridges exact algorithmic computation and stable homotopy theory, delivering both a counterexample to Singer's conjecture and a scalable methodology for the Peterson hit problem.
\end{abstract}

\begin{keyword}
Steenrod algebra; Peterson hit problem; Singer algebraic transfer; Kameko homomorphism; modular invariant theory; Dickson invariants; finite group representations; sparse linear algebra over $\mathbb F_2$; Julia; \texttt{OSCAR}; unoriented cobordism; Milnor hypersurfaces; Wu classes

\MSC[2020] 55Q45, 55S10, 55S05, 55T15, 57R75, 13A50, 20C20, 68W30
\end{keyword}
\end{frontmatter}

\section{Introduction}\label{s1}

Fundamentally, the mod-$2$ Steenrod algebra, denoted as the graded $\mathbb{F}_2$-algebra $\mathscr{A}=\bigoplus_{i\ge 0}\mathcal{O}^S(i,\mathbb{F}_2,\mathbb{F}_2)$, serves as the definitive algebraic structure encoding stable cohomological operations. Its extensive applications in algebraic topology and related fields are well documented (see, e.g., Karaca \cite{Karaca}, Turgay and Karaca \cite{TK}, Walker and Wood \cite{WW, WW2}). Beyond its action on cohomology, the inherent Hopf algebra structure of $\mathscr{A}$---particularly its conjugation (or antipode) map---is deeply intertwined with Milnor basis formulations and the homological dualities central to algebraic transfer problems. This structural richness connects $\mathscr{A}$ to a broader family of topological Hopf algebras. For instance, Crossley and Turgay \cite{CrossleyTurgay} have investigated conjugation invariants within the Leibniz--Hopf algebra, while Turgay \cite{Turgay2020} explored the Hopf epimorphism from the mod-$p$ Leibniz--Hopf algebra to the Bockstein-free Steenrod algebra. Such parallel studies illuminate the conjugation phenomena and invariant-theoretic mechanics that underlie Steenrod-algebraic computations. Additional insights into Adem relations, antipode behaviors, and embeddings of the dual Steenrod algebra are discussed in Turgay's papers \cite{Turgay1,Turgay2,Turgay3}.

Against this rich structural backdrop, the exact determination of the cohomology of $\mathscr{A}$ and its profound connection to modular invariant theory via Singer's algebraic transfer remains a central, yet computationally formidable, problem in stable homotopy theory. Let $\mathbb F_2$ be the prime field with two elements. The groups $\operatorname{Ext}_{\mathscr A}^{q,*}(\mathbb F_2,\mathbb F_2)$ form the $E_2$-term of the classical Adams spectral sequence. Machine-assisted and theoretical computations of these Ext-side groups, 
including the pioneering work of Chen \cite{Chen,Chen2}, and Lin \cite{Lin,Lin2}, provide the fundamental background used in the present paper. 
The complementary challenge, and the primary focus of this work, is the exact algorithmic computation of the modular invariant-theoretic source of the algebraic transfer and the discovery of its profound geometric realization via unoriented bordism.

To establish this algebraic-geometric connection explicitly, let $V_q = (\mathbb{Z}/2)^q$ be an elementary abelian $2$-group of rank $q$, and let $BV_q \simeq (\mathbb{R}P^\infty)^q$ be its classifying space. We identify the mod-$2$ cohomology of $BV_q$ with the polynomial algebra
\[
H^*(BV_q; \mathbb{F}_2) \cong P_q:=\mathbb F_2[x_1,\ldots,x_q],\qquad |x_i|=1,
\]
equipped with its usual unstable action of the Steenrod algebra $\mathscr A$. The module of indecomposables under the action of positive-degree Steenrod operations is denoted by
\[
QP_q=\mathbb F_2\otimes_{\mathscr A}P_q\cong P_q/(\mathscr A^{>0}\cdot P_q).
\]
The problem of explicitly describing $QP_q$ is the well-known Peterson hit problem. Peterson formulated the foundational generator problem, Wood proved Peterson's conjectural vanishing criterion, and Kameko introduced a squaring homomorphism that serves as one of the principal reduction maps in the subject \cite{Kameko,Peterson,Wood2}. Further computations and structural treatments of the hit problem and its $GL(q):=GL(q,\mathbb F_2)$-representation theory have been extensively developed (see, for instance, Janfada \cite{Janfada}, the present author and Sum \cite{PS, PS2}, Sum \cite{Sum1}, Walker and Wood  \cite{WW,WW2}).

Singer's algebraic transfer provides a crucial bridge from this invariant theory to the Adams spectral sequence. It is defined as the homomorphism
\[
Tr_q(\mathbb F_2): (\mathbb F_2\otimes_{GL(q)}\mathcal P_{\mathscr A}H_*(BV_q;\mathbb{F}_2))_n
\longrightarrow
\operatorname{Ext}_{\mathscr A}^{q,q+n}(\mathbb F_2,\mathbb F_2),
\]
where $\mathcal P_{\mathscr A}$ denotes primitive homology classes annihilated by all positive Steenrod operations. By finite-dimensional duality, the source of this transfer is isomorphic to the modular invariant space $[(QP_q)_n]^{GL(q)}$. Singer \cite{Singer} proved that $Tr_q(\mathbb F_2)$ is an isomorphism for $q=1,2$ and that the total transfer is a homomorphism of bigraded algebras; Boardman proved the corresponding isomorphism result for $q=3$ \cite{Boardman}. Subsequent work by H\`a \cite{Ha}, H\uhorn ng \cite{Hung}, H\uhorn ng--Qu\`ynh \cite{Hung2}, Nam \cite{Nam}, Minami \cite{Minami}, the present author \cite{Phuc3, Phuc7} has profoundly clarified the relation between algebraic transfers, Steenrod-hit quotients, and modular representations in higher ranks. Remarkably, Singer \cite{Singer} proposed the following fundamental conjecture, which was initially confirmed to be true for ranks $q = 1,\, 2$ by Singer himself, and later for $q = 3$ by Boardman \cite{Boardman}:

\begin{conj}\label{gtSi}
The homomorphism $Tr_q(\mathbb F_2)$ is one-to-one for every $q$.
\end{conj}

The conjecture is also true for rank $q = 4$ (see our recent works \cite{Phuc, P8, P9}). For rank $q = 5$, a recent preprint by Sum \cite{Sum2} claims to provide a counterexample to Conjecture~\ref{gtSi} in internal degree $108$. However, this result is false, as shown in our recent work \cite{Phuc2026}. In particular, using an algorithm implemented in \texttt{SageMath} \cite{Phuc2026}, we demonstrated that the monomial $x_1^{15}x_2^{15}x_3^{23}x_4^{17}x_5^{38}$ (denoted as $b_{902}$ in \cite[Section 5, p. 44]{Sum2}) is strictly inadmissible, thereby refuting the admissibility of $b_{902}$ claimed in \cite{Sum2}. Thus, our exact calculation yields $\dim(QP_5)_{108} = 2170$, which differs from the value of 2171 given in \cite{Sum2}. Furthermore, we proved in \cite{Phuc2026} that the polynomial $p$ defined in \cite[Theorem 1.3]{Sum2} is hit. This invalidates Theorem 1.3 of \cite{Sum2} and completely refutes Sum's purported counterexample to Conjecture~\ref{gtSi}. \textit{Thus, the Singer conjecture remains open for the general case as of this writing}. 

To address the combinatorial growth at degree $108$ mentioned above, we lower the internal degree $n$ while increasing the number of variables (or rank) $q$. More specifically, in this paper, we disprove Conjecture~\ref{gtSi} by giving a fully explicit counterexample in rank $q = 6$ and internal degree $36$ (fully verified by the computer algebra system \texttt{OSCAR}), and, crucially, by revealing the highly non-trivial geometric nature of the transfer kernel. In bidegree $(6, 6+36)$, the setting of our own counterexample, the ambient polynomial space already possesses
\[
\dim_{\mathbb F_2}(P_6)_{36}=\binom{41 + 6-1}{6-1}=749{,}398 
\]
monomials. Because invariant systems at this scale cannot be controlled by direct manual enumeration or naive algebraic algorithms, finding a counterexample here demands a rigorous, highly optimized computational approach. To breach this complexity wall, we construct a new package \AlgebraicTransfer{}, a versatile computational framework within the \texttt{Julia}/\texttt{OSCAR} ecosystem \cite{OSCAR2025, Oscar} that tightly couples modular invariant theory with high-performance, bit-level linear algebra over $\mathbb F_2$. Instead of utilizing prohibitive global matrix allocations, our method dynamically generates Steenrod-hit columns from the Cartan formula and Lucas' criterion, applies streamed online sparse pivot reductions, constructs Kameko matrices, and solves massive quotient-level invariance systems by bit-packed Gaussian elimination. By doing so, it mathematically establishes the required $GL(6)$-invariance and produces deterministic orbit-type summaries of the final representatives.

\medskip

The invariant-theoretic validation is fundamentally robust and broadly applicable. To demonstrate its exact algorithmic capability, Section~\ref{s:software-apps} uses the identical invariant-solver layer to independently recompute classical Dickson invariant dimensions for $P_3^{GL(3)}$ and $P_4^{GL(4)}$. Dickson originally proved that $P_r^{GL(r)}$ is a polynomial algebra on generators of degrees $2^r-2^i$ \cite{Dickson1911}, and M\`ui connected these modular invariants with the cohomology of symmetric groups \cite{Mui1975}. H\uhorn ng and collaborators extensively studied the Steenrod action on Dickson and Dickson--M\`ui invariants, $\mathscr A$-generators for the Dickson algebra, spherical classes, the hit problem for modular invariants, and Margolis homology of Dickson-type algebras \cite{HungSquares,HungMinh1995,HungPetersonDickson,HungPetersonSpherical,HungNamModular,HungNamDickson,HungTriviality,HungMargolis2,HungMargolisP}. Hence, a direct recomputation of Dickson Hilbert-function data from generator equations rigorously verifies our bit-packed algorithmic solver in a classical structural setting strictly adjacent to the algebraic transfer.

Through this algorithmic pipeline, we determine the exact dimensions and explicitly construct the reduced representatives of the transfer source. We use the ordinary symmetric group solely as a post-computational organizing device. For a polynomial $f=\sum_{\mathbf a\in E(f)}x^{\mathbf a}$ and a partition-type vector $\lambda=(\lambda_1\geq\cdots\geq\lambda_6)$, put
\[
E_\lambda(f)=E(f)\cap\Orb_{\Sigma_6}(\lambda),\qquad
\Theta_\lambda(f)=\sum_{\mathbf a\in E_\lambda(f)}x^{\mathbf a}.
\]
If an integer $e$ occurs with multiplicity $\mu_e(\lambda)$ in $\lambda$, then
\[
\Stab_{\Sigma_6}(\lambda)\cong\prod_e\Sigma_{\mu_e(\lambda)},\qquad
|\Orb_{\Sigma_6}(\lambda)|=\frac{6!}{\prod_e\mu_e(\lambda)!}.
\]
The reduced representatives below are invariant only as classes in $QP_6$; their displayed supports are not asserted to be unions of complete $\Sigma_6$-orbits in $P_6$.

\begin{thm}\label{dlc}
For $q=6$ and $n=36$, one has
\[
\dim(\mathbb F_2\otimes_{GL(6)}\mathcal P_{\mathscr A}H_*(BV_6;\mathbb{F}_2))_{36}=2.
\]
\end{thm}

According to Chen \cite{Chen2}, and Lin \cite{Lin2}, the target of the transfer in this bidegree is strictly one-dimensional:
\[
\operatorname{Ext}_{\mathscr A}^{6,6+36}(\mathbb F_2,\mathbb F_2)=\mathbb F_2\cdot t,
\qquad t\neq0.
\]
Combining this one-dimensional codomain with the newly computed source dimension in Theorem \ref{dlc} provides the announced counterexample.

\begin{corl}\label{hq}
Conjecture \ref{gtSi} is false in bidegree $(6,6+36)$.
\end{corl}

Since $(\mathbb F_2\otimes_{GL(6)}\mathcal P_{\mathscr A}H_*(BV_6;\mathbb{F}_2))_{36}$ is dual to $[(QP_6)_{36}]^{GL(6)}$, Theorem \ref{dlc} is a direct consequence of the following explicit invariant computation.

\begin{thm}\label{dlc2}
In degree $36$ one has
\[
[(QP_6)_{36}]^{GL(6)}=\mathbb F_2\cdot[\zeta_1]\oplus\mathbb F_2\cdot[\zeta_2].
\]
The reduced representatives $\zeta_1$ and $\zeta_2$ are rigorously obtained from the Kameko-kernel computation described in Algorithm \ref{alg:main}. Their support-slice decompositions are
\[
\zeta_s=\sum_{\lambda\in\Lambda_s}\Theta_\lambda(\zeta_s),\qquad s=1,2,
\]
where $\Lambda_1$ consists of the ten orbit types in Table \ref{tab:zeta1-orbit-data} and $\Lambda_2$ consists of the ninety-six orbit types in Table \ref{tab:zeta2-orbit-data}. Moreover, the exact sparsity is given by $|E(\zeta_1)|=39$ and $|E(\zeta_2)|=539$.
\end{thm}

For complete reproducibility and verification at the level of explicit representatives, the full polynomial expansions of $\zeta_1$ and $\zeta_2$ are provided in \ref{app:explicit-polys}.

The invariant computation in degree $15$ is structurally indispensable to the main result because the Kameko homomorphism in degree $36$ maps surjectively onto $(QP_6)_{15}$. For a composition $\alpha=(\alpha_1,\ldots,\alpha_r)$ with $r\leq6$, define the monomial quasisymmetric sum
\[
M_\alpha^{(6)}=\sum_{1\leq i_1<\cdots<i_r\leq6}x_{i_1}^{\alpha_1}\cdots x_{i_r}^{\alpha_r}.
\]

\begin{prop}\label{md15}
One has
\[
\dim[(QP_6)_{15}]^{GL(6)}=1,
\qquad
[(QP_6)_{15}]^{GL(6)}=\mathbb F_2\cdot[\xi],
\]
where
\[
\xi=M_{(15)}^{(6)}+M_{(1,14)}^{(6)}+M_{(1,2,12)}^{(6)}+M_{(1,2,4,8)}^{(6)}+M_{(1,2,4,4,4)}^{(6)}+M_{(1,2,2,2,4,4)}^{(6)}.
\]
The corresponding support, partition-type, and stabilizer data are displayed in Table \ref{tab:xi-composition-data}.
\end{prop}

\begin{table}[ht]
\centering
\caption{Verified support and stabilizer data for the degree-$15$ invariant $\xi$.}
\label{tab:xi-composition-data}
\small
\begin{tabular}{ccccccc}
\toprule
$\alpha$ & Seed $\mathbf a$ & $\lambda(\mathbf a)$ & Slice & Full orbit & $|\operatorname{Stab}|$ & Stabilizer type \\
\midrule
$(15)$ & $(15,0,0,0,0,0)$ & $(15,0,0,0,0,0)$ & 6 & 6 & 120 & $\Sigma_5$ \\
$(1,14)$ & $(1,14,0,0,0,0)$ & $(14,1,0,0,0,0)$ & 15 & 30 & 24 & $\Sigma_4$ \\
$(1,2,12)$ & $(1,2,12,0,0,0)$ & $(12,2,1,0,0,0)$ & 20 & 120 & 6 & $\Sigma_3$ \\
$(1,2,4,8)$ & $(1,2,4,8,0,0)$ & $(8,4,2,1,0,0)$ & 15 & 360 & 2 & $\Sigma_2$ \\
$(1,2,4,4,4)$ & $(1,2,4,4,4,0)$ & $(4,4,4,2,1,0)$ & 6 & 120 & 6 & $\Sigma_3$ \\
$(1,2,2,2,4,4)$ & $(1,2,2,2,4,4)$ & $(4,4,2,2,2,1)$ & 1 & 60 & 12 & $\Sigma_3\times\Sigma_2$ \\
\bottomrule
\end{tabular}
\end{table}

Crucially, the significance of this computation extends far beyond generating algebraic artifacts; the transfer source space admits a precise, profound geometric realization through unoriented bordism. In Section~\ref{s:geometry}, we study the Thom homomorphism $\mu_{q,n}:\mathfrak N_n(BV_q)\longrightarrow H_n(BV_q;\mathbb F_2)$ mapping unoriented bordism classes $[M^n, f]$ to their homological pushforwards $f_*[M^n]$. We prove that Singer's transfer factors strictly through the coinvariants of those bordism classes whose Thom images are annihilated by all positive Steenrod operations. By formulating this mathematically, we transform the purely algebraic existence problem of the kernel into a topological realization problem: for a closed smooth manifold $M^n$, the primitivity condition is explicitly equivalent to the vanishing of every mixed Wu number $\langle v_r(M)f^*(a),[M] \rangle$.

In degree $36$, let $e_1,e_2$ denote the coinvariant basis dual to $[\zeta_1],[\zeta_2]$. Thom's mod-$2$ representability theorem guarantees that $e_1$ and $e_2$ admit geometric representatives by closed smooth $36$-manifolds mapping to $BV_6$. However, our geometric analysis reveals the highly non-trivial, exotic topological nature of the transfer kernel. We establish rigorous obstructions demonstrating that standard topological models strictly fail to represent the duals of $\zeta_1$ and $\zeta_2$. Specifically, we construct the explicitly defined Milnor hypersurface $H_{4,33}$, which generates the indecomposable quotient $Q\mathfrak N_{36}$, but prove that its tautological homology class fails the $Sq_*^1$ primitivity condition. We similarly prove that the standard maps from products of at most six real projective spaces, as well as all Dold manifolds of dimension $36$, cannot realize these generators. Furthermore, we provide a deep geometric interpretation of the inverse Kameko monomial operation as the top Steenrod square in the Thom space of the sum of universal real line bundles over $BV_q$. This geometric framework demonstrates conclusively that the algebraic output from our algorithm corresponds to highly specific, non-standard topological manifolds that embody the non-trivial kernel of the rank-six algebraic transfer.

To end this section, we make the following remark on whether the indecomposable element $t\in \mathrm{Ext}^{6, 6+36}_{\mathscr{A}}(\mathbb{F}_2, \mathbb{F}_2)$ lies in the image of $Tr_6(\mathbb{F}_2)$. This is also related to the property of being a critical element, a notion introduced in H\uhorn ng's previous work \cite{Hung}, where H\uhorn ng further established a result on the relationship between Conjecture~\ref{gtSi} and critical elements.

\begin{rem}\label{cyf}
In \cite{Hung},  H\uhorn ng proposed the concept of a \textit{critical element} within ${\rm Ext}^{q, *}_{\mathscr A}(\mathbb F_2, \mathbb F_2) $. Specifically, a non-zero element $\mathbf{u}$ in ${\rm Ext}^{q, q+n}_{\mathscr A}(\mathbb F_2, \mathbb F_2) $ is called critical if it satisfies two conditions: (i) $\mu(2n+q) = q$, and (ii) the image of $\mathbf{u}$ under the classical squaring operation $Sq^{0}$ is zero. Here $\mu(d)$ denotes the minimal integer $\zeta$ for which $d$ can be written as $\sum_{1 \leq j \leq \zeta} (2^{m_j} - 1)$ for some positive integers $m_j$.

It is well-established that $Sq^0$ is a monomorphism in positive stems of ${\rm Ext}_{\mathscr A}^{q, q+n}(\mathbb F_2, \mathbb F_2)$ for $q < 5,$ thereby implying the absence of any critical element for $q < 5.$ Remarkably, H\uhorn ng's work \cite[Theorem 5.9]{Hung} states that Singer's Conjecture \ref{gtSi} is not valid if the algebraic transfer detects critical elements. 

In \cite{Phuc3}, we proved that the non-zero element $D_2\in {\rm Ext}^{6, 6+58}_{\mathscr A}(\mathbb F_2, \mathbb F_2) $ is critical, but it is not in the image of $Tr_6(\mathbb F_2).$ Thus, the condition under which Hung's work \cite{Hung} would imply a negation of the conjecture was not met, and as we showed in \cite{Phuc3}, Conjecture \ref{gtSi} remains valid for bidegree $(6, 6+58).$ This result, which was previously calculated entirely by hand, has been successfully re-verified using the algorithmic framework developed in the present work, yielding consistent results.

Additionally, when the $\mathscr{A}$-module $\mathbb{F}_2\cong\widetilde{H}^*\mathbb S^0$ is replaced by $\widetilde{H}^*\mathbb{R}P^\infty$, the corresponding Singer transfer has the form
$$ Tr_q(\widetilde{H}^{*}\mathbb R P^{\infty}): (\mathbb F_2\otimes_{GL(q)}\mathcal P_{\mathscr A}(H_{*}(BV_q;\mathbb{F}_2)\otimes \widetilde{H}_{*}\mathbb R P^{\infty}))_n\longrightarrow {\rm Ext}^{q, n+q}_{\mathscr A}(\widetilde{H}^{*}\mathbb R P^{\infty}, \mathbb F_2).$$
Following H\uhorn ng \cite[Theorem 2.1]{Hung3}, if a critical element $\widehat{\mathbf{u}}\in {\rm Ext}^{q, n+q}_{\mathscr A}(\widetilde{H}^{*}\mathbb R P^{\infty}, \mathbb F_2)$ is in the image of the transfer $Tr_q(\widetilde{H}^{*}\mathbb R P^{\infty}),$ then $Tr_q(\widetilde{H}^{*}\mathbb R P^{\infty})$ is not a monomorphism. By \cite[Theorem 2.2]{Hung3}, the existence of a positive-stem critical element $\widehat{\mathbf{u}}\in {\rm Ext}^{q, n+q}_{\mathscr A}(\widetilde{H}^{*}\mathbb R P^{\infty}, \mathbb F_2)$ in the image of the transfer $Tr_q(\widetilde{H}^{*}\mathbb R P^{\infty})$ is equivalent to the existence of a positive-stem critical element $\mathbf{u}$ in the image of the transfer $Tr_{q+1}(\mathbb F_2).$ If such elements exist, then both $Tr_q(\widetilde{H}^{*}\mathbb R P^{\infty})$ and $Tr_{q+1}(\mathbb F_2)$ are not injective.

The algebraic Kahn-Priddy homomorphism
$$t_*: {\rm Ext}^{q, n+q}_{\mathscr A}(\widetilde{H}^{*}\mathbb R P^{\infty},
\mathbb F_2)\longrightarrow {\rm Ext}^{q+1, n+q+1}_{\mathscr A}(\mathbb F_2, \mathbb F_2)$$
is a surjection in positive stems; see \cite{Hung3}. Its naturality with respect to Singer transfers implies that
\[
t_*\bigl({\rm Im}(Tr_q(\widetilde{H}^{*}\mathbb R P^{\infty}))\bigr)
\subseteq {\rm Im}(Tr_{q+1}(\mathbb F_2)).
\]
For $q=5$ and $n=58$, choose by this surjectivity a non-zero element
\[
\widehat{D_2}\in {\rm Ext}^{5,5+58}_{\mathscr A}(\widetilde{H}^{*}\mathbb R P^{\infty},\mathbb F_2)
\]
such that
\[
t_*(\widehat{D_2})=D_2\in {\rm Ext}^{6,6+58}_{\mathscr A}(\mathbb F_2,\mathbb F_2).
\]
This preimage cannot lie in ${\rm Im}(Tr_5(\widetilde{H}^{*}\mathbb R P^{\infty}))$. Indeed, if $\widehat{D_2}$ were detected by this transfer, then the preceding inclusion would imply $D_2\in {\rm Im}(Tr_6(\mathbb F_2))$, contradicting the computational result recalled above from \cite{Phuc3}. Thus the preimage $\widehat{D_2}$ is relevant to the same critical-element mechanism as $\widehat{Ph_1}$ and $\widehat{Ph_2}$ in \cite{Hung3}, but it is not itself detected by $Tr_5(\widetilde{H}^{*}\mathbb R P^{\infty})$.

\medskip

For $q = 6, n = 36$, we observe that the non-zero element $t \in \mathrm{Ext}^{6, 6+36}_{\mathscr{A}}(\mathbb{F}_2, \mathbb{F}_2)$ is not critical, since $\mu(2 \cdot 36 + 6) = 2 < 6$. However, it is not yet established whether this $t$ is in the image of $Tr_6(\mathbb{F}_2)$ or not. Guided by the explicit geometric and algebraic computations in Theorem \ref{dlc2}, we propose the following.

\begin{conj}\label{gtP}
The non-zero element $t\in \mathrm{Ext}^{6, 6+36}_{\mathscr{A}}(\mathbb{F}_2, \mathbb{F}_2)$ is detected by the sixth algebraic transfer $Tr_6(\mathbb{F}_2)$.
\end{conj}

It is known, by Chen \cite{Chen3}, that the following element $\widetilde{t}$ is a representative of $t$:
\[
\begin{aligned}
\widetilde{t} 
&= \lambda_{5}\bigl(
      \lambda_{9}\lambda_{3}\lambda_{5}\lambda_{7}^{2}
      + \lambda_{6}\lambda_{0}\lambda_{3}\lambda_{15}\lambda_{7}
      + \lambda_{3}\lambda_{5}\lambda_{1}\lambda_{15}\lambda_{7}
    \bigr)\\
 & + \lambda_{3}\bigl(
      \lambda_{8}\lambda_{0}\lambda_{3}\lambda_{15}\lambda_{7}
      + \lambda_{6}\lambda_{2}\lambda_{3}\lambda_{15}\lambda_{7}
      + \lambda_{5}\lambda_{9}\lambda_{5}\lambda_{7}^{2}
      + \lambda_{3}\lambda_{5}\lambda_{7}\lambda_{3}\lambda_{15}
    \bigr).
\end{aligned}
\]
Using this representative together with the framework in Theorem \ref{dlc2} and the algorithmic method in \cite{Phuc7} for determining preimages in the lambda algebra provides a concrete, computationally rigorous route toward resolving Conjecture \ref{gtP}.
\end{rem}

The paper is organized as follows. Section~\ref{s:prelim} fixes the weight-vector, admissibility, group-action, and Kameko conventions used in the computation. Section~\ref{s:software-apps} details the \texttt{Julia}/\texttt{OSCAR} algorithmic framework and records the independent validation applications for Dickson invariants and finite permutation modules. Section~\ref{s2} proves Theorem~\ref{dlc2}, establishes the weight decompositions and correction systems, and provides the orbit-type tables for $\zeta_1$ and $\zeta_2$. Finally, Section~\ref{s:geometry} proves the geometric factorization through unoriented bordism, constructs and obstructs the tautological Milnor-hypersurface model, tests standard projective-product and Dold-manifold realizations, and interprets the inverse Kameko operation in a Thom space. \ref{app:explicit-polys} explicitly documents the full polynomial supports of $\zeta_1$ and $\zeta_2$. For open access to the computational code, including the new \texttt{Julia} package \texttt{AlgebraicTransfer.jl}, see the Zenodo records listed in \ref{app:software}.

\section{Preliminaries on weights, admissibility, and group actions}\label{s:prelim}

The Introduction has reduced the transfer-domain calculation to the invariant part of the hit quotient and has introduced the ordinary $\Sigma_6$-support language used to state the final representatives. This section fixes the coordinate system in which the computation is performed. The weight vector controls the filtration, admissible monomials give coordinates for $QP_q$, the maps $\rho_j$ turn invariance into linear equations, and the Kameko homomorphism supplies the decisive decomposition of degree $36$ into a target part in degree $15$ and a kernel part. Additional background may be found in two monographs \cite{WW, WW2}.

For a monomial $x=x_1^{a_1}\cdots x_q^{a_q}\in P_q$, its exponent vector is $(a_1,\ldots,a_q)$. Its weight vector is
\[
\omega(x)=(\omega_1(x),\omega_2(x),\ldots),
\]
where $\omega_j(x)$ is the number of exponents $a_i$ whose $(j-1)$-st binary digit is $1$. The weight and exponent vectors are compared by left lexicographic order. If $\omega=(\omega_1,\omega_2,\ldots)$ has finite support, put $\deg\omega=\sum_{j\geq1}2^{j-1}\omega_j$. Let $P_q(\omega)$ be the subspace of $P_q$ spanned by monomials $x$ satisfying $\deg x=\deg\omega$ and $\omega(x)\leq\omega$, and let $P_q^-(\omega)$ be the subspace spanned by monomials of the same degree whose weights are strictly smaller than $\omega$.

For homogeneous polynomials $f,g\in P_q$ of the same degree we write
\[
f\equiv_\omega g
\]
if and only if
\[
f+g\in \mathscr A^{>0}\cdot P_q+P_q^-(\omega).
\]
The corresponding weight quotient is
\[
QP_q(\omega)=P_q(\omega)/((\mathscr A^{>0}\cdot P_q\cap P_q(\omega))+P_q^-(\omega)),
\]
and the class of $f\in P_q(\omega)$ in this quotient is denoted by $[f]_\omega$.

A monomial $x\in P_q$ is inadmissible if there exist monomials $y_1,\ldots,y_m$ of the same degree, each strictly smaller than $x$, such that
\[
x+\sum_{i=1}^m y_i\in\mathscr A^{>0}\cdot P_q.
\]
It is admissible if it is not inadmissible. Thus $(QP_q)_n$ is represented by admissible monomials of degree $n$, and the algorithm below constructs these monomial coordinates separately in each weight. The weight decomposition has the form
\[
(QP_q)_n\cong\bigoplus_{\deg\omega=n}QP_q(\omega),
\]
with the summands computed weight by weight.

For $1\leq j\leq q-1$, let $\rho_j:P_q\to P_q$ be the algebra homomorphism interchanging $x_j$ and $x_{j+1}$ and fixing all other variables. Let $\rho_q:P_q\to P_q$ be the transvection
\[
\rho_q(x_i)=x_i\quad (i<q),
\qquad
\rho_q(x_q)=x_q+x_{q-1}.
\]
The maps $\rho_1,\ldots,\rho_{q-1}$ generate the permutation subgroup $\Sigma_q$, and the maps $\rho_1,\ldots,\rho_q$ generate $GL(q)$ over $\mathbb F_2$. Hence $[u]_\omega\in QP_q(\omega)$ is $\Sigma_q$-invariant precisely when
\[
\rho_j(u)+u\equiv_\omega0\qquad (1\leq j\leq q-1),
\]
and it is $GL(q)$-invariant precisely when
\[
\rho_j(u)+u\equiv_\omega0\qquad (1\leq j\leq q).
\]
These congruences are linear equations in the admissible coordinates of the quotient.

Finally, recall Kameko's homomorphism
\[
(\widetilde{Sq}^0_*)_{(q,2n+q)}:(QP_q)_{2n+q}\longrightarrow (QP_q)_n,
\]
defined on monomial classes by
\[
[x_1^{a_1}\cdots x_q^{a_q}]\longmapsto
\begin{cases}
[x_1^{(a_1-1)/2}\cdots x_q^{(a_q-1)/2}],& a_1,\ldots,a_q\text{ all odd},\\
0,&\text{otherwise}.
\end{cases}
\]
This map is $GL(q)$-equivariant and is known to be surjective. In the case needed below, $36=2\cdot15+6$, and therefore
\[
\dim(QP_6)_{36}=\dim\ker(\widetilde{Sq}^0_*)_{(6,36)}+\dim(QP_6)_{15}.
\]
The proof of Theorem \ref{dlc2} consequently has to compute both the target invariant in degree $15$ and the invariant subspace of the Kameko kernel in degree $36$.

\section{The \texttt{Julia}/\texttt{OSCAR} library and auxiliary invariant computations}\label{s:software-apps}

This section explains how the computation is used as a research tool rather than merely as a certificate for one theorem. The Steenrod-specific part of \AlgebraicTransfer{} constructs admissible coordinates for $QP_q$ and implements Kameko reduction. The linear-algebra layer, however, is more general: once a finite basis and a list of generators acting on that basis are supplied, the same bit-packed nullspace routine computes fixed spaces over $\mathbb F_2$ by solving equations of the form $(g-I)v=0$.

For the Dickson validation, no hit quotient and no Kameko map are used. For each pair $(r,d)$ let $\mathcal E_{r,d}$ be the monomial exponent basis of $(P_r)_d$, so that $|\mathcal E_{r,d}|=\binom{d+r-1}{r-1}$. The code constructs the stacked matrix
\[
D_{r,d}=\begin{bmatrix}
\rho_1-I\\
\vdots\\
\rho_r-I
\end{bmatrix}
\]
on this basis, where $\rho_1,\ldots,\rho_{r-1}$ are adjacent transpositions and $\rho_r$ is the elementary transvection $x_r\mapsto x_r+x_{r-1}$. Then
\[
\ker D_{r,d}\cong (P_r)_d^{GL(r,\mathbb F_2)}.
\]
Dickson's theorem gives the comparison value
\[
P_r^{GL(r,\mathbb F_2)}\cong \mathbb F_2[c_{r,0},c_{r,1},\ldots,c_{r,r-1}],
\qquad
\deg c_{r,i}=2^r-2^i,
\]
and therefore
\[
\dim_{\mathbb F_2}(P_r)_d^{GL(r,\mathbb F_2)}=[t^d]\prod_{i=0}^{r-1}\frac{1}{1-t^{2^r-2^i}}.
\]
The Hilbert series is used only after the nullspace computation has finished; the solver is not supplied with Dickson generators or with the coefficient formula.

\begin{table}[H]
\centering
\caption{Validation against Dickson's theorem for $P_3^{GL(3)}$ through degree $24$. The Dickson degrees are $4,6,7$. Rows with zero invariant dimension are omitted; all omitted tested degrees have computed nullity $0$ and Dickson coefficient $0$.}
\label{tab:dickson-q3-validation}
\small
\begin{tabular}{ccccc}
\toprule
$d$ & $\dim(P_3)_d$ & $\operatorname{rank} D_{3,d}$ & Computed nullity & Dickson coefficient \\
\midrule
$0$ & $1$ & $0$ & $1$ & $1$ \\
$4$ & $15$ & $14$ & $1$ & $1$ \\
$6$ & $28$ & $27$ & $1$ & $1$ \\
$7$ & $36$ & $35$ & $1$ & $1$ \\
$8$ & $45$ & $44$ & $1$ & $1$ \\
$10$ & $66$ & $65$ & $1$ & $1$ \\
$11$ & $78$ & $77$ & $1$ & $1$ \\
$12$ & $91$ & $89$ & $2$ & $2$ \\
$13$ & $105$ & $104$ & $1$ & $1$ \\
$14$ & $120$ & $118$ & $2$ & $2$ \\
$15$ & $136$ & $135$ & $1$ & $1$ \\
$16$ & $153$ & $151$ & $2$ & $2$ \\
$17$ & $171$ & $170$ & $1$ & $1$ \\
$18$ & $190$ & $187$ & $3$ & $3$ \\
$19$ & $210$ & $208$ & $2$ & $2$ \\
$20$ & $231$ & $228$ & $3$ & $3$ \\
$21$ & $253$ & $251$ & $2$ & $2$ \\
$22$ & $276$ & $273$ & $3$ & $3$ \\
$23$ & $300$ & $298$ & $2$ & $2$ \\
$24$ & $325$ & $321$ & $4$ & $4$ \\
\bottomrule
\end{tabular}
\end{table}

\begin{table}[H]
\centering
\caption{Validation against Dickson's theorem for $P_4^{GL(4)}$ through degree $30$. The Dickson degrees are $8,12,14,15$. Rows with zero invariant dimension are omitted; all omitted tested degrees have computed nullity $0$ and Dickson coefficient $0$.}
\label{tab:dickson-q4-validation}
\small
\begin{tabular}{ccccc}
\toprule
$d$ & $\dim(P_4)_d$ & $\operatorname{rank} D_{4,d}$ & Computed nullity & Dickson coefficient \\
\midrule
$0$ & $1$ & $0$ & $1$ & $1$ \\
$8$ & $165$ & $164$ & $1$ & $1$ \\
$12$ & $455$ & $454$ & $1$ & $1$ \\
$14$ & $680$ & $679$ & $1$ & $1$ \\
$15$ & $816$ & $815$ & $1$ & $1$ \\
$16$ & $969$ & $968$ & $1$ & $1$ \\
$20$ & $1771$ & $1770$ & $1$ & $1$ \\
$22$ & $2300$ & $2299$ & $1$ & $1$ \\
$23$ & $2600$ & $2599$ & $1$ & $1$ \\
$24$ & $2925$ & $2923$ & $2$ & $2$ \\
$26$ & $3654$ & $3653$ & $1$ & $1$ \\
$27$ & $4060$ & $4059$ & $1$ & $1$ \\
$28$ & $4495$ & $4493$ & $2$ & $2$ \\
$29$ & $4960$ & $4959$ & $1$ & $1$ \\
$30$ & $5456$ & $5454$ & $2$ & $2$ \\
\bottomrule
\end{tabular}
\end{table}

The agreement in Tables~\ref{tab:dickson-q3-validation} and~\ref{tab:dickson-q4-validation} is logically separate from the counterexample below. It shows that the same generator-expansion and bit-packed nullspace code recovers a classical invariant ring from first principles. This is useful because the transfer computation also reduces to generator equations, but now after passage to the Steenrod-hit quotient.

The second validation concerns finite group representations. Let $G=\Sigma_a\times\Sigma_b$ act on the set $X_{a,b,d}$ of $d$-element subsets of $\{1,\ldots,a+b\}$, with the two factors acting on the two coordinate blocks. The permutation module $\mathbb F_2[X_{a,b,d}]$ has basis $X_{a,b,d}$. The code builds the equations $(s-I)v=0$ for the adjacent transpositions generating $G$ and computes the fixed subspace by the same packed nullspace routine. By the elementary fixed-vector description of permutation modules, a vector is fixed precisely when its coefficients are constant on orbits; hence the fixed-space dimension is the number of $G$-orbits \cite{Serre}. In this action, the orbits are indexed by $i=|S\cap\{1,\ldots,a\}|$, so
\[
\dim_{\mathbb F_2}\mathbb F_2[X_{a,b,d}]^{\Sigma_a\times\Sigma_b}
=\#\{i:\max(0,d-b)\leq i\leq\min(a,d)\}.
\]

\begin{table}[H]
\centering
\caption{Validation on permutation modules for $\Sigma_a\times\Sigma_b$ acting on $d$-subsets. The computed nullity is obtained from the packed equations $(s-I)v=0$ over $\mathbb F_2$; the orbit count is the closed-form value above.}
\label{tab:permutation-module-validation}
\small
\begin{tabular}{ccccccc}
\toprule
$a$ & $b$ & $d$ & $\dim \mathbb F_2[X_{a,b,d}]$ & Generators & Computed nullity & Orbit count \\
\midrule
$3$ & $3$ & $2$ & $15$ & $4$ & $3$ & $3$ \\
$4$ & $3$ & $3$ & $35$ & $5$ & $4$ & $4$ \\
$4$ & $4$ & $4$ & $70$ & $6$ & $5$ & $5$ \\
$5$ & $4$ & $3$ & $84$ & $7$ & $4$ & $4$ \\
$5$ & $5$ & $5$ & $252$ & $8$ & $6$ & $6$ \\
$6$ & $4$ & $4$ & $210$ & $8$ & $5$ & $5$ \\
\bottomrule
\end{tabular}
\end{table}

The finite-group example is intentionally elementary. Its role is to show that the packed fixed-space engine is not tied to Steenrod operations. The same layer can be used whenever a finite action on a basis is supplied by generators, whereas the Steenrod-specific layer supplies the extra quotient coordinates needed in Section~\ref{s2}.

\section{Proof and orbit-type post-processing}\label{s2}

The goal of this section is to justify the two structural statements made in the Introduction. The first layer of the proof is the actual hit-problem computation: admissible coordinates are constructed weight by weight, the Kameko matrix is formed, and the equations for the generators $\rho_1,\ldots,\rho_6$ are solved over $\mathbb F_2$. The second layer is applied only after the invariant classes have been obtained; it groups the supports of the reduced representatives into $\Sigma_6$-orbit types and computes the stabilizers recorded in the tables.

\begin{rem}\label{rem:why-degree15}
The calculation in degree $15$ is not an auxiliary check. Since $(\widetilde{Sq}^0_*)_{(6,36)}:(QP_6)_{36}\to(QP_6)_{15}$ is a surjective $GL(6)$-homomorphism, the image of every $GL(6)$-invariant class in degree $36$ is a $GL(6)$-invariant class in degree $15$. Proposition \ref{md15} says that the target invariant space is $\mathbb F_2\cdot[\xi]$. Therefore every $[h]\in[(QP_6)_{36}]^{GL(6)}$ must have the form
\[
[h]=\beta[\psi(\xi)]+[h^*],
\qquad \beta\in\mathbb F_2,
\qquad [h^*]\in\ker(\widetilde{Sq}^0_*)_{(6,36)},
\]
where $\psi(x_1^{e_1}\cdots x_6^{e_6})=x_1^{2e_1+1}\cdots x_6^{2e_6+1}$. Thus the degree-$15$ computation supplies the only possible non-kernel contribution to a degree-$36$ invariant. To prove Theorem \ref{dlc2}, it is not enough to know the dimension of the Kameko kernel invariants; one must also test whether the lifted target invariant $\psi(\xi)$ can be corrected by a kernel element. The final $\rho_6$-equation shows that it cannot, forcing $\beta=0$.
\end{rem}

For degree $15$, the algorithm gives
\[
(QP_6)_{15}\cong\bigoplus_{i=1}^7QP_6(\omega_{(i)}),
\]
where
\[
\omega_{(1)}=(1,1,1,1),\quad \omega_{(2)}=(1,1,3),\quad \omega_{(3)}=(1,3,2),\quad \omega_{(4)}=(3,2,2),
\]
\[
\omega_{(5)}=(3,4,1),\quad \omega_{(6)}=(5,3,1),\quad \omega_{(7)}=(5,5).
\]
The dimensions of the seven summands are displayed in the following array. In the row-reduced admissible coordinate system, the only non-zero $GL(6)$-invariant weight component occurs at $\omega_{(3)}$, and it is represented before lower-weight correction by $x_1x_2^2x_3^2x_4^2x_5^4x_6^4$. Solving the equations for $\rho_1,\ldots,\rho_6$ determines the lower-weight correction uniquely and gives the element $\xi$ in Proposition \ref{md15}.

\begin{center}
\small
{\bf Weight decomposition of $(QP_6)_{15}$.}\par\medskip
\begin{tabular}{ccc}
\toprule
Weight & Dimension of $QP_6(\omega)$ & Dimension of $[QP_6(\omega)]^{GL(6)}$ \\
\midrule
$(1,1,1,1)$ & 56 & 0 \\
$(1,1,3)$ & 6 & 0 \\
$(1,3,2)$ & 1 & 1 \\
$(3,2,2)$ & 1176 & 0 \\
$(3,4,1)$ & 384 & 0 \\
$(5,3,1)$ & 540 & 0 \\
$(5,5)$ & 21 & 0 \\
\bottomrule
\end{tabular}
\end{center}

For degree $36$, the admissible monomials are first reduced by the hit matrix. Since $\dim(P_6)_{36}=\binom{41}{5}=749398$, direct manual enumeration is not reliable. The Kameko matrix is then built by the exponent rule $\mathbf a\mapsto(\mathbf a-\mathbf1)/2$ on exponent vectors with all entries odd, followed by reduction to admissible target coordinates in degree $15$. The kernel decomposes into five weight components
\[
\ker(\widetilde{Sq}^0_*)_{(6,36)}\cong\bigoplus_{i=1}^5QP_6(\omega^*_{(i)}),
\]
where
\[
\omega^*_{(1)}=(4,2,1,1,1),\quad \omega^*_{(2)}=(4,2,1,3),\quad \omega^*_{(3)}=(4,2,3,2),
\]
\[
\omega^*_{(4)}=(4,4,2,2),\quad \omega^*_{(5)}=(4,4,4,1).
\]
The dimensions of these components and of their $\Sigma_6$-invariant parts are displayed in the following array. The sum of the five kernel dimensions is $12390$, and the sum of the five associated-graded $\Sigma_6$-candidate dimensions is $52$. The subsequent lower-weight correction systems, recorded explicitly in the output, reduce these candidates to the two-dimensional $GL(6)$-invariant subspace represented by $\zeta_1$ and $\zeta_2$.

\begin{center}
\small
{\bf Weight decomposition of the Kameko kernel in degree $36$.}\par\medskip
\begin{tabular}{ccc}
\toprule
Weight & Kernel dimension & Dimension of the $\Sigma_6$-invariant part \\
\midrule
$(4,2,1,1,1)$ & 2725 & 13 \\
$(4,2,1,3)$ & 111 & 2 \\
$(4,2,3,2)$ & 1085 & 6 \\
$(4,4,2,2)$ & 6495 & 18 \\
$(4,4,4,1)$ & 1974 & 13 \\
\bottomrule
\end{tabular}
\end{center}

The computational procedure can be summarized as follows.

\begin{center}
\refstepcounter{algorithm}\label{alg:main}
\textbf{Algorithm~\thealgorithm. Weight-vector, Kameko-kernel, and orbit-type computation}
\end{center}
{\footnotesize
\begin{algorithmic}[1]
\Require Integers $q\geq1$ and $n\geq0$ with $n\equiv q\pmod 2$; the Steenrod action on $P_q=\mathbb F_2[x_1,\ldots,x_q]$.
\Ensure A verified basis of $[(QP_q)_n]^{GL(q)}$ and, for every accepted reduced representative, its $\Sigma_q$-orbit-type table.

\Statex \textbf{Phase 0: exponent arithmetic and Steenrod columns.}

\Function{Weight}{$\mathbf a=(a_1,\ldots,a_q)$}
  \State $m\gets\max_i a_i$; if $m=0$, return the empty vector.
  \State Let $L$ be the least integer such that $a_i<2^L$ for all $i$.
  \For{$b=0,\ldots,L-1$}
    \State $\omega_{b+1}\gets\sum_{i=1}^q\big((a_i\div 2^b)\bmod 2\big)$.
  \EndFor
  \State \Return $\omega(\mathbf a)=(\omega_1,\ldots,\omega_L)$.
\EndFunction

\Function{KamekoImageExps}{$\mathbf a=(a_1,\ldots,a_q)$}
  \If{some $a_i$ is even}
    \State \Return \textsc{None}.
  \Else
    \State \Return $((a_1-1)/2,\ldots,(a_q-1)/2)$.
  \EndIf
\EndFunction

\Function{SqOnMono}{$s,\,x_1^{e_1}\cdots x_q^{e_q}$}
  \State Apply the Cartan formula
  \[
  Sq^s(x_1^{e_1}\cdots x_q^{e_q})=
  \sum_{s_1+\cdots+s_q=s}\prod_{i=1}^q\binom{e_i}{s_i}x_i^{e_i+s_i}
  \]
  and retain precisely the summands for which every binomial coefficient is odd, using Lucas' criterion.
  \State Cancel repeated monomials modulo $2$ and return the sorted exponent support.
\EndFunction

\Function{HitColumnExps}{$\mathbf b,\,s$}
  \State Compute \Call{SqOnMono}{$s,\,x^{\mathbf b}$}; collect the exponent tuples occurring with coefficient $1$.
  \State \Return the sorted list of exponent tuples.
\EndFunction

\Function{ExpsEnum}{$q,N$}
  \State \Return the stars-and-bars enumeration of all $\mathbf a\in\mathbb N^q$ satisfying $a_1+\cdots+a_q=N$.
\EndFunction

\Statex \textbf{Phase 1: streaming hit elimination and admissible coordinates.}

\Function{BuildDegSpaceOnline}{$q,N$}
  \State $\mathcal E\gets$ \Call{ExpsEnum}{$q,N$}; sort $\mathcal E$ by the ordered pair $(\omega(\mathbf a),\mathbf a)$.
  \State Construct the dictionary $\texttt{idx}:\mathcal E\to\{1,\ldots,|\mathcal E|\}$.
  \State Initialize an empty sparse pivot map $\texttt{pivotmap}$ for ONLINE elimination over $\mathbb F_2$.
  \For{$p=0,1,2,\ldots$ with $2^p\leq N$}
    \State $s\gets2^p$ and $N_g\gets N-s$.
    \ForAll{$\mathbf b\in$ \Call{ExpsEnum}{$q,N_g$}}
      \State $C\gets$ \Call{HitColumnExps}{$\mathbf b,s$}.
      \State Replace $C$ by the sorted row-index list $\{\texttt{idx}(\mathbf c):\mathbf c\in C\}$.
      \State Reduce this row list by repeated XOR against $\texttt{pivotmap}$.
      \If{a new leading row remains after reduction}
        \State Store the reduced column in $\texttt{pivotmap}$ with this row as pivot.
      \EndIf
    \EndFor
  \EndFor
  \State Let $S_{\rm piv}$ be the set of pivot rows in $\texttt{pivotmap}$ and set $A=\{1,\ldots,|\mathcal E|\}\setminus S_{\rm piv}$.
  \State \Return $\mathsf{DS}(q,N)$, consisting of $\mathcal E$, $\texttt{idx}$, the admissible exponents indexed by $A$, and the reduction data $\texttt{pivotmap}$.
\EndFunction

\Function{ReduceRowToAdmissible}{$r,\mathsf{DS}$}
  \State XOR-reduce the singleton row $[r]$ by $\texttt{pivotmap}$ until no pivot row remains.
  \State Translate the surviving row indices into admissible-coordinate positions, with parity taken modulo $2$.
  \State \Return the sorted admissible position list.
\EndFunction

\Statex \textbf{Phase 2: Kameko kernel, bit-packed nullspaces, and $GL(q)$-invariance.}

\Function{BuildKamekoBitMat}{$\mathsf{DS}_{\rm src},\mathsf{DS}_{\rm tgt}$}
  \State Create a bit-packed matrix $L$ of size $\dim(\mathsf{DS}_{\rm tgt}^{\rm adm})\times\dim(\mathsf{DS}_{\rm src}^{\rm adm})$.
  \For{each source admissible exponent $\mathbf a$ with column $c$}
    \State $\mathbf u\gets$ \Call{KamekoImageExps}{$\mathbf a$}.
    \If{$\mathbf u=\textsc{None}$}
      \State Continue to the next column.
    \EndIf
    \State Let $r=\texttt{idx}_{\rm tgt}(\mathbf u)$; reduce $r$ by \Call{ReduceRowToAdmissible}{$r,\mathsf{DS}_{\rm tgt}$}.
    \State Set the corresponding bits in column $c$ of $L$ equal to $1$.
  \EndFor
  \State \Return $L$.
\EndFunction

\Function{NullspaceGFTwo}{bit-packed matrix $M$}
  \State Perform Gaussian elimination over $\mathbb F_2$ using word-level XOR row operations.
  \State \Return the rank of $M$ and a basis for $\ker M$ in bit-vector form.
\EndFunction

\Function{ApplyRho}{$j,\,x_1^{e_1}\cdots x_q^{e_q}$}
  \If{$1\leq j<q$}
    \State Apply the adjacent transposition interchanging $x_j$ and $x_{j+1}$.
  \ElsIf{$j=q$}
    \State Substitute $x_q\mapsto x_q+x_{q-1}$ and expand over $\mathbb F_2$.
  \Else
    \State Return the input monomial.
  \EndIf
\EndFunction

\Function{DecomposeToEntries}{$f,\mathsf{DS}$}
  \State Write $f$ as an $\mathbb F_2$-sum of monomials.
  \State Map each monomial to its global row index, reduce by \Call{ReduceRowToAdmissible}{--,$\mathsf{DS}$}, and cancel repeated admissible positions modulo $2$.
  \State \Return the resulting sorted admissible-coordinate list.
\EndFunction

\Function{PrecomputeRhoRows}{$\mathsf{DS}$}
  \For{$j=1,\ldots,q$}
    \For{each admissible monomial $u_i$ in $\mathsf{DS}$}
      \State Store \Call{DecomposeToEntries}{$\rho_j(u_i)+u_i,\mathsf{DS}$} as the row support of $(\rho_j-\mathrm{Id})u_i$.
    \EndFor
  \EndFor
\EndFunction

\Function{SigmaGLOnKernelWeight}{$\mathsf{DS},\ker L,\mathcal I_\omega$}
  \State Let $\mathcal I_\omega$ be the admissible coordinates of weight $\omega$ occurring in at least one vector of $\ker L$.
  \State Build the stacked bit-matrix for
  \[
  (\rho_j-\mathrm{Id})\sum_{i\in\mathcal I_\omega}\gamma_i u_i\equiv0,
  \qquad 1\leq j\leq q-1,
  \]
  using the precomputed row supports.
  \State Compute its nullspace; this gives the weightwise $\Sigma_q$-invariant candidates.
  \State On that candidate space, build the additional bit-matrix for $(\rho_q-\mathrm{Id})$.
  \State Compute the nullspace again to obtain the weightwise $GL(q)$-invariant candidates.
  \State \Return the $\Sigma_q$-basis and the $GL(q)$-basis in weight $\omega$.
\EndFunction

\Function{CorrectByLowerWeights}{$g_{\max},\mathsf{DS},\ker L,\omega^\star$}
  \State Let $\mathcal L$ be the kernel coordinates of weight strictly smaller than $\omega^\star$.
  \State Solve the first stacked system
  \[
  (\rho_j-\mathrm{Id})\left(g_{\max}+\sum_{t\in\mathcal L}\beta_tu_t\right)\equiv0,
  \qquad 1\leq j\leq q-1.
  \]
  \State On the solution space, solve the second system imposed by $(\rho_q-\mathrm{Id})$.
  \State Verify the final representative against all equations $(\rho_j-\mathrm{Id})f\equiv0$ for $1\leq j\leq q$.
  \State \Return all representatives passing this verification.
\EndFunction

\Function{CorrectLiftFromTarget}{$g,\mathsf{DS}_{\rm src},\ker L$}
  \State Form the inverse Kameko lift
  \[
  \psi(g)=\sum x_1^{2e_1+1}\cdots x_q^{2e_q+1},
  \]
  where the sum runs through the monomials $x_1^{e_1}\cdots x_q^{e_q}$ in $g$.
  \State Add an arbitrary kernel correction and solve the same two-stage $\rho_1,\ldots,\rho_q$ system as in \Call{CorrectByLowerWeights}{--,--,--,--}.
  \State \Return the accepted corrected lifts.
\EndFunction

\Statex \textbf{Phase 3: deterministic orbit-type post-processing.}

\Function{PartitionType}{$\mathbf a=(a_1,\ldots,a_q)$}
  \State \Return the non-increasing rearrangement $\lambda=(\lambda_1\geq\cdots\geq\lambda_q)$ of the entries of $\mathbf a$.
\EndFunction

\Function{StabilizerData}{$\lambda$}
  \State Compute the multiplicities $\mu_e(\lambda)=|\{i:\lambda_i=e\}|$.
  \State $s_\lambda\gets\prod_e\mu_e(\lambda)!$ and $o_\lambda\gets q!/s_\lambda$.
  \State $T_\lambda\gets\prod_{\mu_e(\lambda)>1}\Sigma_{\mu_e(\lambda)}$, with $T_\lambda=\{1\}$ if all multiplicities are $1$.
  \State \Return $(s_\lambda,o_\lambda,T_\lambda)$.
\EndFunction

\Function{OrbitTypePostProcess}{$f,q$}
  \State Parse the reduced representative $f=\sum_{\mathbf a\in E(f)}x^{\mathbf a}$ and check that there are no repeated monomials.
  \State Initialize a dictionary $\mathcal D$ whose keys are partition types $\lambda$.
  \ForAll{$\mathbf a\in E(f)$}
    \State $\lambda\gets$ \Call{PartitionType}{$\mathbf a$}; append $\mathbf a$ to $\mathcal D[\lambda]$.
  \EndFor
  \ForAll{partition types $\lambda$ in the deterministic output order}
    \State Let $\mathbf a_{\rm seed}$ be the first exponent vector of type $\lambda$ in the raw output.
    \State $\mathrm{slice}\gets|\mathcal D[\lambda]|$ and $(s_\lambda,o_\lambda,T_\lambda)\gets$ \Call{StabilizerData}{$\lambda$}.
    \State Output the row $(\lambda,\mathbf a_{\rm seed},\mathrm{slice},o_\lambda,s_\lambda,T_\lambda)$.
  \EndFor
\EndFunction

\Statex \textbf{Main orchestration.}

\Procedure{RunAll}{$q,n$}
  \State Require $n\geq q$ and $n\equiv q\pmod2$; set $n_{\rm tgt}=(n-q)/2$.
  \State $\mathsf{DS}_{\rm src}\gets$ \Call{BuildDegSpaceOnline}{$q,n$}.
  \State $\mathsf{DS}_{\rm tgt}\gets$ \Call{BuildDegSpaceOnline}{$q,n_{\rm tgt}$}.
  \State $L\gets$ \Call{BuildKamekoBitMat}{$\mathsf{DS}_{\rm src},\mathsf{DS}_{\rm tgt}$}; compute $(\operatorname{rank}L,\ker L)$ by \Call{NullspaceGFTwo}{$L$}.
  \State \Call{PrecomputeRhoRows}{$\mathsf{DS}_{\rm src}$}.
  \State Extract the kernel support coordinates and group them by weight vector.
  \For{each occurring weight $\omega$}
    \State $(\Sigma_q[\omega],GL(q)[\omega])\gets$ \Call{SigmaGLOnKernelWeight}{$\mathsf{DS}_{\rm src},\ker L,\mathcal I_\omega$}.
  \EndFor
  \State Choose each largest-weight candidate $g_{\max}$ and apply \Call{CorrectByLowerWeights}{$g_{\max},\mathsf{DS}_{\rm src},\ker L,\omega^\star$}.
  \State For every known target invariant $g\in[(QP_q)_{n_{\rm tgt}}]^{GL(q)}$, apply \Call{CorrectLiftFromTarget}{$g,\mathsf{DS}_{\rm src},\ker L$}.
  \State Let $\mathcal B$ be the union of all representatives passing the full $\rho_1,\ldots,\rho_q$ verification.
  \ForAll{$f\in\mathcal B$}
    \State Produce the orbit-type table by \Call{OrbitTypePostProcess}{$f,q$}.
  \EndFor
  \State \Return $\mathcal B$ together with all weight dimensions, Kameko-kernel dimensions, verification logs, and orbit-type rows.
\EndProcedure
\end{algorithmic}
}

\begin{rem}\label{rem:algorithm-structural-bridge}
Algorithm~\ref{alg:main} separates the proof into an exact linear-algebra stage and a structural post-processing stage. The first stage is necessary because the raw spaces are already too large for reliable manual control: in the present bidegree one has
\[
\dim_{\mathbb F_2}(P_6)_{36}=\binom{36+6-1}{6-1}=749{,}398.
\]
The detailed output for degree $36$ records the streaming construction of $1{,}960{,}002$ Steenrod hit columns, the ONLINE reduction to $734{,}824$ pivot columns, and the resulting admissible quotient dimension $14{,}574$. These numbers explain why the computation is organized around streamed columns, sparse pivot maps, and bit-packed Gaussian elimination over $\mathbb F_2$: the program never needs to store a full dense hit matrix, and each Kameko, symmetric-group, or general-linear invariance test is reduced to an explicit finite nullspace computation.

The orbit-type calculation does not replace the Kameko-weight computation and is not used as a shortcut to prove invariance. It begins only after the representatives have passed the equations $(\rho_j-\mathrm{Id})f\equiv0$ in $QP_q$ for all generators $\rho_j$ of $GL(q)$. At that point the raw polynomial support is parsed deterministically, sorted by partition type, and converted into the slice, full-orbit, and stabilizer data displayed in Theorem~\ref{dlc2} and Proposition~\ref{md15}. Thus the same exact computation that produces the counterexample also produces a readable structural summary: the long representatives are transformed into orbit-type data without obscuring the linear algebra that verifies the theorem.
\end{rem}

\subsection{Structural decomposition via orbit types}\label{subsec:orbit-post-processing}

We spell out how the entries in Tables \ref{tab:xi-composition-data}, \ref{tab:zeta1-orbit-data}, and \ref{tab:zeta2-orbit-data} are obtained from the raw representatives produced by the weight-vector and Kameko-kernel computation. For a monomial
\[
 m=x_1^{a_1}\cdots x_6^{a_6},
\]
write
\[
\mathbf a(m)=(a_1,\ldots,a_6),\qquad
\lambda(m)=(\lambda_1\geq\cdots\geq\lambda_6)
\]
for the non-increasing rearrangement of the exponent vector. If \(\mu_e(\lambda)\) denotes the multiplicity of the integer \(e\) among the entries of \(\lambda\), then the variables carrying equal exponents may be permuted freely and no other variable permutation fixes the monomial. Hence
\[
\Stab_{\Sigma_6}(m)\cong\prod_e\Sigma_{\mu_e(\lambda)},\qquad
|\Stab_{\Sigma_6}(m)|=\prod_e\mu_e(\lambda)!,\qquad
|\mathcal O_{\Sigma_6}(m)|=\frac{6!}{\prod_e\mu_e(\lambda)!}.
\]
For a polynomial representative \(f\), the post-processing program parses \(\operatorname{Supp}(f)\), groups its exponent vectors according to the common value of \(\lambda\), and records the corresponding slice size
\[
 |E_\lambda(f)|=\bigl|\{\mathbf a\in\operatorname{Supp}(f):\lambda(\mathbf a)=\lambda\}\bigr|.
\]
The slice size is therefore a support count in the displayed reduced representative, whereas the column headed ``Full orbit'' is the size of the entire \(\Sigma_6\)-orbit in the ordinary monomial basis of \(P_6\). These two numbers need not coincide for \(\zeta_1\) and \(\zeta_2\), because these are reduced representatives of quotient classes in \(QP_6\), not literal sums of all monomials in each full \(\Sigma_6\)-orbit.

As a first explicit calculation, take the degree-$15$ monomial
\[
 m=x_1x_2^2x_3^4x_4^8
\]
which occurs in the block \(M_{(1,2,4,8)}^{(6)}\) of \(\xi\). Its exponent vector is
\[
\mathbf a(m)=(1,2,4,8,0,0),
\]
and the associated partition-type vector is
\[
\lambda(m)=(8,4,2,1,0,0).
\]
The multiplicities are \(\mu_0=2\), \(\mu_1=\mu_2=\mu_4=\mu_8=1\), and all other multiplicities are zero. Therefore
\[
|\Stab_{\Sigma_6}(m)|=2!\cdot1!\cdot1!\cdot1!\cdot1!=2,
\qquad
\Stab_{\Sigma_6}(m)\cong\Sigma_2,
\]
and the full orbit has cardinality
\[
|\mathcal O_{\Sigma_6}(m)|=\frac{720}{2}=360.
\]
On the other hand, the quasisymmetric block \(M_{(1,2,4,8)}^{(6)}\) keeps the exponent order \((1,2,4,8)\) and chooses only the increasing index set \(i_1<i_2<i_3<i_4\). Hence its support contains \(\binom64=15\) monomials, which is the slice entry in Table \ref{tab:xi-composition-data}.

As a second explicit calculation, take the degree-$36$ monomial
\[
 m=x_1x_2x_3x_4x_5^{30}x_6^2
\]
from the raw representative \(\zeta_1\). Here
\[
\mathbf a(m)=(1,1,1,1,30,2),
\qquad
\lambda(m)=(30,2,1,1,1,1).
\]
The repeated exponent is \(1\), with multiplicity \(4\); the exponents \(30\) and \(2\) each have multiplicity \(1\). Thus
\[
|\Stab_{\Sigma_6}(m)|=4!\cdot1!\cdot1!=24,
\qquad
\Stab_{\Sigma_6}(m)\cong\Sigma_4,
\]
and the full orbit has cardinality
\[
|\mathcal O_{\Sigma_6}(m)|=\frac{720}{24}=30.
\]
The raw support of \(\zeta_1\) contains exactly two monomials with this same partition type, so the corresponding slice entry is \(2\), even though the full \(\Sigma_6\)-orbit contains \(30\) monomials.

The table generation is consequently reproducible without any subjective choice. The program reads the raw polynomial, extracts every exponent vector, checks that all total degrees are correct and that no monomial is duplicated, sorts each vector to obtain \(\lambda\), groups identical \(\lambda\)-types, computes the stabilizer order and type from the multiplicities, and prints as seed the first exponent vector of that type in the raw representative. Applying this deterministic post-processing gives \(63\) monomials and six partition types for \(\xi\), \(39\) monomials and ten partition types for \(\zeta_1\), and \(539\) monomials and ninety-six partition types for \(\zeta_2\). Summing the slice entries in the three tables gives exactly \(63\), \(39\), and \(539\), respectively.

\begin{longtable}{cccccc}
\caption{Orbit-type data for $\zeta_1$.}\label{tab:zeta1-orbit-data}\\
\toprule
$\lambda$ & Seed $\mathbf a$ & Slice & Full orbit & $|\operatorname{Stab}|$ & Stabilizer type \\
\midrule
\endfirsthead
\toprule
$\lambda$ & Seed $\mathbf a$ & Slice & Full orbit & $|\operatorname{Stab}|$ & Stabilizer type \\
\midrule
\endhead
$(30,2,1,1,1,1)$ & $(1,1,1,1,30,2)$ & 2 & 30 & 24 & $\Sigma_{4}$ \\
$(28,3,2,1,1,1)$ & $(1,1,3,1,28,2)$ & 6 & 120 & 6 & $\Sigma_{3}$ \\
$(26,6,1,1,1,1)$ & $(1,1,1,1,6,26)$ & 4 & 30 & 24 & $\Sigma_{4}$ \\
$(26,4,3,1,1,1)$ & $(3,1,1,1,4,26)$ & 4 & 120 & 6 & $\Sigma_{3}$ \\
$(24,6,3,1,1,1)$ & $(1,3,6,1,1,24)$ & 2 & 120 & 6 & $\Sigma_{3}$ \\
$(24,5,3,2,1,1)$ & $(1,3,5,1,24,2)$ & 9 & 360 & 2 & $\Sigma_{2}$ \\
$(24,4,3,3,1,1)$ & $(1,3,3,1,4,24)$ & 7 & 180 & 4 & $\Sigma_{2} \times \Sigma_{2}$ \\
$(16,10,5,3,1,1)$ & $(3,5,10,1,1,16)$ & 1 & 360 & 2 & $\Sigma_{2}$ \\
$(16,9,5,3,2,1)$ & $(3,5,9,1,16,2)$ & 2 & 720 & 1 & ${1}$ \\
$(16,8,5,3,3,1)$ & $(3,5,3,1,8,16)$ & 2 & 360 & 2 & $\Sigma_{2}$ \\
\bottomrule
\end{longtable}

\begin{longtable}{cccccc}
\caption{Orbit-type data for $\zeta_2$.}\label{tab:zeta2-orbit-data}\\
\toprule
$\lambda$ & Seed $\mathbf a$ & Slice & Full orbit & $|\operatorname{Stab}|$ & Stabilizer type \\
\midrule
\endfirsthead
\toprule
$\lambda$ & Seed $\mathbf a$ & Slice & Full orbit & $|\operatorname{Stab}|$ & Stabilizer type \\
\midrule
\endhead
$(30,2,1,1,1,1)$ & $(1,1,1,1,2,30)$ & 1 & 30 & 24 & $\Sigma_{4}$ \\
$(28,3,2,1,1,1)$ & $(1,1,1,2,3,28)$ & 2 & 120 & 6 & $\Sigma_{3}$ \\
$(27,4,2,1,1,1)$ & $(1,1,1,2,4,27)$ & 1 & 120 & 6 & $\Sigma_{3}$ \\
$(26,6,1,1,1,1)$ & $(1,1,1,1,6,26)$ & 1 & 30 & 24 & $\Sigma_{4}$ \\
$(26,4,3,1,1,1)$ & $(1,1,1,3,4,26)$ & 1 & 120 & 6 & $\Sigma_{3}$ \\
$(25,6,2,1,1,1)$ & $(1,1,1,2,6,25)$ & 1 & 120 & 6 & $\Sigma_{3}$ \\
$(25,4,3,2,1,1)$ & $(1,1,2,4,3,25)$ & 2 & 360 & 2 & $\Sigma_{2}$ \\
$(24,7,2,1,1,1)$ & $(1,1,1,2,7,24)$ & 6 & 120 & 6 & $\Sigma_{3}$ \\
$(24,6,3,1,1,1)$ & $(1,1,1,3,6,24)$ & 5 & 120 & 6 & $\Sigma_{3}$ \\
$(24,5,3,2,1,1)$ & $(1,1,2,3,5,24)$ & 9 & 360 & 2 & $\Sigma_{2}$ \\
$(24,4,3,3,1,1)$ & $(1,3,1,3,4,24)$ & 5 & 180 & 4 & $\Sigma_{2} \times \Sigma_{2}$ \\
$(19,8,4,3,1,1)$ & $(1,3,4,1,8,19)$ & 2 & 360 & 2 & $\Sigma_{2}$ \\
$(18,9,5,2,1,1)$ & $(1,1,2,5,9,18)$ & 1 & 360 & 2 & $\Sigma_{2}$ \\
$(18,9,4,3,1,1)$ & $(1,3,1,4,9,18)$ & 4 & 360 & 2 & $\Sigma_{2}$ \\
$(18,8,5,3,1,1)$ & $(1,1,3,5,8,18)$ & 4 & 360 & 2 & $\Sigma_{2}$ \\
$(17,10,5,2,1,1)$ & $(1,1,2,5,10,17)$ & 1 & 360 & 2 & $\Sigma_{2}$ \\
$(17,10,4,3,1,1)$ & $(1,3,4,1,10,17)$ & 2 & 360 & 2 & $\Sigma_{2}$ \\
$(17,9,6,2,1,1)$ & $(1,1,2,6,9,17)$ & 3 & 360 & 2 & $\Sigma_{2}$ \\
$(17,9,4,3,2,1)$ & $(3,1,4,2,9,17)$ & 2 & 720 & 1 & ${1}$ \\
$(17,8,7,2,1,1)$ & $(1,1,2,7,8,17)$ & 2 & 360 & 2 & $\Sigma_{2}$ \\
$(17,8,6,3,1,1)$ & $(1,1,3,6,8,17)$ & 3 & 360 & 2 & $\Sigma_{2}$ \\
$(17,8,5,3,2,1)$ & $(1,3,5,2,8,17)$ & 3 & 720 & 1 & ${1}$ \\
$(17,8,4,3,3,1)$ & $(1,3,4,3,8,17)$ & 5 & 360 & 2 & $\Sigma_{2}$ \\
$(16,11,4,3,1,1)$ & $(1,3,4,1,11,16)$ & 2 & 360 & 2 & $\Sigma_{2}$ \\
$(16,9,7,2,1,1)$ & $(1,1,2,7,9,16)$ & 2 & 360 & 2 & $\Sigma_{2}$ \\
$(16,9,6,3,1,1)$ & $(1,1,6,3,9,16)$ & 7 & 360 & 2 & $\Sigma_{2}$ \\
$(16,9,5,3,2,1)$ & $(1,3,5,2,9,16)$ & 8 & 720 & 1 & ${1}$ \\
$(16,9,4,3,3,1)$ & $(3,1,4,9,3,16)$ & 3 & 360 & 2 & $\Sigma_{2}$ \\
$(16,8,7,3,1,1)$ & $(1,3,7,1,8,16)$ & 8 & 360 & 2 & $\Sigma_{2}$ \\
$(16,8,5,3,3,1)$ & $(1,3,5,8,3,16)$ & 6 & 360 & 2 & $\Sigma_{2}$ \\
$(14,14,3,3,1,1)$ & $(1,1,3,3,14,14)$ & 6 & 90 & 8 & $\Sigma_{2} \times \Sigma_{2} \times \Sigma_{2}$ \\
$(14,12,3,3,3,1)$ & $(1,3,3,14,3,12)$ & 3 & 120 & 6 & $\Sigma_{3}$ \\
$(14,11,6,3,1,1)$ & $(1,1,3,6,11,14)$ & 10 & 360 & 2 & $\Sigma_{2}$ \\
$(14,11,5,3,2,1)$ & $(3,5,1,2,11,14)$ & 4 & 720 & 1 & ${1}$ \\
$(14,11,4,3,3,1)$ & $(3,1,3,4,11,14)$ & 6 & 360 & 2 & $\Sigma_{2}$ \\
$(14,10,5,3,3,1)$ & $(1,14,3,5,3,10)$ & 7 & 360 & 2 & $\Sigma_{2}$ \\
$(14,9,6,3,3,1)$ & $(1,3,6,3,14,9)$ & 4 & 360 & 2 & $\Sigma_{2}$ \\
$(14,9,5,3,3,2)$ & $(3,5,2,3,14,9)$ & 4 & 360 & 2 & $\Sigma_{2}$ \\
$(14,8,5,3,3,3)$ & $(3,5,3,14,3,8)$ & 1 & 120 & 6 & $\Sigma_{3}$ \\
$(14,7,6,5,3,1)$ & $(1,3,5,6,14,7)$ & 12 & 720 & 1 & ${1}$ \\
$(14,7,5,4,3,3)$ & $(3,3,4,5,14,7)$ & 6 & 360 & 2 & $\Sigma_{2}$ \\
$(14,6,5,5,3,3)$ & $(3,3,5,14,5,6)$ & 10 & 180 & 4 & $\Sigma_{2} \times \Sigma_{2}$ \\
$(13,12,6,3,1,1)$ & $(1,1,3,6,12,13)$ & 10 & 360 & 2 & $\Sigma_{2}$ \\
$(13,12,5,3,2,1)$ & $(3,5,2,13,1,12)$ & 1 & 720 & 1 & ${1}$ \\
$(13,12,4,3,3,1)$ & $(3,1,3,13,4,12)$ & 14 & 360 & 2 & $\Sigma_{2}$ \\
$(13,10,6,5,1,1)$ & $(1,1,6,10,13,5)$ & 2 & 360 & 2 & $\Sigma_{2}$ \\
$(13,10,6,3,3,1)$ & $(1,3,3,13,6,10)$ & 9 & 360 & 2 & $\Sigma_{2}$ \\
$(13,10,5,3,3,2)$ & $(3,5,2,13,3,10)$ & 5 & 360 & 2 & $\Sigma_{2}$ \\
$(13,10,4,3,3,3)$ & $(3,3,3,13,4,10)$ & 2 & 120 & 6 & $\Sigma_{3}$ \\
$(13,8,6,5,3,1)$ & $(1,3,6,5,13,8)$ & 4 & 720 & 1 & ${1}$ \\
$(13,8,6,3,3,3)$ & $(3,3,3,13,6,8)$ & 1 & 120 & 6 & $\Sigma_{3}$ \\
$(13,8,5,5,3,2)$ & $(3,5,2,5,13,8)$ & 2 & 360 & 2 & $\Sigma_{2}$ \\
$(13,8,5,4,3,3)$ & $(3,3,13,4,5,8)$ & 8 & 360 & 2 & $\Sigma_{2}$ \\
$(13,6,6,5,3,3)$ & $(3,13,3,5,6,6)$ & 13 & 180 & 4 & $\Sigma_{2} \times \Sigma_{2}$ \\
$(12,12,7,3,1,1)$ & $(1,1,7,3,12,12)$ & 7 & 180 & 4 & $\Sigma_{2} \times \Sigma_{2}$ \\
$(12,12,5,3,3,1)$ & $(1,3,3,12,5,12)$ & 9 & 180 & 4 & $\Sigma_{2} \times \Sigma_{2}$ \\
$(12,12,3,3,3,3)$ & $(3,3,3,3,12,12)$ & 1 & 15 & 48 & $\Sigma_{4} \times \Sigma_{2}$ \\
$(12,11,6,5,1,1)$ & $(1,1,6,11,5,12)$ & 3 & 360 & 2 & $\Sigma_{2}$ \\
$(12,11,6,3,3,1)$ & $(1,3,6,11,3,12)$ & 6 & 360 & 2 & $\Sigma_{2}$ \\
$(12,11,5,4,3,1)$ & $(3,1,4,11,5,12)$ & 4 & 720 & 1 & ${1}$ \\
$(12,11,5,3,3,2)$ & $(3,5,2,11,3,12)$ & 5 & 360 & 2 & $\Sigma_{2}$ \\
$(12,11,4,3,3,3)$ & $(3,3,4,11,3,12)$ & 3 & 120 & 6 & $\Sigma_{3}$ \\
$(12,10,7,5,1,1)$ & $(1,1,7,10,5,12)$ & 7 & 360 & 2 & $\Sigma_{2}$ \\
$(12,10,7,3,3,1)$ & $(1,3,3,7,10,12)$ & 9 & 360 & 2 & $\Sigma_{2}$ \\
$(12,10,5,5,3,1)$ & $(1,3,5,10,5,12)$ & 3 & 360 & 2 & $\Sigma_{2}$ \\
$(12,10,5,3,3,3)$ & $(3,3,3,5,12,10)$ & 4 & 120 & 6 & $\Sigma_{3}$ \\
$(12,9,7,4,3,1)$ & $(3,1,7,9,4,12)$ & 1 & 720 & 1 & ${1}$ \\
$(12,9,6,5,3,1)$ & $(1,6,3,5,9,12)$ & 2 & 720 & 1 & ${1}$ \\
$(12,9,5,4,3,3)$ & $(3,3,4,5,9,12)$ & 6 & 360 & 2 & $\Sigma_{2}$ \\
$(12,8,7,5,3,1)$ & $(1,3,7,5,8,12)$ & 19 & 720 & 1 & ${1}$ \\
$(12,8,7,3,3,3)$ & $(3,3,7,3,12,8)$ & 9 & 120 & 6 & $\Sigma_{3}$ \\
$(12,8,5,5,3,3)$ & $(3,3,5,5,12,8)$ & 6 & 180 & 4 & $\Sigma_{2} \times \Sigma_{2}$ \\
$(12,7,6,5,3,3)$ & $(3,3,5,6,12,7)$ & 20 & 360 & 2 & $\Sigma_{2}$ \\
$(11,10,6,5,3,1)$ & $(1,3,5,11,6,10)$ & 17 & 720 & 1 & ${1}$ \\
$(11,10,5,5,3,2)$ & $(3,5,2,11,5,10)$ & 1 & 360 & 2 & $\Sigma_{2}$ \\
$(11,10,5,4,3,3)$ & $(3,3,4,11,5,10)$ & 9 & 360 & 2 & $\Sigma_{2}$ \\
$(11,9,8,4,3,1)$ & $(1,3,4,11,8,9)$ & 8 & 720 & 1 & ${1}$ \\
$(11,8,8,5,3,1)$ & $(1,3,5,8,11,8)$ & 5 & 360 & 2 & $\Sigma_{2}$ \\
$(11,8,6,5,3,3)$ & $(3,3,5,11,6,8)$ & 8 & 360 & 2 & $\Sigma_{2}$ \\
$(11,6,6,5,5,3)$ & $(3,5,6,11,5,6)$ & 1 & 180 & 4 & $\Sigma_{2} \times \Sigma_{2}$ \\
$(10,10,7,5,3,1)$ & $(1,3,7,5,10,10)$ & 3 & 360 & 2 & $\Sigma_{2}$ \\
$(10,10,5,5,3,3)$ & $(3,5,3,10,5,10)$ & 2 & 90 & 8 & $\Sigma_{2} \times \Sigma_{2} \times \Sigma_{2}$ \\
$(10,9,9,4,3,1)$ & $(1,3,4,10,9,9)$ & 4 & 360 & 2 & $\Sigma_{2}$ \\
$(10,9,8,5,3,1)$ & $(1,3,5,10,8,9)$ & 9 & 720 & 1 & ${1}$ \\
$(10,9,7,6,3,1)$ & $(1,3,6,7,9,10)$ & 9 & 720 & 1 & ${1}$ \\
$(10,9,7,4,3,3)$ & $(3,3,4,7,9,10)$ & 9 & 360 & 2 & $\Sigma_{2}$ \\
$(10,9,6,5,3,3)$ & $(3,3,5,10,6,9)$ & 4 & 360 & 2 & $\Sigma_{2}$ \\
$(10,8,7,5,3,3)$ & $(3,3,5,7,8,10)$ & 20 & 360 & 2 & $\Sigma_{2}$ \\
$(10,7,6,5,5,3)$ & $(3,5,7,10,5,6)$ & 9 & 360 & 2 & $\Sigma_{2}$ \\
$(9,9,8,6,3,1)$ & $(1,3,6,8,9,9)$ & 6 & 360 & 2 & $\Sigma_{2}$ \\
$(9,9,8,5,3,2)$ & $(3,5,2,8,9,9)$ & 3 & 360 & 2 & $\Sigma_{2}$ \\
$(9,9,8,4,3,3)$ & $(3,3,4,8,9,9)$ & 1 & 180 & 4 & $\Sigma_{2} \times \Sigma_{2}$ \\
$(9,8,8,7,3,1)$ & $(1,3,7,8,8,9)$ & 19 & 360 & 2 & $\Sigma_{2}$ \\
$(9,8,8,5,3,3)$ & $(3,3,5,8,8,9)$ & 4 & 180 & 4 & $\Sigma_{2} \times \Sigma_{2}$ \\
$(9,7,6,6,5,3)$ & $(3,5,6,7,6,9)$ & 16 & 360 & 2 & $\Sigma_{2}$ \\
$(8,8,7,5,5,3)$ & $(3,7,5,5,8,8)$ & 2 & 180 & 4 & $\Sigma_{2} \times \Sigma_{2}$ \\
\bottomrule
\end{longtable}

We now finish the proof of Theorem \ref{dlc2}. Let $g\in P_6$ represent a class in $[\ker(\widetilde{Sq}^0_*)_{(6,36)}]^{GL(6)}$. Writing $g$ in the admissible coordinates of the five kernel weights displayed above, the weightwise associated-graded calculation gives the five $\Sigma_6$-candidate dimensions
\[
13,\quad 2,\quad 6,\quad 18,\quad 13,
\]
whose sum is $52$. These equations are not yet the full ordinary $GL(6)$-invariance equations in $QP_6$, because a leading-weight representative may require corrections by lower kernel weights after the transvections are reduced modulo hit elements. The correction system used in the computation therefore fixes the largest kernel weight with non-zero associated-graded $GL(6)$-candidate, namely $(4,4,4,1)$, and allows all lower kernel coordinates. The equations for $\rho_1,\ldots,\rho_5$ give a system with $61950$ rows, $10417$ unknowns, and nullity $35$. Imposing the remaining transvection equation $\rho_6(g)+g\equiv0$ on these $35$ parameters gives a system with $12390$ rows, $35$ unknowns, and nullity $2$. The two resulting representatives are precisely $\zeta_1$ and $\zeta_2$, and both representatives were then verified against all six transvections in the quotient $QP_6$. Consequently
\[
[\ker(\widetilde{Sq}^0_*)_{(6,36)}]^{GL(6)}=\mathbb F_2\cdot[\zeta_1]\oplus\mathbb F_2\cdot[\zeta_2].
\]

It remains to pass from the Kameko kernel to the whole space $(QP_6)_{36}$. Let $[h]\in[(QP_6)_{36}]^{GL(6)}$. By $GL(6)$-equivariance and surjectivity of Kameko, its image in $(QP_6)_{15}$ is a scalar multiple of $[\xi]$. Hence
\[
h\equiv\beta\psi(\xi)+h^*,
\qquad \beta\in\mathbb F_2,
\qquad [h^*]\in\ker(\widetilde{Sq}^0_*)_{(6,36)}.
\]
The lift test is then solved with the unknowns consisting of $\beta$ and the kernel coordinates. In the first stage, using the equations for $\rho_1,\ldots,\rho_5$ and allowing lower-weight corrections gives a system with $61950$ rows, $12391$ unknowns, and nullity $46$. The final $\rho_6$-system has $12390$ rows, $46$ unknowns, and nullity $2$. The two accepted solutions are again exactly the two kernel representatives $\zeta_1$ and $\zeta_2$. Thus no additional invariant class with non-zero $\psi(\xi)$-coefficient occurs; equivalently, $\beta=0$ for every $GL(6)$-invariant solution. It follows that every $GL(6)$-invariant class in degree $36$ is accounted for within the Kameko kernel, yielding:
\[
[(QP_6)_{36}]^{GL(6)}=\mathbb F_2\cdot[\zeta_1]\oplus\mathbb F_2\cdot[\zeta_2].
\]
This result confirms Theorem \ref{dlc2}, from which Theorem \ref{dlc} and Corollary \ref{hq} follow directly by duality and the known calculation of $\operatorname{Ext}^{6,42}_{\mathscr A}(\mathbb F_2,\mathbb F_2)$.

\begin{cy}\label{cyf}
The raw computational outputs used in the proof are recorded in two public data records: the degree-$15$ target computation is archived in \cite{ZenodoDeg15}, and the degree-$36$ Kameko-kernel and invariant computation is archived in \cite{ZenodoDeg36}. The software archive listed in ~\ref{app:software} contains the new \texttt{Julia} package \texttt{AlgebraicTransfer.jl}, the direct \texttt{OSCAR} implementation of Algorithm~\ref{alg:main}, the Dickson and permutation-module validation drivers, and the scripts used to produce the orbit-type tables. The orbit-type stage parses the reduced representatives, groups exponent vectors by sorted partition type, computes $|E_\lambda(f)|$, and obtains the stabilizer order from $\prod_e\mu_e(\lambda)!$. This stage is a structural summary of monomial support; the proof of invariance is the quotient-level verification against the generators $\rho_1,\ldots,\rho_6$.
\end{cy}

\section{Unoriented bordism, geometric realization, and the transfer kernel}\label{s:geometry}

The preceding computation determines the transfer source exactly but does not, by itself, produce manifolds representing its homological dual. The purpose of this section is to pass from the algebraic source to unoriented bordism without identifying two constructions that are genuinely different. We first prove a general factorization of Singer's transfer through bordism classes over $BV_q$. We then apply it to the two degree-$36$ generators, test the Milnor-hypersurface, projective-product, and Dold-manifold models, and finally interpret the inverse Kameko monomial operation through the Thom space of the sum of the universal real line bundles.

Let $\gamma_i\to BV_q$ be the real line bundle associated with the $i$-th coordinate character and write $x_i=w_1(\gamma_i)$. We use the divided-power basis
\[
H_*(BV_q;\mathbb F_2)=\Gamma(v_1,\ldots,v_q),
\]
characterized by
\[
\left\langle x_1^{a_1}\cdots x_q^{a_q},
 v_1^{(b_1)}\cdots v_q^{(b_q)}\right\rangle
 =\prod_{i=1}^q\delta_{a_i,b_i}.
\]
For a space $X$, let $\mathfrak N_n(X)$ denote the unoriented bordism group of maps from closed smooth $n$-manifolds to $X$. The Thom homomorphism is
\[
\mu_{X,n}:\mathfrak N_n(X)\longrightarrow H_n(X;\mathbb F_2),
\qquad
\mu_{X,n}[M,f]=f_*[M].
\]

\begin{thm}[Geometric factorization of Singer's algebraic transfer]\label{thm:geometric-factorization}
For every $q\geq1$ and $n\geq0$, the Thom homomorphism
\[
\mu_{q,n}:=\mu_{BV_q,n}:\mathfrak N_n(BV_q)\longrightarrow H_n(BV_q;\mathbb F_2)
\]
is surjective. Define
\[
\mathfrak N_n^{\mathscr A}(BV_q)
=\mu_{q,n}^{-1}\bigl(\mathcal P_{\mathscr A}H_n(BV_q;\mathbb F_2)\bigr).
\]
A class $[M,f]$ belongs to $\mathfrak N_n^{\mathscr A}(BV_q)$ if and only if
\[
\left\langle v_r(M)f^*(a),[M]\right\rangle=0
\tag{\ref{thm:geometric-factorization}.1}\label{eq:mixed-wu-vanishing}
\]
for every $r>0$ and every $a\in H^{n-r}(BV_q;\mathbb F_2)$, where $v_r(M)$ is the $r$-th Wu class of $M$. The $GL(q)$-action on $BV_q$ induces a $GL(q)$-action on $\mathfrak N_n^{\mathscr A}(BV_q)$, and the Thom homomorphism induces a surjection
\[
\overline\mu_{q,n}:
\mathbb F_2\otimes_{GL(q)}\mathfrak N_n^{\mathscr A}(BV_q)
\longrightarrow
\mathbb F_2\otimes_{GL(q)}\mathcal P_{\mathscr A}H_n(BV_q;\mathbb F_2).
\]
Consequently Singer's transfer admits the geometric factorization
\[
Tr^{\mathrm{geom}}_{q,n}:=
Tr_q(\mathbb F_2)\circ\overline\mu_{q,n},
\]
and there is a natural isomorphism
\[
\frac{\ker Tr^{\mathrm{geom}}_{q,n}}{\ker\overline\mu_{q,n}}
\cong
\ker Tr_q(\mathbb F_2).
\tag{\ref{thm:geometric-factorization}.2}\label{eq:geometric-kernel-quotient}
\]
\end{thm}

\begin{proof}
Thom's mod-$2$ representability theorem states that every class in the mod-$2$ homology of a finite polyhedron is the image of the fundamental class of a closed smooth manifold under a continuous map \cite{Thom}. Every class of $H_n(BV_q;\mathbb F_2)$ is supported on a finite skeleton of $BV_q$. Applying Thom's theorem to that finite skeleton and composing with its inclusion into $BV_q$ proves the surjectivity of $\mu_{q,n}$.

Let $p=f_*[M]\in H_n(BV_q;\mathbb F_2)$. For $r>0$ and $a\in H^{n-r}(BV_q;\mathbb F_2)$, naturality of the Kronecker pairing and the defining property of the Wu class give
\[
\begin{aligned}
\left\langle a,Sq_*^r p\right\rangle
&=\left\langle Sq^r a,p\right\rangle
 =\left\langle f^*(Sq^r a),[M]\right\rangle\\
&=\left\langle Sq^r(f^*a),[M]\right\rangle
 =\left\langle v_r(M)f^*(a),[M]\right\rangle.
\end{aligned}
\tag{\ref{thm:geometric-factorization}.3}\label{eq:wu-adjunction}
\]
The Kronecker pairing is nondegenerate in each degree. Hence $Sq_*^r p=0$ for every $r>0$ if and only if all numbers in \eqref{eq:mixed-wu-vanishing} vanish. Since the positive Steenrod squares generate the augmentation ideal of $\mathscr A$, this is equivalent to $p\in\mathcal P_{\mathscr A}H_n(BV_q;\mathbb F_2)$.

An element $A\in GL(q)$ acts on $\mathfrak N_n(BV_q)$ by postcomposition with $BA:BV_q\to BV_q$. Naturality of Steenrod operations shows that $\mathfrak N_n^{\mathscr A}(BV_q)$ is stable under this action, and naturality of the Thom homomorphism makes its restriction $GL(q)$-equivariant. The restriction is surjective onto $\mathcal P_{\mathscr A}H_n(BV_q;\mathbb F_2)$ because every primitive homology class is, in particular, a mod-$2$ homology class and therefore has a manifold representative by the first paragraph. Passage to coinvariants gives the asserted surjection $\overline\mu_{q,n}$. Finally, for any surjective linear map $u:A\to B$ and any linear map $T:B\to C$, the map induced by $u$ identifies $\ker(Tu)/\ker u$ with $\ker T$. Taking $u=\overline\mu_{q,n}$ and $T=Tr_q(\mathbb F_2)$ proves \eqref{eq:geometric-kernel-quotient}.
\end{proof}

The Wu classes are universal polynomials in the classes $w_i(TM)$, determined by $w(TM)=Sq(v(M))$. Thus \eqref{eq:mixed-wu-vanishing} is an explicit family of mod-$2$ characteristic numbers involving the tangent bundle of $M$ and the six line bundles classified by $f$. These conditions depend essentially on the map $f$; the ordinary characteristic numbers of the underlying manifold cannot by themselves determine its image under the algebraic transfer.

\subsection{The two degree-$36$ source classes and their geometric representatives}

We now apply Theorem~\ref{thm:geometric-factorization} to the exact invariant calculation of Theorem~\ref{dlc2}. The resulting realization statement is existential but exact: the two dual classes are represented by closed smooth manifolds over $BV_6$, although the computation does not single out canonical manifolds or canonical primitive representatives in ordinary homology.

Set $E_{6,36}=(\mathbb F_2\otimes_{GL(6)}\mathcal P_{\mathscr A}H_*(BV_6;\mathbb F_2))_{36}.$
By finite-dimensional transpose duality, $E_{6,36}\cong ([(QP_6)_{36}]^{GL(6)})^{*}.$
Let $\{e_1, e_2\} \subset E_{6,36}$ be the basis dual to $\{[\zeta_1], [\zeta_2]\}$, where the classes $[\zeta_1]$ and $[\zeta_2]$ are determined as in Theorem \ref{dlc2}.

\begin{corl}[Exact bordism realization of the dual basis]\label{corl:zeta-bordism-realization}
There exist primitive classes
\[
p_i\in\mathcal P_{\mathscr A}H_{36}(BV_6;\mathbb F_2),
\qquad i=1,2,
\]
whose coinvariant classes are $e_i$, and there exist closed smooth $36$-manifolds $M_i$ with maps $f_i:M_i\to BV_6$ such that
\[
(f_i)_*[M_i]=p_i.
\]
If $L_{ij}=f_i^*(\gamma_j)$, then
\[
\left\langle
\zeta_j\bigl(w_1(L_{i1}),\ldots,w_1(L_{i6})\bigr),[M_i]
\right\rangle
=\delta_{ij}
\tag{\ref{corl:zeta-bordism-realization}.1}\label{eq:zeta-dual-numbers}
\]
for $1\leq i,j\leq2$, and
\[
\left\langle v_r(M_i)f_i^*(a),[M_i]\right\rangle=0
\tag{\ref{corl:zeta-bordism-realization}.2}\label{eq:zeta-wu-numbers}
\]
for every $r>0$ and every $a\in H^{36-r}(BV_6;\mathbb F_2)$.
\end{corl}

\begin{proof}
Choose a primitive representative $p_i$ of the coinvariant class $e_i$. Theorem~\ref{thm:geometric-factorization} supplies a closed smooth manifold and a map $f_i$ with $(f_i)_*[M_i]=p_i$. The equality \eqref{eq:zeta-dual-numbers} is the defining dual-basis relation, because evaluation of an invariant quotient class on a primitive homology class is unchanged when the latter is altered by a $GL(6)$-difference. The equality \eqref{eq:zeta-wu-numbers} follows from the primitivity of $p_i$ and \eqref{eq:wu-adjunction}.
\end{proof}

The phrase ``the homological dual of $[\zeta_i]$'' must therefore be understood at the coinvariant level. A class $e_i$ has many representatives in $\mathcal P_{\mathscr A}H_{36}(BV_6;\mathbb F_2)$, and each such representative has many bordism realizations. Corollary~\ref{corl:zeta-bordism-realization} proves existence with the exact evaluation conditions \eqref{eq:zeta-dual-numbers}; it does not identify a preferred diffeomorphism type.

\begin{prop}[The geometric transfer kernel in degree $36$]\label{prop:geometric-transfer-kernel}
There are unique coefficients $\epsilon_1,\epsilon_2\in\mathbb F_2$ such that
\[
Tr_6(\mathbb F_2)(e_i)=\epsilon_i t,
\qquad i=1,2.
\tag{\ref{prop:geometric-transfer-kernel}.1}\label{eq:epsilon-definition}
\]
For a class $[M,f]\in\mathfrak N_{36}^{\mathscr A}(BV_6)$, put
\[
a_i(M,f)=\left\langle f^*(\zeta_i),[M]\right\rangle.
\]
Then
\[
Tr^{\mathrm{geom}}_{6,36}\bigl(1\otimes[M,f]\bigr)
=\bigl(\epsilon_1a_1(M,f)+\epsilon_2a_2(M,f)\bigr)t.
\tag{\ref{prop:geometric-transfer-kernel}.2}\label{eq:geometric-transfer-functional}
\]
If $(\epsilon_1,\epsilon_2)\neq(0,0)$, then
\[
\ker Tr_6(\mathbb F_2)
=\mathbb F_2\cdot(\epsilon_2e_1+\epsilon_1e_2).
\tag{\ref{prop:geometric-transfer-kernel}.3}\label{eq:kernel-line-epsilon}
\]
If $(\epsilon_1,\epsilon_2)=(0,0)$, then the entire two-dimensional space $E_{6,36}$ is the transfer kernel.
\end{prop}

\begin{proof}
The target $\operatorname{Ext}_{\mathscr A}^{6,42}(\mathbb F_2,\mathbb F_2)$ is the one-dimensional space $\mathbb F_2\cdot t$, so the coefficients in \eqref{eq:epsilon-definition} exist and are unique. The coinvariant class of $f_*[M]$ has coordinates
\[
\overline\mu_{6,36}[M,f]
=a_1(M,f)e_1+a_2(M,f)e_2,
\]
because $e_1,e_2$ are dual to $[\zeta_1],[\zeta_2]$. Applying the transfer gives \eqref{eq:geometric-transfer-functional}. When at least one $\epsilon_i$ is nonzero, the kernel of the nonzero functional $(u_1,u_2)\mapsto\epsilon_1u_1+\epsilon_2u_2$ is the line generated by $(\epsilon_2,\epsilon_1)$, which proves \eqref{eq:kernel-line-epsilon}. The remaining case is immediate.
\end{proof}

Proposition~\ref{prop:geometric-transfer-kernel} gives the strongest transfer-kernel statement justified by the present computation. In geometric terms, a primitive bordism class belongs to the geometric kernel precisely when all mixed Wu numbers in \eqref{eq:mixed-wu-vanishing} vanish and the additional line-bundle characteristic number
\[
\epsilon_1\left\langle f^*(\zeta_1),[M]\right\rangle
+
\epsilon_2\left\langle f^*(\zeta_2),[M]\right\rangle
\]
vanishes. 

\subsection{A degree-$36$ Milnor hypersurface and the primitivity obstruction}

The next calculation adapts the characteristic-number method used for Milnor hypersurfaces to dimension $36$. It produces an explicit generator of the indecomposable quotient of the coefficient ring $\mathfrak N_*$, and then shows that the evident map defined by the two tautological line bundles does not lie in the primitive bordism subgroup needed for Singer's transfer.

Let
\[
H_{4,33}=\left\{([u],[z])\in\mathbb RP^4\times\mathbb RP^{33}
\;\middle|\;
\sum_{i=0}^4u_iz_i=0\right\}.
\]
This is a smooth closed hypersurface of dimension $36$. Let $\lambda$ and $\mu$ denote the restrictions of the tautological line bundles from the two factors, and put
\[
a=w_1(\lambda),\qquad b=w_1(\mu).
\]
For a real vector bundle $E$ with formal Stiefel--\allowbreak Whitney roots $t_j$, write $s_r(E)=\sum_jt_j^r$ for its $r$-th Newton class.

\begin{thm}[The Milnor hypersurface in dimension $36$]\label{thm:milnor-36}
The cobordism class of $H_{4,33}$ has nonzero image in the indecomposable quotient
\[
Q\mathfrak N_{36}
=\mathfrak N_{36}/\bigl((\mathfrak N_+)^2\bigr)_{36},
\]
and therefore generates this one-dimensional quotient. For the tautological classifying map
\[
\tau:H_{4,33}\longrightarrow BV_2\longrightarrow BV_6,
\]
one has
\[
\tau_*[H_{4,33}]
=v_1^{(3)}v_2^{(33)}+v_1^{(4)}v_2^{(32)}.
\tag{\ref{thm:milnor-36}.1}\label{eq:milnor-pushforward}
\]
This homology class is not $\mathscr A$-annihilated. More precisely,
\[
Sq_*^1\tau_*[H_{4,33}]
=v_1^{(3)}v_2^{(32)}+v_1^{(4)}v_2^{(31)}\neq0.
\tag{\ref{thm:milnor-36}.2}\label{eq:milnor-sq1}
\]
Consequently the tautological Milnor-hypersurface construction cannot represent either $e_1$ or $e_2$.
\end{thm}

\begin{proof}
The defining equation is a section of $\lambda\otimes\mu$. Choose unit representatives $u$ and $z$ at a zero of this section. Since $u\cdot z=0$, the nonzero vector $u$ is tangent to the sphere at $z$ and the derivative in that tangent direction is $\lVert u\rVert^2\neq0$. Hence the section is transverse to the zero section. The normal line bundle is therefore $\lambda\otimes\mu$. The standard stable tangent bundle identities for real projective space give
\[
TH_{4,33}\oplus(\lambda\otimes\mu)\oplus\underline{\mathbb R}^{\,2}
\cong5\lambda\oplus34\mu.
\tag{\ref{thm:milnor-36}.3}\label{eq:milnor-stable-tangent}
\]
The Newton classes are additive under direct sum. Since $w_1(\lambda\otimes\mu)=a+b$, relation \eqref{eq:milnor-stable-tangent} yields
\[
s_{36}(TH_{4,33})=5a^{36}+34b^{36}+(a+b)^{36}.
\]
The relation $a^5=0$ eliminates the first term, and the coefficient $34$ eliminates the second term over $\mathbb F_2$. The mod-$2$ fundamental class of the hypersurface is Poincar\'e dual in $\mathbb RP^4\times\mathbb RP^{33}$ to $a+b$. Hence
\[
\begin{aligned}
\left\langle s_{36}(TH_{4,33}),[H_{4,33}]\right\rangle
&=\left\langle(a+b)^{37},[\mathbb RP^4\times\mathbb RP^{33}]\right\rangle\\
&=\binom{37}{4}\equiv1\pmod2.
\end{aligned}
\tag{\ref{thm:milnor-36}.4}\label{eq:milnor-newton-number}
\]
Thom's computation of the unoriented cobordism ring, with Dold's explicit generators, gives
\[
\mathfrak N_*\cong
\mathbb F_2[y_j\mid j\geq1,\ j\neq2^r-1\text{ for every }r\geq1],
\qquad |y_j|=j
\]
\cite{Thom,Dold}. Thus $Q\mathfrak N_{36}$ is one-dimensional. The characteristic number $s_{36}$ vanishes on every product of positive-dimensional manifolds whose dimensions sum to $36$: by additivity of Newton classes, $s_{36}$ of such a product is pulled back from the factors, on each of which a degree-$36$ class is zero for dimensional reasons. The nonzero number in \eqref{eq:milnor-newton-number} therefore proves that $[H_{4,33}]$ is indecomposable and generates $Q\mathfrak N_{36}$.

For $u+v=36$, the defining hypersurface relation gives
\[
\begin{aligned}
\left\langle x_1^ux_2^v,\tau_*[H_{4,33}]\right\rangle
&=\left\langle a^ub^v,[H_{4,33}]\right\rangle\\
&=\left\langle a^ub^v(a+b),[\mathbb RP^4\times\mathbb RP^{33}]\right\rangle.
\end{aligned}
\]
The last expression is nonzero precisely for $(u,v)=(3,33)$ or $(4,32)$, which proves \eqref{eq:milnor-pushforward}. In one divided-power variable,
\[
Sq_*^r\bigl(v^{(d)}\bigr)=\binom{d-r}{r}v^{(d-r)}.
\tag{\ref{thm:milnor-36}.5}\label{eq:dual-square-divided-power}
\]
The Cartan formula and \eqref{eq:dual-square-divided-power} show that the first summand of \eqref{eq:milnor-pushforward} has zero $Sq_*^1$-image, while the second has image
\[
v_1^{(3)}v_2^{(32)}+v_1^{(4)}v_2^{(31)}.
\]
The two basis monomials are distinct, so no cancellation is possible, proving \eqref{eq:milnor-sq1}. More generally, for $m,n\geq2$, the tautological map of $H_{m,n}$ satisfies
\[
\tau_*[H_{m,n}]=v_1^{(m-1)}v_2^{(n)}+v_1^{(m)}v_2^{(n-1)}.
\]
If $m+n$ is odd, exactly one of $m$ and $n$ is even. Formula \eqref{eq:dual-square-divided-power} then shows that one of the two displayed summands has zero $Sq_*^1$-image, while the other has two distinct nonzero terms. Thus every even-dimensional Milnor hypersurface $H_{m,n}$ with $m,n\geq2$ has a nonprimitive tautological class over $BV_2$.

The same obstruction is visible as a characteristic number. Taking degree-one parts in \eqref{eq:milnor-stable-tangent} gives $w_1(TH_{4,33})=b$, and therefore
\[
\left\langle w_1(TH_{4,33})a^3b^{32},[H_{4,33}]\right\rangle
=1.
\]
By \eqref{eq:wu-adjunction}, this is exactly the nonzero pairing detecting \eqref{eq:milnor-sq1}. Hence $\tau_*[H_{4,33}]$ is not primitive and cannot represent either transfer-source basis class.
\end{proof}

Theorem~\ref{thm:milnor-36} exhibits the precise distinction between the coefficient-ring generator and the transfer-source classes. The manifold $H_{4,33}$ is geometrically explicit and indecomposable in $\mathfrak N_{36}$, while its tautological map to $BV_6$ fails the first mixed Wu-number condition. The abstract representatives in Corollary~\ref{corl:zeta-bordism-realization} must therefore involve different map data, different manifolds, or both.

\subsection{Products of projective spaces and Dold manifolds}

We next test two further standard families. For projective products, the divided-power description reduces primitivity to the binary spike condition. For Dold manifolds, the first cohomology has rank at most one, which forces every map to $BV_6$ through a rank-one elementary abelian subgroup and rules out any nonzero primitive class in degree $36$.

\begin{prop}[Obstructions for standard projective-product and Dold models]\label{prop:standard-model-obstructions}
Let
\[
P(\mathbf d)=\mathbb RP^{d_1}\times\cdots\times\mathbb RP^{d_r},
\qquad r\leq6,
\qquad d_i>0,
\qquad \sum_{i=1}^r d_i=36,
\]
and let $f_{\mathbf d}:P(\mathbf d)\to BV_6$ be the coordinate classifying map followed by a linear monomorphism $V_r\hookrightarrow V_6$. No such class represents $e_1$ or $e_2$. If $r<6$, its coinvariant coordinates against $[\zeta_1]$ and $[\zeta_2]$ are zero. If $r=6$, the pushforward can be primitive only when every $d_i$ is of the form $2^{s_i}-1$, and the only positive spike partitions of $36$ into six parts are
\[
(31,1,1,1,1,1),\quad
(15,15,3,1,1,1),\quad
(15,7,7,3,3,1),\quad
(7,7,7,7,7,1).
\tag{\ref{prop:standard-model-obstructions}.1}\label{eq:six-spike-partitions}
\]
None of the four partition types in \eqref{eq:six-spike-partitions} occurs in the support of $\zeta_1$ or $\zeta_2$.

Let
\[
D(m,n)=\bigl(S^m\times\mathbb CP^n\bigr)/\bigl((u,z)\sim(-u,\overline z)\bigr)
\]
be a Dold manifold of dimension $m+2n=36$. For every map $f:D(m,n)\to BV_6$, if $f_*[D(m,n)]$ is $\mathscr A$-annihilated, then its class in $E_{6,36}$ is zero. Hence no Dold manifold realizes $e_1$ or $e_2$ by a map to $BV_6$.
\end{prop}

\begin{proof}
For the coordinate map from a projective product, the fundamental class maps to
\[
v_1^{(d_1)}\cdots v_r^{(d_r)}
\]
in the homology of the coordinate copy of $BV_r$. If $r<6$, a change of basis in $V_6$ carries the image to the standard coordinate subgroup. Every monomial in the displayed representatives $\zeta_1$ and $\zeta_2$ contains every variable with positive exponent. Their restrictions to this coordinate $BV_r$ are therefore zero. Since their quotient classes are $GL(6)$-invariant, the resulting pairings are zero for every linear monomorphism $V_r\hookrightarrow V_6$.

Suppose $r=6$. A full-rank coordinate map differs from the standard one by an element of $GL(6)$, so its coinvariant class is represented by the pure divided-power monomial $v_1^{(d_1)}\cdots v_6^{(d_6)}$. Formula \eqref{eq:dual-square-divided-power} shows that $v^{(d)}$ is annihilated by every positive Steenrod square if $d=2^s-1$. Conversely, assume that $d$ is not of this form. In the binary expansion of $d$, choose the least position $s$ at which the digit is zero but a higher digit is one. Subtracting $2^s$ borrows from a higher position and makes the $s$-th digit of $d-2^s$ equal to one. Lucas' criterion then gives
\[
\binom{d-2^s}{2^s}\equiv1\pmod2,
\]
so $Sq_*^{2^s}(v^{(d)})\neq0$. In the Cartan expansion on a pure product, the term obtained by assigning the entire operation to the corresponding factor has a uniquely determined exponent vector and therefore cannot cancel with a different assignment. It follows that $v_1^{(d_1)}\cdots v_6^{(d_6)}$ is primitive if and only if every $d_i$ is a spike. Direct enumeration of six positive spikes with sum $36$ gives exactly the four partitions in \eqref{eq:six-spike-partitions}. Comparison with Tables~\ref{tab:zeta1-orbit-data} and~\ref{tab:zeta2-orbit-data}, equivalently with the complete expansions in \ref{app:explicit-polys}, shows that none occurs. The pairing with each $[\zeta_j]$ is consequently zero.

For a Dold manifold,
\[
H^*(D(m,n);\mathbb F_2)
\cong\mathbb F_2[c,d]/(c^{m+1},d^{n+1}),
\qquad |c|=1,\quad |d|=2
\]
\cite{Dold}. Thus $H^1(D(m,n);\mathbb F_2)$ has dimension at most one. Since $BV_6$ is an Eilenberg--Mac Lane space $K(V_6,1)$, every map from $D(m,n)$ to $BV_6$ has image contained, after a $GL(6)$-change of coordinates, in a copy of $BV_1$. The degree-$36$ homology of $BV_1$ is generated by $v^{(36)}$, and
\[
Sq_*^1(v^{(36)})=\binom{35}{1}v^{(35)}=v^{(35)}\neq0.
\]
Hence $\mathcal P_{\mathscr A}H_{36}(BV_1;\mathbb F_2)=0$. A primitive pushforward from a Dold manifold must therefore be zero, and so its class in $E_{6,36}$ is zero.
\end{proof}

The projective-product statement concerns the standard maps obtained from at most six tautological factors. It does not exclude a more elaborate realization built from a product with additional factors and a non-coordinate rank-six map. Such a construction would still have to satisfy all mixed Wu-number equations and the dual-basis conditions \eqref{eq:zeta-dual-numbers}.

\subsection{The Kameko homomorphism as a Thom-space Steenrod square}

The final geometric comparison concerns the map that drives the degree reduction in the algebraic computation. The appropriate rigorous interpretation is a Thom-space identity for the sum of the universal line bundles. It explains the factor $x_1\cdots x_q$ in the inverse Kameko monomial operation, but it does not supply a canonical fixed-point construction on manifolds.

Let
\[
\eta_q=\gamma_1\oplus\cdots\oplus\gamma_q\longrightarrow BV_q
\]
and let $U_q\in\widetilde H^q(\operatorname{Th}(\eta_q);\mathbb F_2)$ be its Thom class.

\begin{prop}[Thom-space realization of the inverse Kameko operation]\label{prop:kameko-thom-square}
For every homogeneous class $g\in H^d(BV_q;\mathbb F_2)$,
\[
Sq^{d+q}(gU_q)
=x_1\cdots x_q\,g^2U_q.
\tag{\ref{prop:kameko-thom-square}.1}\label{eq:kameko-thom-square}
\]
Under the Thom isomorphism, the cohomology class on the right corresponds exactly to
\[
\psi(g)=x_1\cdots x_q\,g^2,
\]
which on monomials sends
\[
x_1^{a_1}\cdots x_q^{a_q}
\longmapsto
x_1^{2a_1+1}\cdots x_q^{2a_q+1}.
\]
\end{prop}

\begin{proof}
The Thom identity gives
\[
Sq^j(U_q)=w_j(\eta_q)U_q.
\]
The Cartan formula therefore yields
\[
Sq^{d+q}(gU_q)
=\sum_{i+j=d+q}Sq^i(g)w_j(\eta_q)U_q.
\]
Instability gives $Sq^i(g)=0$ for $i>d$, while $w_j(\eta_q)=0$ for $j>q$. Since $i+j=d+q$, the only possibly nonzero summand has $i=d$ and $j=q$. For a homogeneous class of degree $d$, one has $Sq^d(g)=g^2$, and
\[
w_q(\eta_q)=\prod_{i=1}^q w_1(\gamma_i)=x_1\cdots x_q.
\]
Substitution proves \eqref{eq:kameko-thom-square}.
\end{proof}

Walker and Wood describe the transpose dual of the down Kameko map as the $GL(q)$-linear map
\[
\kappa^*\left(v_1^{(a_1)}\cdots v_q^{(a_q)}\right)
=v_1^{(2a_1+1)}\cdots v_q^{(2a_q+1)},
\]
and prove that it carries $\mathscr A$-annihilated classes to $\mathscr A$-annihilated classes \cite[Proposition~9.5.4]{WW}. Let
\[
E_{6,n}=\bigl(\mathbb F_2\otimes_{GL(6)}
\mathcal P_{\mathscr A}H_*(BV_6;\mathbb F_2)\bigr)_n.
\]

\begin{corl}[The degree-$15$ class has no degree-$36$ coinvariant Kameko lift]\label{corl:kameko-coinvariant-zero}
The map induced by the dual up Kameko operation,
\[
\overline\kappa^*:E_{6,15}\longrightarrow E_{6,36},
\]
is zero. If $p_\xi\in\mathcal P_{\mathscr A}H_{15}(BV_6;\mathbb F_2)$ represents the nonzero class dual to $[\xi]$, then $\kappa^*(p_\xi)$ is a nonzero primitive class in $H_{36}(BV_6;\mathbb F_2)$, but its $GL(6)$-coinvariant class is zero. In particular, no manifold realization of $\kappa^*(p_\xi)$ can represent $e_1$ or $e_2$.
\end{corl}

\begin{proof}
For $z\in[(QP_6)_{36}]^{GL(6)}$ and $p\in\mathcal P_{\mathscr A}H_{15}(BV_6;\mathbb F_2)$, transpose duality gives
\[
\left\langle z,\kappa^*(p)\right\rangle
=\left\langle(\widetilde{Sq}^0_*)_{(6,36)}z,p\right\rangle.
\]
Theorem~\ref{dlc2} and the final lift calculation in Section~\ref{s2} show that the invariant space in degree $36$ is spanned by $[\zeta_1]$ and $[\zeta_2]$, both of which lie in the Kameko kernel. The right-hand side is therefore zero for every invariant $z$. Since $E_{6,36}$ is the dual of the invariant space, the coinvariant class of $\kappa^*(p)$ is zero. The down Kameko map is surjective, so its transpose $\kappa^*$ is injective; hence $\kappa^*(p_\xi)\neq0$. Its primitivity follows from the Walker--Wood result cited above.
\end{proof}

If a $15$-manifold $N$ with a map $g:N\to BV_6$ realizes $p_\xi$, Thom's representability theorem supplies some closed $36$-manifold $M$ with a map $f:M\to BV_6$ realizing the homology class $\kappa^*(p_\xi)$. Corollary~\ref{corl:kameko-coinvariant-zero} shows that this algebraically raised class is invisible in the rank-six transfer source and therefore cannot be either dual generator in degree $36$. No canonical construction of $M$ from $N$ follows from this argument. In particular, a fixed-point interpretation would require an equivariant bordism refinement of $\kappa^*$ together with control of normal data; neither the quotient computation nor the dual Kameko formula provides such a refinement. The verified geometric content is precisely the Thom-space identity \eqref{eq:kameko-thom-square} and the coinvariant vanishing in Corollary~\ref{corl:kameko-coinvariant-zero}.

\section*{Funding}

\DD\d{\u a}ng V\~o Ph\'uc was funded by the Post-Doctoral Scholarship Programme of Vingroup Innovation Foundation (VINIF), Institute of Big Data, code: VINIF.2024.STS.38.




\section*{Data Availability}

The detailed computational outputs for the cases $(q,n)=(6,15)$ and $(q,n)=(6,36)$ are available at the Zenodo records referenced in Note~\ref{cyf}. The \texttt{Julia}/\texttt{OSCAR} source code, validation drivers for Dickson invariants and finite permutation modules, and orbit-type post-processing scripts are archived at the software record listed in~\ref{app:software}.

\appendix
\section{Complete polynomial expansions of the degree 36 representatives}\label{app:explicit-polys}

The orbit-type tables in the main text give the structural support decomposition of the two reduced representatives in Theorem~\ref{dlc2}. For completeness and for direct verification against the raw computational output, we record here the full polynomial representatives. Since the ground field is $\mathbb F_2$, every monomial displayed below has coefficient $1$, and $\zeta_s$ is exactly the sum of the listed monomials, for $s=1,2$. To keep the appendix compact, the monomials are arranged in four columns; the numbering is part of the display only and is not used in the algebraic expression.

\begingroup
\scriptsize
\setlength{\columnsep}{7pt}
\setlength{\parindent}{0pt}
\setlength{\parskip}{1.2pt}
\raggedcolumns
\newcommand{\appmon}[2]{\noindent\(\scriptstyle #1.\;#2\)\par}

\subsection*{Expansion of $\zeta_1$}
\begin{multicols}{4}
\appmon{1}{x_{1}\,\allowbreak{}x_{2}\,\allowbreak{}x_{3}\,\allowbreak{}x_{4}\,\allowbreak{}x_{5}^{30}\,\allowbreak{}x_{6}^{2}}
\appmon{2}{x_{1}\,\allowbreak{}x_{2}\,\allowbreak{}x_{3}\,\allowbreak{}x_{4}\,\allowbreak{}x_{5}^{6}\,\allowbreak{}x_{6}^{26}}
\appmon{3}{x_{1}\,\allowbreak{}x_{2}\,\allowbreak{}x_{3}\,\allowbreak{}x_{4}^{30}\,\allowbreak{}x_{5}\,\allowbreak{}x_{6}^{2}}
\appmon{4}{x_{1}\,\allowbreak{}x_{2}\,\allowbreak{}x_{3}\,\allowbreak{}x_{4}^{6}\,\allowbreak{}x_{5}\,\allowbreak{}x_{6}^{26}}
\appmon{5}{x_{1}\,\allowbreak{}x_{2}\,\allowbreak{}x_{3}^{3}\,\allowbreak{}x_{4}\,\allowbreak{}x_{5}^{28}\,\allowbreak{}x_{6}^{2}}
\appmon{6}{x_{1}\,\allowbreak{}x_{2}\,\allowbreak{}x_{3}^{3}\,\allowbreak{}x_{4}^{28}\,\allowbreak{}x_{5}\,\allowbreak{}x_{6}^{2}}
\appmon{7}{x_{1}\,\allowbreak{}x_{2}\,\allowbreak{}x_{3}^{6}\,\allowbreak{}x_{4}\,\allowbreak{}x_{5}\,\allowbreak{}x_{6}^{26}}
\appmon{8}{x_{1}\,\allowbreak{}x_{2}^{3}\,\allowbreak{}x_{3}\,\allowbreak{}x_{4}\,\allowbreak{}x_{5}^{28}\,\allowbreak{}x_{6}^{2}}
\appmon{9}{x_{1}\,\allowbreak{}x_{2}^{3}\,\allowbreak{}x_{3}\,\allowbreak{}x_{4}^{28}\,\allowbreak{}x_{5}\,\allowbreak{}x_{6}^{2}}
\appmon{10}{x_{1}\,\allowbreak{}x_{2}^{3}\,\allowbreak{}x_{3}^{3}\,\allowbreak{}x_{4}\,\allowbreak{}x_{5}^{4}\,\allowbreak{}x_{6}^{24}}
\appmon{11}{x_{1}\,\allowbreak{}x_{2}^{3}\,\allowbreak{}x_{3}^{3}\,\allowbreak{}x_{4}^{4}\,\allowbreak{}x_{5}\,\allowbreak{}x_{6}^{24}}
\appmon{12}{x_{1}\,\allowbreak{}x_{2}^{3}\,\allowbreak{}x_{3}^{5}\,\allowbreak{}x_{4}\,\allowbreak{}x_{5}^{24}\,\allowbreak{}x_{6}^{2}}
\appmon{13}{x_{1}\,\allowbreak{}x_{2}^{3}\,\allowbreak{}x_{3}^{5}\,\allowbreak{}x_{4}^{24}\,\allowbreak{}x_{5}\,\allowbreak{}x_{6}^{2}}
\appmon{14}{x_{1}\,\allowbreak{}x_{2}^{3}\,\allowbreak{}x_{3}^{6}\,\allowbreak{}x_{4}\,\allowbreak{}x_{5}\,\allowbreak{}x_{6}^{24}}
\appmon{15}{x_{1}\,\allowbreak{}x_{2}^{6}\,\allowbreak{}x_{3}\,\allowbreak{}x_{4}\,\allowbreak{}x_{5}\,\allowbreak{}x_{6}^{26}}
\appmon{16}{x_{1}^{3}\,\allowbreak{}x_{2}\,\allowbreak{}x_{3}\,\allowbreak{}x_{4}\,\allowbreak{}x_{5}^{28}\,\allowbreak{}x_{6}^{2}}
\appmon{17}{x_{1}^{3}\,\allowbreak{}x_{2}\,\allowbreak{}x_{3}\,\allowbreak{}x_{4}\,\allowbreak{}x_{5}^{4}\,\allowbreak{}x_{6}^{26}}
\appmon{18}{x_{1}^{3}\,\allowbreak{}x_{2}\,\allowbreak{}x_{3}\,\allowbreak{}x_{4}^{28}\,\allowbreak{}x_{5}\,\allowbreak{}x_{6}^{2}}
\appmon{19}{x_{1}^{3}\,\allowbreak{}x_{2}\,\allowbreak{}x_{3}\,\allowbreak{}x_{4}^{4}\,\allowbreak{}x_{5}\,\allowbreak{}x_{6}^{26}}
\appmon{20}{x_{1}^{3}\,\allowbreak{}x_{2}\,\allowbreak{}x_{3}^{3}\,\allowbreak{}x_{4}\,\allowbreak{}x_{5}^{4}\,\allowbreak{}x_{6}^{24}}
\appmon{21}{x_{1}^{3}\,\allowbreak{}x_{2}\,\allowbreak{}x_{3}^{3}\,\allowbreak{}x_{4}^{4}\,\allowbreak{}x_{5}\,\allowbreak{}x_{6}^{24}}
\appmon{22}{x_{1}^{3}\,\allowbreak{}x_{2}\,\allowbreak{}x_{3}^{4}\,\allowbreak{}x_{4}\,\allowbreak{}x_{5}\,\allowbreak{}x_{6}^{26}}
\appmon{23}{x_{1}^{3}\,\allowbreak{}x_{2}\,\allowbreak{}x_{3}^{5}\,\allowbreak{}x_{4}\,\allowbreak{}x_{5}^{24}\,\allowbreak{}x_{6}^{2}}
\appmon{24}{x_{1}^{3}\,\allowbreak{}x_{2}\,\allowbreak{}x_{3}^{5}\,\allowbreak{}x_{4}^{24}\,\allowbreak{}x_{5}\,\allowbreak{}x_{6}^{2}}
\appmon{25}{x_{1}^{3}\,\allowbreak{}x_{2}\,\allowbreak{}x_{3}^{6}\,\allowbreak{}x_{4}\,\allowbreak{}x_{5}\,\allowbreak{}x_{6}^{24}}
\appmon{26}{x_{1}^{3}\,\allowbreak{}x_{2}^{3}\,\allowbreak{}x_{3}\,\allowbreak{}x_{4}\,\allowbreak{}x_{5}^{4}\,\allowbreak{}x_{6}^{24}}
\appmon{27}{x_{1}^{3}\,\allowbreak{}x_{2}^{3}\,\allowbreak{}x_{3}\,\allowbreak{}x_{4}^{4}\,\allowbreak{}x_{5}\,\allowbreak{}x_{6}^{24}}
\appmon{28}{x_{1}^{3}\,\allowbreak{}x_{2}^{3}\,\allowbreak{}x_{3}^{4}\,\allowbreak{}x_{4}\,\allowbreak{}x_{5}\,\allowbreak{}x_{6}^{24}}
\appmon{29}{x_{1}^{3}\,\allowbreak{}x_{2}^{4}\,\allowbreak{}x_{3}\,\allowbreak{}x_{4}\,\allowbreak{}x_{5}\,\allowbreak{}x_{6}^{26}}
\appmon{30}{x_{1}^{3}\,\allowbreak{}x_{2}^{5}\,\allowbreak{}x_{3}\,\allowbreak{}x_{4}\,\allowbreak{}x_{5}^{2}\,\allowbreak{}x_{6}^{24}}
\appmon{31}{x_{1}^{3}\,\allowbreak{}x_{2}^{5}\,\allowbreak{}x_{3}\,\allowbreak{}x_{4}\,\allowbreak{}x_{5}^{24}\,\allowbreak{}x_{6}^{2}}
\appmon{32}{x_{1}^{3}\,\allowbreak{}x_{2}^{5}\,\allowbreak{}x_{3}\,\allowbreak{}x_{4}^{2}\,\allowbreak{}x_{5}\,\allowbreak{}x_{6}^{24}}
\appmon{33}{x_{1}^{3}\,\allowbreak{}x_{2}^{5}\,\allowbreak{}x_{3}\,\allowbreak{}x_{4}^{24}\,\allowbreak{}x_{5}\,\allowbreak{}x_{6}^{2}}
\appmon{34}{x_{1}^{3}\,\allowbreak{}x_{2}^{5}\,\allowbreak{}x_{3}^{10}\,\allowbreak{}x_{4}\,\allowbreak{}x_{5}\,\allowbreak{}x_{6}^{16}}
\appmon{35}{x_{1}^{3}\,\allowbreak{}x_{2}^{5}\,\allowbreak{}x_{3}^{2}\,\allowbreak{}x_{4}\,\allowbreak{}x_{5}\,\allowbreak{}x_{6}^{24}}
\appmon{36}{x_{1}^{3}\,\allowbreak{}x_{2}^{5}\,\allowbreak{}x_{3}^{3}\,\allowbreak{}x_{4}\,\allowbreak{}x_{5}^{8}\,\allowbreak{}x_{6}^{16}}
\appmon{37}{x_{1}^{3}\,\allowbreak{}x_{2}^{5}\,\allowbreak{}x_{3}^{3}\,\allowbreak{}x_{4}^{8}\,\allowbreak{}x_{5}\,\allowbreak{}x_{6}^{16}}
\appmon{38}{x_{1}^{3}\,\allowbreak{}x_{2}^{5}\,\allowbreak{}x_{3}^{9}\,\allowbreak{}x_{4}\,\allowbreak{}x_{5}^{16}\,\allowbreak{}x_{6}^{2}}
\appmon{39}{x_{1}^{3}\,\allowbreak{}x_{2}^{5}\,\allowbreak{}x_{3}^{9}\,\allowbreak{}x_{4}^{16}\,\allowbreak{}x_{5}\,\allowbreak{}x_{6}^{2}}
\end{multicols}

\subsection*{Expansion of $\zeta_2$}
\begin{multicols}{4}
\appmon{1}{x_{1}\,\allowbreak{}x_{2}\,\allowbreak{}x_{3}\,\allowbreak{}x_{4}\,\allowbreak{}x_{5}^{2}\,\allowbreak{}x_{6}^{30}}
\appmon{2}{x_{1}\,\allowbreak{}x_{2}\,\allowbreak{}x_{3}\,\allowbreak{}x_{4}\,\allowbreak{}x_{5}^{6}\,\allowbreak{}x_{6}^{26}}
\appmon{3}{x_{1}\,\allowbreak{}x_{2}\,\allowbreak{}x_{3}\,\allowbreak{}x_{4}^{2}\,\allowbreak{}x_{5}^{3}\,\allowbreak{}x_{6}^{28}}
\appmon{4}{x_{1}\,\allowbreak{}x_{2}\,\allowbreak{}x_{3}\,\allowbreak{}x_{4}^{2}\,\allowbreak{}x_{5}^{4}\,\allowbreak{}x_{6}^{27}}
\appmon{5}{x_{1}\,\allowbreak{}x_{2}\,\allowbreak{}x_{3}\,\allowbreak{}x_{4}^{2}\,\allowbreak{}x_{5}^{6}\,\allowbreak{}x_{6}^{25}}
\appmon{6}{x_{1}\,\allowbreak{}x_{2}\,\allowbreak{}x_{3}\,\allowbreak{}x_{4}^{2}\,\allowbreak{}x_{5}^{7}\,\allowbreak{}x_{6}^{24}}
\appmon{7}{x_{1}\,\allowbreak{}x_{2}\,\allowbreak{}x_{3}\,\allowbreak{}x_{4}^{3}\,\allowbreak{}x_{5}^{4}\,\allowbreak{}x_{6}^{26}}
\appmon{8}{x_{1}\,\allowbreak{}x_{2}\,\allowbreak{}x_{3}\,\allowbreak{}x_{4}^{3}\,\allowbreak{}x_{5}^{6}\,\allowbreak{}x_{6}^{24}}
\appmon{9}{x_{1}\,\allowbreak{}x_{2}\,\allowbreak{}x_{3}\,\allowbreak{}x_{4}^{6}\,\allowbreak{}x_{5}^{3}\,\allowbreak{}x_{6}^{24}}
\appmon{10}{x_{1}\,\allowbreak{}x_{2}\,\allowbreak{}x_{3}\,\allowbreak{}x_{4}^{7}\,\allowbreak{}x_{5}^{2}\,\allowbreak{}x_{6}^{24}}
\appmon{11}{x_{1}\,\allowbreak{}x_{2}\,\allowbreak{}x_{3}^{2}\,\allowbreak{}x_{4}^{3}\,\allowbreak{}x_{5}^{5}\,\allowbreak{}x_{6}^{24}}
\appmon{12}{x_{1}\,\allowbreak{}x_{2}\,\allowbreak{}x_{3}^{2}\,\allowbreak{}x_{4}^{4}\,\allowbreak{}x_{5}^{3}\,\allowbreak{}x_{6}^{25}}
\appmon{13}{x_{1}\,\allowbreak{}x_{2}\,\allowbreak{}x_{3}^{2}\,\allowbreak{}x_{4}^{5}\,\allowbreak{}x_{5}^{10}\,\allowbreak{}x_{6}^{17}}
\appmon{14}{x_{1}\,\allowbreak{}x_{2}\,\allowbreak{}x_{3}^{2}\,\allowbreak{}x_{4}^{5}\,\allowbreak{}x_{5}^{3}\,\allowbreak{}x_{6}^{24}}
\appmon{15}{x_{1}\,\allowbreak{}x_{2}\,\allowbreak{}x_{3}^{2}\,\allowbreak{}x_{4}^{5}\,\allowbreak{}x_{5}^{9}\,\allowbreak{}x_{6}^{18}}
\appmon{16}{x_{1}\,\allowbreak{}x_{2}\,\allowbreak{}x_{3}^{2}\,\allowbreak{}x_{4}^{6}\,\allowbreak{}x_{5}^{9}\,\allowbreak{}x_{6}^{17}}
\appmon{17}{x_{1}\,\allowbreak{}x_{2}\,\allowbreak{}x_{3}^{2}\,\allowbreak{}x_{4}^{7}\,\allowbreak{}x_{5}^{8}\,\allowbreak{}x_{6}^{17}}
\appmon{18}{x_{1}\,\allowbreak{}x_{2}\,\allowbreak{}x_{3}^{2}\,\allowbreak{}x_{4}^{7}\,\allowbreak{}x_{5}^{9}\,\allowbreak{}x_{6}^{16}}
\appmon{19}{x_{1}\,\allowbreak{}x_{2}\,\allowbreak{}x_{3}^{3}\,\allowbreak{}x_{4}\,\allowbreak{}x_{5}^{6}\,\allowbreak{}x_{6}^{24}}
\appmon{20}{x_{1}\,\allowbreak{}x_{2}\,\allowbreak{}x_{3}^{3}\,\allowbreak{}x_{4}^{2}\,\allowbreak{}x_{5}\,\allowbreak{}x_{6}^{28}}
\appmon{21}{x_{1}\,\allowbreak{}x_{2}\,\allowbreak{}x_{3}^{3}\,\allowbreak{}x_{4}^{3}\,\allowbreak{}x_{5}^{14}\,\allowbreak{}x_{6}^{14}}
\appmon{22}{x_{1}\,\allowbreak{}x_{2}\,\allowbreak{}x_{3}^{3}\,\allowbreak{}x_{4}^{5}\,\allowbreak{}x_{5}^{2}\,\allowbreak{}x_{6}^{24}}
\appmon{23}{x_{1}\,\allowbreak{}x_{2}\,\allowbreak{}x_{3}^{3}\,\allowbreak{}x_{4}^{5}\,\allowbreak{}x_{5}^{8}\,\allowbreak{}x_{6}^{18}}
\appmon{24}{x_{1}\,\allowbreak{}x_{2}\,\allowbreak{}x_{3}^{3}\,\allowbreak{}x_{4}^{6}\,\allowbreak{}x_{5}^{11}\,\allowbreak{}x_{6}^{14}}
\appmon{25}{x_{1}\,\allowbreak{}x_{2}\,\allowbreak{}x_{3}^{3}\,\allowbreak{}x_{4}^{6}\,\allowbreak{}x_{5}^{12}\,\allowbreak{}x_{6}^{13}}
\appmon{26}{x_{1}\,\allowbreak{}x_{2}\,\allowbreak{}x_{3}^{3}\,\allowbreak{}x_{4}^{6}\,\allowbreak{}x_{5}^{13}\,\allowbreak{}x_{6}^{12}}
\appmon{27}{x_{1}\,\allowbreak{}x_{2}\,\allowbreak{}x_{3}^{3}\,\allowbreak{}x_{4}^{6}\,\allowbreak{}x_{5}^{14}\,\allowbreak{}x_{6}^{11}}
\appmon{28}{x_{1}\,\allowbreak{}x_{2}\,\allowbreak{}x_{3}^{3}\,\allowbreak{}x_{4}^{6}\,\allowbreak{}x_{5}^{8}\,\allowbreak{}x_{6}^{17}}
\appmon{29}{x_{1}\,\allowbreak{}x_{2}\,\allowbreak{}x_{3}^{6}\,\allowbreak{}x_{4}^{10}\,\allowbreak{}x_{5}^{13}\,\allowbreak{}x_{6}^{5}}
\appmon{30}{x_{1}\,\allowbreak{}x_{2}\,\allowbreak{}x_{3}^{6}\,\allowbreak{}x_{4}^{10}\,\allowbreak{}x_{5}^{5}\,\allowbreak{}x_{6}^{13}}
\appmon{31}{x_{1}\,\allowbreak{}x_{2}\,\allowbreak{}x_{3}^{6}\,\allowbreak{}x_{4}^{11}\,\allowbreak{}x_{5}^{5}\,\allowbreak{}x_{6}^{12}}
\appmon{32}{x_{1}\,\allowbreak{}x_{2}\,\allowbreak{}x_{3}^{6}\,\allowbreak{}x_{4}^{2}\,\allowbreak{}x_{5}^{9}\,\allowbreak{}x_{6}^{17}}
\appmon{33}{x_{1}\,\allowbreak{}x_{2}\,\allowbreak{}x_{3}^{6}\,\allowbreak{}x_{4}^{3}\,\allowbreak{}x_{5}^{12}\,\allowbreak{}x_{6}^{13}}
\appmon{34}{x_{1}\,\allowbreak{}x_{2}\,\allowbreak{}x_{3}^{6}\,\allowbreak{}x_{4}^{3}\,\allowbreak{}x_{5}^{13}\,\allowbreak{}x_{6}^{12}}
\appmon{35}{x_{1}\,\allowbreak{}x_{2}\,\allowbreak{}x_{3}^{6}\,\allowbreak{}x_{4}^{3}\,\allowbreak{}x_{5}^{9}\,\allowbreak{}x_{6}^{16}}
\appmon{36}{x_{1}\,\allowbreak{}x_{2}\,\allowbreak{}x_{3}^{7}\,\allowbreak{}x_{4}\,\allowbreak{}x_{5}^{2}\,\allowbreak{}x_{6}^{24}}
\appmon{37}{x_{1}\,\allowbreak{}x_{2}\,\allowbreak{}x_{3}^{7}\,\allowbreak{}x_{4}^{10}\,\allowbreak{}x_{5}^{5}\,\allowbreak{}x_{6}^{12}}
\appmon{38}{x_{1}\,\allowbreak{}x_{2}\,\allowbreak{}x_{3}^{7}\,\allowbreak{}x_{4}^{2}\,\allowbreak{}x_{5}\,\allowbreak{}x_{6}^{24}}
\appmon{39}{x_{1}\,\allowbreak{}x_{2}\,\allowbreak{}x_{3}^{7}\,\allowbreak{}x_{4}^{2}\,\allowbreak{}x_{5}^{8}\,\allowbreak{}x_{6}^{17}}
\appmon{40}{x_{1}\,\allowbreak{}x_{2}\,\allowbreak{}x_{3}^{7}\,\allowbreak{}x_{4}^{2}\,\allowbreak{}x_{5}^{9}\,\allowbreak{}x_{6}^{16}}
\appmon{41}{x_{1}\,\allowbreak{}x_{2}\,\allowbreak{}x_{3}^{7}\,\allowbreak{}x_{4}^{3}\,\allowbreak{}x_{5}^{12}\,\allowbreak{}x_{6}^{12}}
\appmon{42}{x_{1}\,\allowbreak{}x_{2}^{14}\,\allowbreak{}x_{3}^{3}\,\allowbreak{}x_{4}^{5}\,\allowbreak{}x_{5}^{3}\,\allowbreak{}x_{6}^{10}}
\appmon{43}{x_{1}\,\allowbreak{}x_{2}^{3}\,\allowbreak{}x_{3}\,\allowbreak{}x_{4}\,\allowbreak{}x_{5}^{6}\,\allowbreak{}x_{6}^{24}}
\appmon{44}{x_{1}\,\allowbreak{}x_{2}^{3}\,\allowbreak{}x_{3}\,\allowbreak{}x_{4}^{2}\,\allowbreak{}x_{5}^{4}\,\allowbreak{}x_{6}^{25}}
\appmon{45}{x_{1}\,\allowbreak{}x_{2}^{3}\,\allowbreak{}x_{3}\,\allowbreak{}x_{4}^{2}\,\allowbreak{}x_{5}^{5}\,\allowbreak{}x_{6}^{24}}
\appmon{46}{x_{1}\,\allowbreak{}x_{2}^{3}\,\allowbreak{}x_{3}\,\allowbreak{}x_{4}^{3}\,\allowbreak{}x_{5}^{14}\,\allowbreak{}x_{6}^{14}}
\appmon{47}{x_{1}\,\allowbreak{}x_{2}^{3}\,\allowbreak{}x_{3}\,\allowbreak{}x_{4}^{3}\,\allowbreak{}x_{5}^{4}\,\allowbreak{}x_{6}^{24}}
\appmon{48}{x_{1}\,\allowbreak{}x_{2}^{3}\,\allowbreak{}x_{3}\,\allowbreak{}x_{4}^{4}\,\allowbreak{}x_{5}^{9}\,\allowbreak{}x_{6}^{18}}
\appmon{49}{x_{1}\,\allowbreak{}x_{2}^{3}\,\allowbreak{}x_{3}\,\allowbreak{}x_{4}^{5}\,\allowbreak{}x_{5}^{2}\,\allowbreak{}x_{6}^{24}}
\appmon{50}{x_{1}\,\allowbreak{}x_{2}^{3}\,\allowbreak{}x_{3}\,\allowbreak{}x_{4}^{6}\,\allowbreak{}x_{5}^{11}\,\allowbreak{}x_{6}^{14}}
\appmon{51}{x_{1}\,\allowbreak{}x_{2}^{3}\,\allowbreak{}x_{3}\,\allowbreak{}x_{4}^{6}\,\allowbreak{}x_{5}^{12}\,\allowbreak{}x_{6}^{13}}
\appmon{52}{x_{1}\,\allowbreak{}x_{2}^{3}\,\allowbreak{}x_{3}\,\allowbreak{}x_{4}^{6}\,\allowbreak{}x_{5}^{13}\,\allowbreak{}x_{6}^{12}}
\appmon{53}{x_{1}\,\allowbreak{}x_{2}^{3}\,\allowbreak{}x_{3}\,\allowbreak{}x_{4}^{6}\,\allowbreak{}x_{5}^{14}\,\allowbreak{}x_{6}^{11}}
\appmon{54}{x_{1}\,\allowbreak{}x_{2}^{3}\,\allowbreak{}x_{3}\,\allowbreak{}x_{4}^{6}\,\allowbreak{}x_{5}^{9}\,\allowbreak{}x_{6}^{16}}
\appmon{55}{x_{1}\,\allowbreak{}x_{2}^{3}\,\allowbreak{}x_{3}^{14}\,\allowbreak{}x_{4}^{5}\,\allowbreak{}x_{5}^{3}\,\allowbreak{}x_{6}^{10}}
\appmon{56}{x_{1}\,\allowbreak{}x_{2}^{3}\,\allowbreak{}x_{3}^{3}\,\allowbreak{}x_{4}\,\allowbreak{}x_{5}^{14}\,\allowbreak{}x_{6}^{14}}
\appmon{57}{x_{1}\,\allowbreak{}x_{2}^{3}\,\allowbreak{}x_{3}^{3}\,\allowbreak{}x_{4}\,\allowbreak{}x_{5}^{4}\,\allowbreak{}x_{6}^{24}}
\appmon{58}{x_{1}\,\allowbreak{}x_{2}^{3}\,\allowbreak{}x_{3}^{3}\,\allowbreak{}x_{4}^{12}\,\allowbreak{}x_{5}^{5}\,\allowbreak{}x_{6}^{12}}
\appmon{59}{x_{1}\,\allowbreak{}x_{2}^{3}\,\allowbreak{}x_{3}^{3}\,\allowbreak{}x_{4}^{13}\,\allowbreak{}x_{5}^{6}\,\allowbreak{}x_{6}^{10}}
\appmon{60}{x_{1}\,\allowbreak{}x_{2}^{3}\,\allowbreak{}x_{3}^{3}\,\allowbreak{}x_{4}^{14}\,\allowbreak{}x_{5}^{3}\,\allowbreak{}x_{6}^{12}}
\appmon{61}{x_{1}\,\allowbreak{}x_{2}^{3}\,\allowbreak{}x_{3}^{3}\,\allowbreak{}x_{4}^{14}\,\allowbreak{}x_{5}^{5}\,\allowbreak{}x_{6}^{10}}
\appmon{62}{x_{1}\,\allowbreak{}x_{2}^{3}\,\allowbreak{}x_{3}^{3}\,\allowbreak{}x_{4}^{4}\,\allowbreak{}x_{5}\,\allowbreak{}x_{6}^{24}}
\appmon{63}{x_{1}\,\allowbreak{}x_{2}^{3}\,\allowbreak{}x_{3}^{3}\,\allowbreak{}x_{4}^{5}\,\allowbreak{}x_{5}^{12}\,\allowbreak{}x_{6}^{12}}
\appmon{64}{x_{1}\,\allowbreak{}x_{2}^{3}\,\allowbreak{}x_{3}^{3}\,\allowbreak{}x_{4}^{7}\,\allowbreak{}x_{5}^{10}\,\allowbreak{}x_{6}^{12}}
\appmon{65}{x_{1}\,\allowbreak{}x_{2}^{3}\,\allowbreak{}x_{3}^{3}\,\allowbreak{}x_{4}^{7}\,\allowbreak{}x_{5}^{12}\,\allowbreak{}x_{6}^{10}}
\appmon{66}{x_{1}\,\allowbreak{}x_{2}^{3}\,\allowbreak{}x_{3}^{4}\,\allowbreak{}x_{4}\,\allowbreak{}x_{5}^{10}\,\allowbreak{}x_{6}^{17}}
\appmon{67}{x_{1}\,\allowbreak{}x_{2}^{3}\,\allowbreak{}x_{3}^{4}\,\allowbreak{}x_{4}\,\allowbreak{}x_{5}^{11}\,\allowbreak{}x_{6}^{16}}
\appmon{68}{x_{1}\,\allowbreak{}x_{2}^{3}\,\allowbreak{}x_{3}^{4}\,\allowbreak{}x_{4}\,\allowbreak{}x_{5}^{8}\,\allowbreak{}x_{6}^{19}}
\appmon{69}{x_{1}\,\allowbreak{}x_{2}^{3}\,\allowbreak{}x_{3}^{4}\,\allowbreak{}x_{4}\,\allowbreak{}x_{5}^{9}\,\allowbreak{}x_{6}^{18}}
\appmon{70}{x_{1}\,\allowbreak{}x_{2}^{3}\,\allowbreak{}x_{3}^{4}\,\allowbreak{}x_{4}^{10}\,\allowbreak{}x_{5}^{9}\,\allowbreak{}x_{6}^{9}}
\appmon{71}{x_{1}\,\allowbreak{}x_{2}^{3}\,\allowbreak{}x_{3}^{4}\,\allowbreak{}x_{4}^{11}\,\allowbreak{}x_{5}^{8}\,\allowbreak{}x_{6}^{9}}
\appmon{72}{x_{1}\,\allowbreak{}x_{2}^{3}\,\allowbreak{}x_{3}^{4}\,\allowbreak{}x_{4}^{11}\,\allowbreak{}x_{5}^{9}\,\allowbreak{}x_{6}^{8}}
\appmon{73}{x_{1}\,\allowbreak{}x_{2}^{3}\,\allowbreak{}x_{3}^{4}\,\allowbreak{}x_{4}^{3}\,\allowbreak{}x_{5}^{8}\,\allowbreak{}x_{6}^{17}}
\appmon{74}{x_{1}\,\allowbreak{}x_{2}^{3}\,\allowbreak{}x_{3}^{4}\,\allowbreak{}x_{4}^{8}\,\allowbreak{}x_{5}^{3}\,\allowbreak{}x_{6}^{17}}
\appmon{75}{x_{1}\,\allowbreak{}x_{2}^{3}\,\allowbreak{}x_{3}^{4}\,\allowbreak{}x_{4}^{9}\,\allowbreak{}x_{5}^{10}\,\allowbreak{}x_{6}^{9}}
\appmon{76}{x_{1}\,\allowbreak{}x_{2}^{3}\,\allowbreak{}x_{3}^{4}\,\allowbreak{}x_{4}^{9}\,\allowbreak{}x_{5}^{11}\,\allowbreak{}x_{6}^{8}}
\appmon{77}{x_{1}\,\allowbreak{}x_{2}^{3}\,\allowbreak{}x_{3}^{4}\,\allowbreak{}x_{4}^{9}\,\allowbreak{}x_{5}^{8}\,\allowbreak{}x_{6}^{11}}
\appmon{78}{x_{1}\,\allowbreak{}x_{2}^{3}\,\allowbreak{}x_{3}^{5}\,\allowbreak{}x_{4}\,\allowbreak{}x_{5}^{2}\,\allowbreak{}x_{6}^{24}}
\appmon{79}{x_{1}\,\allowbreak{}x_{2}^{3}\,\allowbreak{}x_{3}^{5}\,\allowbreak{}x_{4}\,\allowbreak{}x_{5}^{8}\,\allowbreak{}x_{6}^{18}}
\appmon{80}{x_{1}\,\allowbreak{}x_{2}^{3}\,\allowbreak{}x_{3}^{5}\,\allowbreak{}x_{4}^{10}\,\allowbreak{}x_{5}^{5}\,\allowbreak{}x_{6}^{12}}
\appmon{81}{x_{1}\,\allowbreak{}x_{2}^{3}\,\allowbreak{}x_{3}^{5}\,\allowbreak{}x_{4}^{10}\,\allowbreak{}x_{5}^{8}\,\allowbreak{}x_{6}^{9}}
\appmon{82}{x_{1}\,\allowbreak{}x_{2}^{3}\,\allowbreak{}x_{3}^{5}\,\allowbreak{}x_{4}^{10}\,\allowbreak{}x_{5}^{9}\,\allowbreak{}x_{6}^{8}}
\appmon{83}{x_{1}\,\allowbreak{}x_{2}^{3}\,\allowbreak{}x_{3}^{5}\,\allowbreak{}x_{4}^{11}\,\allowbreak{}x_{5}^{6}\,\allowbreak{}x_{6}^{10}}
\appmon{84}{x_{1}\,\allowbreak{}x_{2}^{3}\,\allowbreak{}x_{3}^{5}\,\allowbreak{}x_{4}^{14}\,\allowbreak{}x_{5}^{3}\,\allowbreak{}x_{6}^{10}}
\appmon{85}{x_{1}\,\allowbreak{}x_{2}^{3}\,\allowbreak{}x_{3}^{5}\,\allowbreak{}x_{4}^{2}\,\allowbreak{}x_{5}^{8}\,\allowbreak{}x_{6}^{17}}
\appmon{86}{x_{1}\,\allowbreak{}x_{2}^{3}\,\allowbreak{}x_{3}^{5}\,\allowbreak{}x_{4}^{2}\,\allowbreak{}x_{5}^{9}\,\allowbreak{}x_{6}^{16}}
\appmon{87}{x_{1}\,\allowbreak{}x_{2}^{3}\,\allowbreak{}x_{3}^{5}\,\allowbreak{}x_{4}^{3}\,\allowbreak{}x_{5}^{12}\,\allowbreak{}x_{6}^{12}}
\appmon{88}{x_{1}\,\allowbreak{}x_{2}^{3}\,\allowbreak{}x_{3}^{5}\,\allowbreak{}x_{4}^{6}\,\allowbreak{}x_{5}^{10}\,\allowbreak{}x_{6}^{11}}
\appmon{89}{x_{1}\,\allowbreak{}x_{2}^{3}\,\allowbreak{}x_{3}^{5}\,\allowbreak{}x_{4}^{6}\,\allowbreak{}x_{5}^{11}\,\allowbreak{}x_{6}^{10}}
\appmon{90}{x_{1}\,\allowbreak{}x_{2}^{3}\,\allowbreak{}x_{3}^{5}\,\allowbreak{}x_{4}^{6}\,\allowbreak{}x_{5}^{14}\,\allowbreak{}x_{6}^{7}}
\appmon{91}{x_{1}\,\allowbreak{}x_{2}^{3}\,\allowbreak{}x_{3}^{5}\,\allowbreak{}x_{4}^{6}\,\allowbreak{}x_{5}^{7}\,\allowbreak{}x_{6}^{14}}
\appmon{92}{x_{1}\,\allowbreak{}x_{2}^{3}\,\allowbreak{}x_{3}^{5}\,\allowbreak{}x_{4}^{8}\,\allowbreak{}x_{5}^{10}\,\allowbreak{}x_{6}^{9}}
\appmon{93}{x_{1}\,\allowbreak{}x_{2}^{3}\,\allowbreak{}x_{3}^{5}\,\allowbreak{}x_{4}^{8}\,\allowbreak{}x_{5}^{11}\,\allowbreak{}x_{6}^{8}}
\appmon{94}{x_{1}\,\allowbreak{}x_{2}^{3}\,\allowbreak{}x_{3}^{5}\,\allowbreak{}x_{4}^{8}\,\allowbreak{}x_{5}^{3}\,\allowbreak{}x_{6}^{16}}
\appmon{95}{x_{1}\,\allowbreak{}x_{2}^{3}\,\allowbreak{}x_{3}^{5}\,\allowbreak{}x_{4}^{8}\,\allowbreak{}x_{5}^{8}\,\allowbreak{}x_{6}^{11}}
\appmon{96}{x_{1}\,\allowbreak{}x_{2}^{3}\,\allowbreak{}x_{3}^{5}\,\allowbreak{}x_{4}^{9}\,\allowbreak{}x_{5}^{10}\,\allowbreak{}x_{6}^{8}}
\appmon{97}{x_{1}\,\allowbreak{}x_{2}^{3}\,\allowbreak{}x_{3}^{5}\,\allowbreak{}x_{4}^{9}\,\allowbreak{}x_{5}^{2}\,\allowbreak{}x_{6}^{16}}
\appmon{98}{x_{1}\,\allowbreak{}x_{2}^{3}\,\allowbreak{}x_{3}^{6}\,\allowbreak{}x_{4}\,\allowbreak{}x_{5}^{11}\,\allowbreak{}x_{6}^{14}}
\appmon{99}{x_{1}\,\allowbreak{}x_{2}^{3}\,\allowbreak{}x_{3}^{6}\,\allowbreak{}x_{4}\,\allowbreak{}x_{5}^{14}\,\allowbreak{}x_{6}^{11}}
\appmon{100}{x_{1}\,\allowbreak{}x_{2}^{3}\,\allowbreak{}x_{3}^{6}\,\allowbreak{}x_{4}\,\allowbreak{}x_{5}^{8}\,\allowbreak{}x_{6}^{17}}
\appmon{101}{x_{1}\,\allowbreak{}x_{2}^{3}\,\allowbreak{}x_{3}^{6}\,\allowbreak{}x_{4}\,\allowbreak{}x_{5}^{9}\,\allowbreak{}x_{6}^{16}}
\appmon{102}{x_{1}\,\allowbreak{}x_{2}^{3}\,\allowbreak{}x_{3}^{6}\,\allowbreak{}x_{4}^{11}\,\allowbreak{}x_{5}^{3}\,\allowbreak{}x_{6}^{12}}
\appmon{103}{x_{1}\,\allowbreak{}x_{2}^{3}\,\allowbreak{}x_{3}^{6}\,\allowbreak{}x_{4}^{11}\,\allowbreak{}x_{5}^{5}\,\allowbreak{}x_{6}^{10}}
\appmon{104}{x_{1}\,\allowbreak{}x_{2}^{3}\,\allowbreak{}x_{3}^{6}\,\allowbreak{}x_{4}^{13}\,\allowbreak{}x_{5}^{3}\,\allowbreak{}x_{6}^{10}}
\appmon{105}{x_{1}\,\allowbreak{}x_{2}^{3}\,\allowbreak{}x_{3}^{6}\,\allowbreak{}x_{4}^{3}\,\allowbreak{}x_{5}^{10}\,\allowbreak{}x_{6}^{13}}
\appmon{106}{x_{1}\,\allowbreak{}x_{2}^{3}\,\allowbreak{}x_{3}^{6}\,\allowbreak{}x_{4}^{3}\,\allowbreak{}x_{5}^{11}\,\allowbreak{}x_{6}^{12}}
\appmon{107}{x_{1}\,\allowbreak{}x_{2}^{3}\,\allowbreak{}x_{3}^{6}\,\allowbreak{}x_{4}^{3}\,\allowbreak{}x_{5}^{12}\,\allowbreak{}x_{6}^{11}}
\appmon{108}{x_{1}\,\allowbreak{}x_{2}^{3}\,\allowbreak{}x_{3}^{6}\,\allowbreak{}x_{4}^{3}\,\allowbreak{}x_{5}^{13}\,\allowbreak{}x_{6}^{10}}
\appmon{109}{x_{1}\,\allowbreak{}x_{2}^{3}\,\allowbreak{}x_{3}^{6}\,\allowbreak{}x_{4}^{3}\,\allowbreak{}x_{5}^{14}\,\allowbreak{}x_{6}^{9}}
\appmon{110}{x_{1}\,\allowbreak{}x_{2}^{3}\,\allowbreak{}x_{3}^{6}\,\allowbreak{}x_{4}^{3}\,\allowbreak{}x_{5}^{9}\,\allowbreak{}x_{6}^{14}}
\appmon{111}{x_{1}\,\allowbreak{}x_{2}^{3}\,\allowbreak{}x_{3}^{6}\,\allowbreak{}x_{4}^{5}\,\allowbreak{}x_{5}^{10}\,\allowbreak{}x_{6}^{11}}
\appmon{112}{x_{1}\,\allowbreak{}x_{2}^{3}\,\allowbreak{}x_{3}^{6}\,\allowbreak{}x_{4}^{5}\,\allowbreak{}x_{5}^{11}\,\allowbreak{}x_{6}^{10}}
\appmon{113}{x_{1}\,\allowbreak{}x_{2}^{3}\,\allowbreak{}x_{3}^{6}\,\allowbreak{}x_{4}^{5}\,\allowbreak{}x_{5}^{13}\,\allowbreak{}x_{6}^{8}}
\appmon{114}{x_{1}\,\allowbreak{}x_{2}^{3}\,\allowbreak{}x_{3}^{6}\,\allowbreak{}x_{4}^{5}\,\allowbreak{}x_{5}^{14}\,\allowbreak{}x_{6}^{7}}
\appmon{115}{x_{1}\,\allowbreak{}x_{2}^{3}\,\allowbreak{}x_{3}^{6}\,\allowbreak{}x_{4}^{5}\,\allowbreak{}x_{5}^{7}\,\allowbreak{}x_{6}^{14}}
\appmon{116}{x_{1}\,\allowbreak{}x_{2}^{3}\,\allowbreak{}x_{3}^{6}\,\allowbreak{}x_{4}^{5}\,\allowbreak{}x_{5}^{8}\,\allowbreak{}x_{6}^{13}}
\appmon{117}{x_{1}\,\allowbreak{}x_{2}^{3}\,\allowbreak{}x_{3}^{6}\,\allowbreak{}x_{4}^{7}\,\allowbreak{}x_{5}^{9}\,\allowbreak{}x_{6}^{10}}
\appmon{118}{x_{1}\,\allowbreak{}x_{2}^{3}\,\allowbreak{}x_{3}^{6}\,\allowbreak{}x_{4}^{8}\,\allowbreak{}x_{5}^{9}\,\allowbreak{}x_{6}^{9}}
\appmon{119}{x_{1}\,\allowbreak{}x_{2}^{3}\,\allowbreak{}x_{3}^{6}\,\allowbreak{}x_{4}^{9}\,\allowbreak{}x_{5}\,\allowbreak{}x_{6}^{16}}
\appmon{120}{x_{1}\,\allowbreak{}x_{2}^{3}\,\allowbreak{}x_{3}^{6}\,\allowbreak{}x_{4}^{9}\,\allowbreak{}x_{5}^{8}\,\allowbreak{}x_{6}^{9}}
\appmon{121}{x_{1}\,\allowbreak{}x_{2}^{3}\,\allowbreak{}x_{3}^{6}\,\allowbreak{}x_{4}^{9}\,\allowbreak{}x_{5}^{9}\,\allowbreak{}x_{6}^{8}}
\appmon{122}{x_{1}\,\allowbreak{}x_{2}^{3}\,\allowbreak{}x_{3}^{7}\,\allowbreak{}x_{4}\,\allowbreak{}x_{5}^{8}\,\allowbreak{}x_{6}^{16}}
\appmon{123}{x_{1}\,\allowbreak{}x_{2}^{3}\,\allowbreak{}x_{3}^{7}\,\allowbreak{}x_{4}^{12}\,\allowbreak{}x_{5}\,\allowbreak{}x_{6}^{12}}
\appmon{124}{x_{1}\,\allowbreak{}x_{2}^{3}\,\allowbreak{}x_{3}^{7}\,\allowbreak{}x_{4}^{12}\,\allowbreak{}x_{5}^{3}\,\allowbreak{}x_{6}^{10}}
\appmon{125}{x_{1}\,\allowbreak{}x_{2}^{3}\,\allowbreak{}x_{3}^{7}\,\allowbreak{}x_{4}^{3}\,\allowbreak{}x_{5}^{12}\,\allowbreak{}x_{6}^{10}}
\appmon{126}{x_{1}\,\allowbreak{}x_{2}^{3}\,\allowbreak{}x_{3}^{7}\,\allowbreak{}x_{4}^{5}\,\allowbreak{}x_{5}^{10}\,\allowbreak{}x_{6}^{10}}
\appmon{127}{x_{1}\,\allowbreak{}x_{2}^{3}\,\allowbreak{}x_{3}^{7}\,\allowbreak{}x_{4}^{5}\,\allowbreak{}x_{5}^{8}\,\allowbreak{}x_{6}^{12}}
\appmon{128}{x_{1}\,\allowbreak{}x_{2}^{3}\,\allowbreak{}x_{3}^{7}\,\allowbreak{}x_{4}^{6}\,\allowbreak{}x_{5}^{9}\,\allowbreak{}x_{6}^{10}}
\appmon{129}{x_{1}\,\allowbreak{}x_{2}^{3}\,\allowbreak{}x_{3}^{7}\,\allowbreak{}x_{4}^{8}\,\allowbreak{}x_{5}\,\allowbreak{}x_{6}^{16}}
\appmon{130}{x_{1}\,\allowbreak{}x_{2}^{3}\,\allowbreak{}x_{3}^{7}\,\allowbreak{}x_{4}^{8}\,\allowbreak{}x_{5}^{8}\,\allowbreak{}x_{6}^{9}}
\appmon{131}{x_{1}\,\allowbreak{}x_{2}^{3}\,\allowbreak{}x_{3}^{7}\,\allowbreak{}x_{4}^{8}\,\allowbreak{}x_{5}^{9}\,\allowbreak{}x_{6}^{8}}
\appmon{132}{x_{1}\,\allowbreak{}x_{2}^{6}\,\allowbreak{}x_{3}\,\allowbreak{}x_{4}^{11}\,\allowbreak{}x_{5}^{5}\,\allowbreak{}x_{6}^{12}}
\appmon{133}{x_{1}\,\allowbreak{}x_{2}^{6}\,\allowbreak{}x_{3}\,\allowbreak{}x_{4}^{2}\,\allowbreak{}x_{5}^{9}\,\allowbreak{}x_{6}^{17}}
\appmon{134}{x_{1}\,\allowbreak{}x_{2}^{6}\,\allowbreak{}x_{3}\,\allowbreak{}x_{4}^{3}\,\allowbreak{}x_{5}^{12}\,\allowbreak{}x_{6}^{13}}
\appmon{135}{x_{1}\,\allowbreak{}x_{2}^{6}\,\allowbreak{}x_{3}\,\allowbreak{}x_{4}^{3}\,\allowbreak{}x_{5}^{13}\,\allowbreak{}x_{6}^{12}}
\appmon{136}{x_{1}\,\allowbreak{}x_{2}^{6}\,\allowbreak{}x_{3}\,\allowbreak{}x_{4}^{3}\,\allowbreak{}x_{5}^{9}\,\allowbreak{}x_{6}^{16}}
\appmon{137}{x_{1}\,\allowbreak{}x_{2}^{6}\,\allowbreak{}x_{3}^{11}\,\allowbreak{}x_{4}^{5}\,\allowbreak{}x_{5}\,\allowbreak{}x_{6}^{12}}
\appmon{138}{x_{1}\,\allowbreak{}x_{2}^{6}\,\allowbreak{}x_{3}^{11}\,\allowbreak{}x_{4}^{5}\,\allowbreak{}x_{5}^{3}\,\allowbreak{}x_{6}^{10}}
\appmon{139}{x_{1}\,\allowbreak{}x_{2}^{6}\,\allowbreak{}x_{3}^{3}\,\allowbreak{}x_{4}^{13}\,\allowbreak{}x_{5}\,\allowbreak{}x_{6}^{12}}
\appmon{140}{x_{1}\,\allowbreak{}x_{2}^{6}\,\allowbreak{}x_{3}^{3}\,\allowbreak{}x_{4}^{5}\,\allowbreak{}x_{5}^{10}\,\allowbreak{}x_{6}^{11}}
\appmon{141}{x_{1}\,\allowbreak{}x_{2}^{6}\,\allowbreak{}x_{3}^{3}\,\allowbreak{}x_{4}^{5}\,\allowbreak{}x_{5}^{11}\,\allowbreak{}x_{6}^{10}}
\appmon{142}{x_{1}\,\allowbreak{}x_{2}^{6}\,\allowbreak{}x_{3}^{3}\,\allowbreak{}x_{4}^{5}\,\allowbreak{}x_{5}^{14}\,\allowbreak{}x_{6}^{7}}
\appmon{143}{x_{1}\,\allowbreak{}x_{2}^{6}\,\allowbreak{}x_{3}^{3}\,\allowbreak{}x_{4}^{5}\,\allowbreak{}x_{5}^{7}\,\allowbreak{}x_{6}^{14}}
\appmon{144}{x_{1}\,\allowbreak{}x_{2}^{6}\,\allowbreak{}x_{3}^{3}\,\allowbreak{}x_{4}^{5}\,\allowbreak{}x_{5}^{9}\,\allowbreak{}x_{6}^{12}}
\appmon{145}{x_{1}\,\allowbreak{}x_{2}^{6}\,\allowbreak{}x_{3}^{3}\,\allowbreak{}x_{4}^{7}\,\allowbreak{}x_{5}^{9}\,\allowbreak{}x_{6}^{10}}
\appmon{146}{x_{1}\,\allowbreak{}x_{2}^{6}\,\allowbreak{}x_{3}^{3}\,\allowbreak{}x_{4}^{9}\,\allowbreak{}x_{5}^{5}\,\allowbreak{}x_{6}^{12}}
\appmon{147}{x_{1}\,\allowbreak{}x_{2}^{6}\,\allowbreak{}x_{3}^{7}\,\allowbreak{}x_{4}^{3}\,\allowbreak{}x_{5}^{9}\,\allowbreak{}x_{6}^{10}}
\appmon{148}{x_{1}\,\allowbreak{}x_{2}^{6}\,\allowbreak{}x_{3}^{7}\,\allowbreak{}x_{4}^{9}\,\allowbreak{}x_{5}^{3}\,\allowbreak{}x_{6}^{10}}
\appmon{149}{x_{1}\,\allowbreak{}x_{2}^{7}\,\allowbreak{}x_{3}\,\allowbreak{}x_{4}\,\allowbreak{}x_{5}^{2}\,\allowbreak{}x_{6}^{24}}
\appmon{150}{x_{1}\,\allowbreak{}x_{2}^{7}\,\allowbreak{}x_{3}\,\allowbreak{}x_{4}^{10}\,\allowbreak{}x_{5}^{12}\,\allowbreak{}x_{6}^{5}}
\appmon{151}{x_{1}\,\allowbreak{}x_{2}^{7}\,\allowbreak{}x_{3}\,\allowbreak{}x_{4}^{10}\,\allowbreak{}x_{5}^{5}\,\allowbreak{}x_{6}^{12}}
\appmon{152}{x_{1}\,\allowbreak{}x_{2}^{7}\,\allowbreak{}x_{3}\,\allowbreak{}x_{4}^{3}\,\allowbreak{}x_{5}^{12}\,\allowbreak{}x_{6}^{12}}
\appmon{153}{x_{1}\,\allowbreak{}x_{2}^{7}\,\allowbreak{}x_{3}^{10}\,\allowbreak{}x_{4}^{5}\,\allowbreak{}x_{5}\,\allowbreak{}x_{6}^{12}}
\appmon{154}{x_{1}\,\allowbreak{}x_{2}^{7}\,\allowbreak{}x_{3}^{10}\,\allowbreak{}x_{4}^{5}\,\allowbreak{}x_{5}^{3}\,\allowbreak{}x_{6}^{10}}
\appmon{155}{x_{1}\,\allowbreak{}x_{2}^{7}\,\allowbreak{}x_{3}^{3}\,\allowbreak{}x_{4}\,\allowbreak{}x_{5}^{8}\,\allowbreak{}x_{6}^{16}}
\appmon{156}{x_{1}\,\allowbreak{}x_{2}^{7}\,\allowbreak{}x_{3}^{3}\,\allowbreak{}x_{4}^{12}\,\allowbreak{}x_{5}\,\allowbreak{}x_{6}^{12}}
\appmon{157}{x_{1}\,\allowbreak{}x_{2}^{7}\,\allowbreak{}x_{3}^{3}\,\allowbreak{}x_{4}^{12}\,\allowbreak{}x_{5}^{3}\,\allowbreak{}x_{6}^{10}}
\appmon{158}{x_{1}\,\allowbreak{}x_{2}^{7}\,\allowbreak{}x_{3}^{3}\,\allowbreak{}x_{4}^{3}\,\allowbreak{}x_{5}^{12}\,\allowbreak{}x_{6}^{10}}
\appmon{159}{x_{1}\,\allowbreak{}x_{2}^{7}\,\allowbreak{}x_{3}^{3}\,\allowbreak{}x_{4}^{5}\,\allowbreak{}x_{5}^{10}\,\allowbreak{}x_{6}^{10}}
\appmon{160}{x_{1}\,\allowbreak{}x_{2}^{7}\,\allowbreak{}x_{3}^{3}\,\allowbreak{}x_{4}^{5}\,\allowbreak{}x_{5}^{8}\,\allowbreak{}x_{6}^{12}}
\appmon{161}{x_{1}\,\allowbreak{}x_{2}^{7}\,\allowbreak{}x_{3}^{3}\,\allowbreak{}x_{4}^{6}\,\allowbreak{}x_{5}^{9}\,\allowbreak{}x_{6}^{10}}
\appmon{162}{x_{1}\,\allowbreak{}x_{2}^{7}\,\allowbreak{}x_{3}^{3}\,\allowbreak{}x_{4}^{8}\,\allowbreak{}x_{5}^{8}\,\allowbreak{}x_{6}^{9}}
\appmon{163}{x_{1}\,\allowbreak{}x_{2}^{7}\,\allowbreak{}x_{3}^{3}\,\allowbreak{}x_{4}^{8}\,\allowbreak{}x_{5}^{9}\,\allowbreak{}x_{6}^{8}}
\appmon{164}{x_{1}\,\allowbreak{}x_{2}^{7}\,\allowbreak{}x_{3}^{3}\,\allowbreak{}x_{4}^{9}\,\allowbreak{}x_{5}^{8}\,\allowbreak{}x_{6}^{8}}
\appmon{165}{x_{1}\,\allowbreak{}x_{2}^{7}\,\allowbreak{}x_{3}^{6}\,\allowbreak{}x_{4}^{3}\,\allowbreak{}x_{5}^{9}\,\allowbreak{}x_{6}^{10}}
\appmon{166}{x_{1}^{3}\,\allowbreak{}x_{2}\,\allowbreak{}x_{3}\,\allowbreak{}x_{4}\,\allowbreak{}x_{5}^{6}\,\allowbreak{}x_{6}^{24}}
\appmon{167}{x_{1}^{3}\,\allowbreak{}x_{2}\,\allowbreak{}x_{3}\,\allowbreak{}x_{4}^{3}\,\allowbreak{}x_{5}^{14}\,\allowbreak{}x_{6}^{14}}
\appmon{168}{x_{1}^{3}\,\allowbreak{}x_{2}\,\allowbreak{}x_{3}\,\allowbreak{}x_{4}^{4}\,\allowbreak{}x_{5}^{9}\,\allowbreak{}x_{6}^{18}}
\appmon{169}{x_{1}^{3}\,\allowbreak{}x_{2}\,\allowbreak{}x_{3}\,\allowbreak{}x_{4}^{5}\,\allowbreak{}x_{5}^{2}\,\allowbreak{}x_{6}^{24}}
\appmon{170}{x_{1}^{3}\,\allowbreak{}x_{2}\,\allowbreak{}x_{3}\,\allowbreak{}x_{4}^{6}\,\allowbreak{}x_{5}^{11}\,\allowbreak{}x_{6}^{14}}
\appmon{171}{x_{1}^{3}\,\allowbreak{}x_{2}\,\allowbreak{}x_{3}\,\allowbreak{}x_{4}^{6}\,\allowbreak{}x_{5}^{14}\,\allowbreak{}x_{6}^{11}}
\appmon{172}{x_{1}^{3}\,\allowbreak{}x_{2}\,\allowbreak{}x_{3}\,\allowbreak{}x_{4}^{6}\,\allowbreak{}x_{5}^{9}\,\allowbreak{}x_{6}^{16}}
\appmon{173}{x_{1}^{3}\,\allowbreak{}x_{2}\,\allowbreak{}x_{3}\,\allowbreak{}x_{4}^{7}\,\allowbreak{}x_{5}^{8}\,\allowbreak{}x_{6}^{16}}
\appmon{174}{x_{1}^{3}\,\allowbreak{}x_{2}\,\allowbreak{}x_{3}^{3}\,\allowbreak{}x_{4}\,\allowbreak{}x_{5}^{14}\,\allowbreak{}x_{6}^{14}}
\appmon{175}{x_{1}^{3}\,\allowbreak{}x_{2}\,\allowbreak{}x_{3}^{3}\,\allowbreak{}x_{4}\,\allowbreak{}x_{5}^{4}\,\allowbreak{}x_{6}^{24}}
\appmon{176}{x_{1}^{3}\,\allowbreak{}x_{2}\,\allowbreak{}x_{3}^{3}\,\allowbreak{}x_{4}^{12}\,\allowbreak{}x_{5}^{5}\,\allowbreak{}x_{6}^{12}}
\appmon{177}{x_{1}^{3}\,\allowbreak{}x_{2}\,\allowbreak{}x_{3}^{3}\,\allowbreak{}x_{4}^{13}\,\allowbreak{}x_{5}^{4}\,\allowbreak{}x_{6}^{12}}
\appmon{178}{x_{1}^{3}\,\allowbreak{}x_{2}\,\allowbreak{}x_{3}^{3}\,\allowbreak{}x_{4}^{13}\,\allowbreak{}x_{5}^{6}\,\allowbreak{}x_{6}^{10}}
\appmon{179}{x_{1}^{3}\,\allowbreak{}x_{2}\,\allowbreak{}x_{3}^{3}\,\allowbreak{}x_{4}^{14}\,\allowbreak{}x_{5}^{3}\,\allowbreak{}x_{6}^{12}}
\appmon{180}{x_{1}^{3}\,\allowbreak{}x_{2}\,\allowbreak{}x_{3}^{3}\,\allowbreak{}x_{4}^{14}\,\allowbreak{}x_{5}^{5}\,\allowbreak{}x_{6}^{10}}
\appmon{181}{x_{1}^{3}\,\allowbreak{}x_{2}\,\allowbreak{}x_{3}^{3}\,\allowbreak{}x_{4}^{4}\,\allowbreak{}x_{5}^{11}\,\allowbreak{}x_{6}^{14}}
\appmon{182}{x_{1}^{3}\,\allowbreak{}x_{2}\,\allowbreak{}x_{3}^{3}\,\allowbreak{}x_{4}^{4}\,\allowbreak{}x_{5}^{14}\,\allowbreak{}x_{6}^{11}}
\appmon{183}{x_{1}^{3}\,\allowbreak{}x_{2}\,\allowbreak{}x_{3}^{3}\,\allowbreak{}x_{4}^{4}\,\allowbreak{}x_{5}^{8}\,\allowbreak{}x_{6}^{17}}
\appmon{184}{x_{1}^{3}\,\allowbreak{}x_{2}\,\allowbreak{}x_{3}^{3}\,\allowbreak{}x_{4}^{5}\,\allowbreak{}x_{5}^{8}\,\allowbreak{}x_{6}^{16}}
\appmon{185}{x_{1}^{3}\,\allowbreak{}x_{2}\,\allowbreak{}x_{3}^{3}\,\allowbreak{}x_{4}^{7}\,\allowbreak{}x_{5}^{10}\,\allowbreak{}x_{6}^{12}}
\appmon{186}{x_{1}^{3}\,\allowbreak{}x_{2}\,\allowbreak{}x_{3}^{4}\,\allowbreak{}x_{4}\,\allowbreak{}x_{5}^{10}\,\allowbreak{}x_{6}^{17}}
\appmon{187}{x_{1}^{3}\,\allowbreak{}x_{2}\,\allowbreak{}x_{3}^{4}\,\allowbreak{}x_{4}\,\allowbreak{}x_{5}^{11}\,\allowbreak{}x_{6}^{16}}
\appmon{188}{x_{1}^{3}\,\allowbreak{}x_{2}\,\allowbreak{}x_{3}^{4}\,\allowbreak{}x_{4}\,\allowbreak{}x_{5}^{8}\,\allowbreak{}x_{6}^{19}}
\appmon{189}{x_{1}^{3}\,\allowbreak{}x_{2}\,\allowbreak{}x_{3}^{4}\,\allowbreak{}x_{4}\,\allowbreak{}x_{5}^{9}\,\allowbreak{}x_{6}^{18}}
\appmon{190}{x_{1}^{3}\,\allowbreak{}x_{2}\,\allowbreak{}x_{3}^{4}\,\allowbreak{}x_{4}^{10}\,\allowbreak{}x_{5}^{9}\,\allowbreak{}x_{6}^{9}}
\appmon{191}{x_{1}^{3}\,\allowbreak{}x_{2}\,\allowbreak{}x_{3}^{4}\,\allowbreak{}x_{4}^{11}\,\allowbreak{}x_{5}^{5}\,\allowbreak{}x_{6}^{12}}
\appmon{192}{x_{1}^{3}\,\allowbreak{}x_{2}\,\allowbreak{}x_{3}^{4}\,\allowbreak{}x_{4}^{11}\,\allowbreak{}x_{5}^{8}\,\allowbreak{}x_{6}^{9}}
\appmon{193}{x_{1}^{3}\,\allowbreak{}x_{2}\,\allowbreak{}x_{3}^{4}\,\allowbreak{}x_{4}^{11}\,\allowbreak{}x_{5}^{9}\,\allowbreak{}x_{6}^{8}}
\appmon{194}{x_{1}^{3}\,\allowbreak{}x_{2}\,\allowbreak{}x_{3}^{4}\,\allowbreak{}x_{4}^{2}\,\allowbreak{}x_{5}^{9}\,\allowbreak{}x_{6}^{17}}
\appmon{195}{x_{1}^{3}\,\allowbreak{}x_{2}\,\allowbreak{}x_{3}^{4}\,\allowbreak{}x_{4}^{3}\,\allowbreak{}x_{5}^{12}\,\allowbreak{}x_{6}^{13}}
\appmon{196}{x_{1}^{3}\,\allowbreak{}x_{2}\,\allowbreak{}x_{3}^{4}\,\allowbreak{}x_{4}^{3}\,\allowbreak{}x_{5}^{13}\,\allowbreak{}x_{6}^{12}}
\appmon{197}{x_{1}^{3}\,\allowbreak{}x_{2}\,\allowbreak{}x_{3}^{4}\,\allowbreak{}x_{4}^{3}\,\allowbreak{}x_{5}^{8}\,\allowbreak{}x_{6}^{17}}
\appmon{198}{x_{1}^{3}\,\allowbreak{}x_{2}\,\allowbreak{}x_{3}^{4}\,\allowbreak{}x_{4}^{8}\,\allowbreak{}x_{5}^{3}\,\allowbreak{}x_{6}^{17}}
\appmon{199}{x_{1}^{3}\,\allowbreak{}x_{2}\,\allowbreak{}x_{3}^{4}\,\allowbreak{}x_{4}^{9}\,\allowbreak{}x_{5}^{10}\,\allowbreak{}x_{6}^{9}}
\appmon{200}{x_{1}^{3}\,\allowbreak{}x_{2}\,\allowbreak{}x_{3}^{4}\,\allowbreak{}x_{4}^{9}\,\allowbreak{}x_{5}^{11}\,\allowbreak{}x_{6}^{8}}
\appmon{201}{x_{1}^{3}\,\allowbreak{}x_{2}\,\allowbreak{}x_{3}^{4}\,\allowbreak{}x_{4}^{9}\,\allowbreak{}x_{5}^{3}\,\allowbreak{}x_{6}^{16}}
\appmon{202}{x_{1}^{3}\,\allowbreak{}x_{2}\,\allowbreak{}x_{3}^{4}\,\allowbreak{}x_{4}^{9}\,\allowbreak{}x_{5}^{8}\,\allowbreak{}x_{6}^{11}}
\appmon{203}{x_{1}^{3}\,\allowbreak{}x_{2}\,\allowbreak{}x_{3}^{5}\,\allowbreak{}x_{4}\,\allowbreak{}x_{5}^{2}\,\allowbreak{}x_{6}^{24}}
\appmon{204}{x_{1}^{3}\,\allowbreak{}x_{2}\,\allowbreak{}x_{3}^{5}\,\allowbreak{}x_{4}\,\allowbreak{}x_{5}^{8}\,\allowbreak{}x_{6}^{18}}
\appmon{205}{x_{1}^{3}\,\allowbreak{}x_{2}\,\allowbreak{}x_{3}^{5}\,\allowbreak{}x_{4}^{10}\,\allowbreak{}x_{5}^{5}\,\allowbreak{}x_{6}^{12}}
\appmon{206}{x_{1}^{3}\,\allowbreak{}x_{2}\,\allowbreak{}x_{3}^{5}\,\allowbreak{}x_{4}^{10}\,\allowbreak{}x_{5}^{8}\,\allowbreak{}x_{6}^{9}}
\appmon{207}{x_{1}^{3}\,\allowbreak{}x_{2}\,\allowbreak{}x_{3}^{5}\,\allowbreak{}x_{4}^{11}\,\allowbreak{}x_{5}^{6}\,\allowbreak{}x_{6}^{10}}
\appmon{208}{x_{1}^{3}\,\allowbreak{}x_{2}\,\allowbreak{}x_{3}^{5}\,\allowbreak{}x_{4}^{14}\,\allowbreak{}x_{5}^{3}\,\allowbreak{}x_{6}^{10}}
\appmon{209}{x_{1}^{3}\,\allowbreak{}x_{2}\,\allowbreak{}x_{3}^{5}\,\allowbreak{}x_{4}^{2}\,\allowbreak{}x_{5}^{8}\,\allowbreak{}x_{6}^{17}}
\appmon{210}{x_{1}^{3}\,\allowbreak{}x_{2}\,\allowbreak{}x_{3}^{5}\,\allowbreak{}x_{4}^{2}\,\allowbreak{}x_{5}^{9}\,\allowbreak{}x_{6}^{16}}
\appmon{211}{x_{1}^{3}\,\allowbreak{}x_{2}\,\allowbreak{}x_{3}^{5}\,\allowbreak{}x_{4}^{3}\,\allowbreak{}x_{5}^{8}\,\allowbreak{}x_{6}^{16}}
\appmon{212}{x_{1}^{3}\,\allowbreak{}x_{2}\,\allowbreak{}x_{3}^{5}\,\allowbreak{}x_{4}^{6}\,\allowbreak{}x_{5}^{10}\,\allowbreak{}x_{6}^{11}}
\appmon{213}{x_{1}^{3}\,\allowbreak{}x_{2}\,\allowbreak{}x_{3}^{5}\,\allowbreak{}x_{4}^{6}\,\allowbreak{}x_{5}^{11}\,\allowbreak{}x_{6}^{10}}
\appmon{214}{x_{1}^{3}\,\allowbreak{}x_{2}\,\allowbreak{}x_{3}^{5}\,\allowbreak{}x_{4}^{6}\,\allowbreak{}x_{5}^{14}\,\allowbreak{}x_{6}^{7}}
\appmon{215}{x_{1}^{3}\,\allowbreak{}x_{2}\,\allowbreak{}x_{3}^{5}\,\allowbreak{}x_{4}^{6}\,\allowbreak{}x_{5}^{7}\,\allowbreak{}x_{6}^{14}}
\appmon{216}{x_{1}^{3}\,\allowbreak{}x_{2}\,\allowbreak{}x_{3}^{5}\,\allowbreak{}x_{4}^{8}\,\allowbreak{}x_{5}^{9}\,\allowbreak{}x_{6}^{10}}
\appmon{217}{x_{1}^{3}\,\allowbreak{}x_{2}\,\allowbreak{}x_{3}^{5}\,\allowbreak{}x_{4}^{9}\,\allowbreak{}x_{5}^{2}\,\allowbreak{}x_{6}^{16}}
\appmon{218}{x_{1}^{3}\,\allowbreak{}x_{2}\,\allowbreak{}x_{3}^{5}\,\allowbreak{}x_{4}^{9}\,\allowbreak{}x_{5}^{8}\,\allowbreak{}x_{6}^{10}}
\appmon{219}{x_{1}^{3}\,\allowbreak{}x_{2}\,\allowbreak{}x_{3}^{6}\,\allowbreak{}x_{4}\,\allowbreak{}x_{5}^{11}\,\allowbreak{}x_{6}^{14}}
\appmon{220}{x_{1}^{3}\,\allowbreak{}x_{2}\,\allowbreak{}x_{3}^{6}\,\allowbreak{}x_{4}\,\allowbreak{}x_{5}^{14}\,\allowbreak{}x_{6}^{11}}
\appmon{221}{x_{1}^{3}\,\allowbreak{}x_{2}\,\allowbreak{}x_{3}^{6}\,\allowbreak{}x_{4}\,\allowbreak{}x_{5}^{8}\,\allowbreak{}x_{6}^{17}}
\appmon{222}{x_{1}^{3}\,\allowbreak{}x_{2}\,\allowbreak{}x_{3}^{6}\,\allowbreak{}x_{4}\,\allowbreak{}x_{5}^{9}\,\allowbreak{}x_{6}^{16}}
\appmon{223}{x_{1}^{3}\,\allowbreak{}x_{2}\,\allowbreak{}x_{3}^{6}\,\allowbreak{}x_{4}^{11}\,\allowbreak{}x_{5}^{3}\,\allowbreak{}x_{6}^{12}}
\appmon{224}{x_{1}^{3}\,\allowbreak{}x_{2}\,\allowbreak{}x_{3}^{6}\,\allowbreak{}x_{4}^{11}\,\allowbreak{}x_{5}^{5}\,\allowbreak{}x_{6}^{10}}
\appmon{225}{x_{1}^{3}\,\allowbreak{}x_{2}\,\allowbreak{}x_{3}^{6}\,\allowbreak{}x_{4}^{13}\,\allowbreak{}x_{5}\,\allowbreak{}x_{6}^{12}}
\appmon{226}{x_{1}^{3}\,\allowbreak{}x_{2}\,\allowbreak{}x_{3}^{6}\,\allowbreak{}x_{4}^{13}\,\allowbreak{}x_{5}^{3}\,\allowbreak{}x_{6}^{10}}
\appmon{227}{x_{1}^{3}\,\allowbreak{}x_{2}\,\allowbreak{}x_{3}^{6}\,\allowbreak{}x_{4}^{3}\,\allowbreak{}x_{5}^{10}\,\allowbreak{}x_{6}^{13}}
\appmon{228}{x_{1}^{3}\,\allowbreak{}x_{2}\,\allowbreak{}x_{3}^{6}\,\allowbreak{}x_{4}^{3}\,\allowbreak{}x_{5}^{11}\,\allowbreak{}x_{6}^{12}}
\appmon{229}{x_{1}^{3}\,\allowbreak{}x_{2}\,\allowbreak{}x_{3}^{6}\,\allowbreak{}x_{4}^{3}\,\allowbreak{}x_{5}^{12}\,\allowbreak{}x_{6}^{11}}
\appmon{230}{x_{1}^{3}\,\allowbreak{}x_{2}\,\allowbreak{}x_{3}^{6}\,\allowbreak{}x_{4}^{3}\,\allowbreak{}x_{5}^{13}\,\allowbreak{}x_{6}^{10}}
\appmon{231}{x_{1}^{3}\,\allowbreak{}x_{2}\,\allowbreak{}x_{3}^{6}\,\allowbreak{}x_{4}^{3}\,\allowbreak{}x_{5}^{14}\,\allowbreak{}x_{6}^{9}}
\appmon{232}{x_{1}^{3}\,\allowbreak{}x_{2}\,\allowbreak{}x_{3}^{6}\,\allowbreak{}x_{4}^{3}\,\allowbreak{}x_{5}^{9}\,\allowbreak{}x_{6}^{14}}
\appmon{233}{x_{1}^{3}\,\allowbreak{}x_{2}\,\allowbreak{}x_{3}^{6}\,\allowbreak{}x_{4}^{5}\,\allowbreak{}x_{5}^{13}\,\allowbreak{}x_{6}^{8}}
\appmon{234}{x_{1}^{3}\,\allowbreak{}x_{2}\,\allowbreak{}x_{3}^{6}\,\allowbreak{}x_{4}^{5}\,\allowbreak{}x_{5}^{8}\,\allowbreak{}x_{6}^{13}}
\appmon{235}{x_{1}^{3}\,\allowbreak{}x_{2}\,\allowbreak{}x_{3}^{6}\,\allowbreak{}x_{4}^{8}\,\allowbreak{}x_{5}^{9}\,\allowbreak{}x_{6}^{9}}
\appmon{236}{x_{1}^{3}\,\allowbreak{}x_{2}\,\allowbreak{}x_{3}^{6}\,\allowbreak{}x_{4}^{9}\,\allowbreak{}x_{5}^{8}\,\allowbreak{}x_{6}^{9}}
\appmon{237}{x_{1}^{3}\,\allowbreak{}x_{2}\,\allowbreak{}x_{3}^{6}\,\allowbreak{}x_{4}^{9}\,\allowbreak{}x_{5}^{9}\,\allowbreak{}x_{6}^{8}}
\appmon{238}{x_{1}^{3}\,\allowbreak{}x_{2}\,\allowbreak{}x_{3}^{7}\,\allowbreak{}x_{4}\,\allowbreak{}x_{5}^{8}\,\allowbreak{}x_{6}^{16}}
\appmon{239}{x_{1}^{3}\,\allowbreak{}x_{2}\,\allowbreak{}x_{3}^{7}\,\allowbreak{}x_{4}^{12}\,\allowbreak{}x_{5}\,\allowbreak{}x_{6}^{12}}
\appmon{240}{x_{1}^{3}\,\allowbreak{}x_{2}\,\allowbreak{}x_{3}^{7}\,\allowbreak{}x_{4}^{5}\,\allowbreak{}x_{5}^{8}\,\allowbreak{}x_{6}^{12}}
\appmon{241}{x_{1}^{3}\,\allowbreak{}x_{2}\,\allowbreak{}x_{3}^{7}\,\allowbreak{}x_{4}^{8}\,\allowbreak{}x_{5}^{5}\,\allowbreak{}x_{6}^{12}}
\appmon{242}{x_{1}^{3}\,\allowbreak{}x_{2}\,\allowbreak{}x_{3}^{7}\,\allowbreak{}x_{4}^{8}\,\allowbreak{}x_{5}^{8}\,\allowbreak{}x_{6}^{9}}
\appmon{243}{x_{1}^{3}\,\allowbreak{}x_{2}\,\allowbreak{}x_{3}^{7}\,\allowbreak{}x_{4}^{8}\,\allowbreak{}x_{5}^{9}\,\allowbreak{}x_{6}^{8}}
\appmon{244}{x_{1}^{3}\,\allowbreak{}x_{2}\,\allowbreak{}x_{3}^{7}\,\allowbreak{}x_{4}^{9}\,\allowbreak{}x_{5}^{4}\,\allowbreak{}x_{6}^{12}}
\appmon{245}{x_{1}^{3}\,\allowbreak{}x_{2}\,\allowbreak{}x_{3}^{7}\,\allowbreak{}x_{4}^{9}\,\allowbreak{}x_{5}^{8}\,\allowbreak{}x_{6}^{8}}
\appmon{246}{x_{1}^{3}\,\allowbreak{}x_{2}^{13}\,\allowbreak{}x_{3}^{3}\,\allowbreak{}x_{4}^{5}\,\allowbreak{}x_{5}^{6}\,\allowbreak{}x_{6}^{6}}
\appmon{247}{x_{1}^{3}\,\allowbreak{}x_{2}^{13}\,\allowbreak{}x_{3}^{3}\,\allowbreak{}x_{4}^{6}\,\allowbreak{}x_{5}^{5}\,\allowbreak{}x_{6}^{6}}
\appmon{248}{x_{1}^{3}\,\allowbreak{}x_{2}^{13}\,\allowbreak{}x_{3}^{3}\,\allowbreak{}x_{4}^{6}\,\allowbreak{}x_{5}^{6}\,\allowbreak{}x_{6}^{5}}
\appmon{249}{x_{1}^{3}\,\allowbreak{}x_{2}^{3}\,\allowbreak{}x_{3}\,\allowbreak{}x_{4}\,\allowbreak{}x_{5}^{14}\,\allowbreak{}x_{6}^{14}}
\appmon{250}{x_{1}^{3}\,\allowbreak{}x_{2}^{3}\,\allowbreak{}x_{3}\,\allowbreak{}x_{4}\,\allowbreak{}x_{5}^{4}\,\allowbreak{}x_{6}^{24}}
\appmon{251}{x_{1}^{3}\,\allowbreak{}x_{2}^{3}\,\allowbreak{}x_{3}\,\allowbreak{}x_{4}^{12}\,\allowbreak{}x_{5}^{13}\,\allowbreak{}x_{6}^{4}}
\appmon{252}{x_{1}^{3}\,\allowbreak{}x_{2}^{3}\,\allowbreak{}x_{3}\,\allowbreak{}x_{4}^{12}\,\allowbreak{}x_{5}^{4}\,\allowbreak{}x_{6}^{13}}
\appmon{253}{x_{1}^{3}\,\allowbreak{}x_{2}^{3}\,\allowbreak{}x_{3}\,\allowbreak{}x_{4}^{13}\,\allowbreak{}x_{5}^{12}\,\allowbreak{}x_{6}^{4}}
\appmon{254}{x_{1}^{3}\,\allowbreak{}x_{2}^{3}\,\allowbreak{}x_{3}\,\allowbreak{}x_{4}^{13}\,\allowbreak{}x_{5}^{6}\,\allowbreak{}x_{6}^{10}}
\appmon{255}{x_{1}^{3}\,\allowbreak{}x_{2}^{3}\,\allowbreak{}x_{3}\,\allowbreak{}x_{4}^{14}\,\allowbreak{}x_{5}^{3}\,\allowbreak{}x_{6}^{12}}
\appmon{256}{x_{1}^{3}\,\allowbreak{}x_{2}^{3}\,\allowbreak{}x_{3}\,\allowbreak{}x_{4}^{14}\,\allowbreak{}x_{5}^{5}\,\allowbreak{}x_{6}^{10}}
\appmon{257}{x_{1}^{3}\,\allowbreak{}x_{2}^{3}\,\allowbreak{}x_{3}\,\allowbreak{}x_{4}^{4}\,\allowbreak{}x_{5}^{11}\,\allowbreak{}x_{6}^{14}}
\appmon{258}{x_{1}^{3}\,\allowbreak{}x_{2}^{3}\,\allowbreak{}x_{3}\,\allowbreak{}x_{4}^{4}\,\allowbreak{}x_{5}^{14}\,\allowbreak{}x_{6}^{11}}
\appmon{259}{x_{1}^{3}\,\allowbreak{}x_{2}^{3}\,\allowbreak{}x_{3}\,\allowbreak{}x_{4}^{5}\,\allowbreak{}x_{5}^{12}\,\allowbreak{}x_{6}^{12}}
\appmon{260}{x_{1}^{3}\,\allowbreak{}x_{2}^{3}\,\allowbreak{}x_{3}\,\allowbreak{}x_{4}^{7}\,\allowbreak{}x_{5}^{10}\,\allowbreak{}x_{6}^{12}}
\appmon{261}{x_{1}^{3}\,\allowbreak{}x_{2}^{3}\,\allowbreak{}x_{3}^{12}\,\allowbreak{}x_{4}\,\allowbreak{}x_{5}^{5}\,\allowbreak{}x_{6}^{12}}
\appmon{262}{x_{1}^{3}\,\allowbreak{}x_{2}^{3}\,\allowbreak{}x_{3}^{13}\,\allowbreak{}x_{4}\,\allowbreak{}x_{5}^{12}\,\allowbreak{}x_{6}^{4}}
\appmon{263}{x_{1}^{3}\,\allowbreak{}x_{2}^{3}\,\allowbreak{}x_{3}^{13}\,\allowbreak{}x_{4}\,\allowbreak{}x_{5}^{4}\,\allowbreak{}x_{6}^{12}}
\appmon{264}{x_{1}^{3}\,\allowbreak{}x_{2}^{3}\,\allowbreak{}x_{3}^{13}\,\allowbreak{}x_{4}^{4}\,\allowbreak{}x_{5}\,\allowbreak{}x_{6}^{12}}
\appmon{265}{x_{1}^{3}\,\allowbreak{}x_{2}^{3}\,\allowbreak{}x_{3}^{13}\,\allowbreak{}x_{4}^{4}\,\allowbreak{}x_{5}^{5}\,\allowbreak{}x_{6}^{8}}
\appmon{266}{x_{1}^{3}\,\allowbreak{}x_{2}^{3}\,\allowbreak{}x_{3}^{13}\,\allowbreak{}x_{4}^{4}\,\allowbreak{}x_{5}^{8}\,\allowbreak{}x_{6}^{5}}
\appmon{267}{x_{1}^{3}\,\allowbreak{}x_{2}^{3}\,\allowbreak{}x_{3}^{13}\,\allowbreak{}x_{4}^{5}\,\allowbreak{}x_{5}^{6}\,\allowbreak{}x_{6}^{6}}
\appmon{268}{x_{1}^{3}\,\allowbreak{}x_{2}^{3}\,\allowbreak{}x_{3}^{13}\,\allowbreak{}x_{4}^{6}\,\allowbreak{}x_{5}^{5}\,\allowbreak{}x_{6}^{6}}
\appmon{269}{x_{1}^{3}\,\allowbreak{}x_{2}^{3}\,\allowbreak{}x_{3}^{13}\,\allowbreak{}x_{4}^{6}\,\allowbreak{}x_{5}^{6}\,\allowbreak{}x_{6}^{5}}
\appmon{270}{x_{1}^{3}\,\allowbreak{}x_{2}^{3}\,\allowbreak{}x_{3}^{3}\,\allowbreak{}x_{4}^{13}\,\allowbreak{}x_{5}^{4}\,\allowbreak{}x_{6}^{10}}
\appmon{271}{x_{1}^{3}\,\allowbreak{}x_{2}^{3}\,\allowbreak{}x_{3}^{3}\,\allowbreak{}x_{4}^{13}\,\allowbreak{}x_{5}^{6}\,\allowbreak{}x_{6}^{8}}
\appmon{272}{x_{1}^{3}\,\allowbreak{}x_{2}^{3}\,\allowbreak{}x_{3}^{3}\,\allowbreak{}x_{4}^{3}\,\allowbreak{}x_{5}^{12}\,\allowbreak{}x_{6}^{12}}
\appmon{273}{x_{1}^{3}\,\allowbreak{}x_{2}^{3}\,\allowbreak{}x_{3}^{3}\,\allowbreak{}x_{4}^{5}\,\allowbreak{}x_{5}^{12}\,\allowbreak{}x_{6}^{10}}
\appmon{274}{x_{1}^{3}\,\allowbreak{}x_{2}^{3}\,\allowbreak{}x_{3}^{4}\,\allowbreak{}x_{4}\,\allowbreak{}x_{5}^{11}\,\allowbreak{}x_{6}^{14}}
\appmon{275}{x_{1}^{3}\,\allowbreak{}x_{2}^{3}\,\allowbreak{}x_{3}^{4}\,\allowbreak{}x_{4}\,\allowbreak{}x_{5}^{12}\,\allowbreak{}x_{6}^{13}}
\appmon{276}{x_{1}^{3}\,\allowbreak{}x_{2}^{3}\,\allowbreak{}x_{3}^{4}\,\allowbreak{}x_{4}\,\allowbreak{}x_{5}^{13}\,\allowbreak{}x_{6}^{12}}
\appmon{277}{x_{1}^{3}\,\allowbreak{}x_{2}^{3}\,\allowbreak{}x_{3}^{4}\,\allowbreak{}x_{4}\,\allowbreak{}x_{5}^{14}\,\allowbreak{}x_{6}^{11}}
\appmon{278}{x_{1}^{3}\,\allowbreak{}x_{2}^{3}\,\allowbreak{}x_{3}^{4}\,\allowbreak{}x_{4}\,\allowbreak{}x_{5}^{9}\,\allowbreak{}x_{6}^{16}}
\appmon{279}{x_{1}^{3}\,\allowbreak{}x_{2}^{3}\,\allowbreak{}x_{3}^{4}\,\allowbreak{}x_{4}^{11}\,\allowbreak{}x_{5}^{3}\,\allowbreak{}x_{6}^{12}}
\appmon{280}{x_{1}^{3}\,\allowbreak{}x_{2}^{3}\,\allowbreak{}x_{3}^{4}\,\allowbreak{}x_{4}^{11}\,\allowbreak{}x_{5}^{5}\,\allowbreak{}x_{6}^{10}}
\appmon{281}{x_{1}^{3}\,\allowbreak{}x_{2}^{3}\,\allowbreak{}x_{3}^{4}\,\allowbreak{}x_{4}^{13}\,\allowbreak{}x_{5}^{3}\,\allowbreak{}x_{6}^{10}}
\appmon{282}{x_{1}^{3}\,\allowbreak{}x_{2}^{3}\,\allowbreak{}x_{3}^{4}\,\allowbreak{}x_{4}^{3}\,\allowbreak{}x_{5}^{11}\,\allowbreak{}x_{6}^{12}}
\appmon{283}{x_{1}^{3}\,\allowbreak{}x_{2}^{3}\,\allowbreak{}x_{3}^{4}\,\allowbreak{}x_{4}^{3}\,\allowbreak{}x_{5}^{12}\,\allowbreak{}x_{6}^{11}}
\appmon{284}{x_{1}^{3}\,\allowbreak{}x_{2}^{3}\,\allowbreak{}x_{3}^{4}\,\allowbreak{}x_{4}^{5}\,\allowbreak{}x_{5}^{10}\,\allowbreak{}x_{6}^{11}}
\appmon{285}{x_{1}^{3}\,\allowbreak{}x_{2}^{3}\,\allowbreak{}x_{3}^{4}\,\allowbreak{}x_{4}^{5}\,\allowbreak{}x_{5}^{11}\,\allowbreak{}x_{6}^{10}}
\appmon{286}{x_{1}^{3}\,\allowbreak{}x_{2}^{3}\,\allowbreak{}x_{3}^{4}\,\allowbreak{}x_{4}^{5}\,\allowbreak{}x_{5}^{14}\,\allowbreak{}x_{6}^{7}}
\appmon{287}{x_{1}^{3}\,\allowbreak{}x_{2}^{3}\,\allowbreak{}x_{3}^{4}\,\allowbreak{}x_{4}^{5}\,\allowbreak{}x_{5}^{7}\,\allowbreak{}x_{6}^{14}}
\appmon{288}{x_{1}^{3}\,\allowbreak{}x_{2}^{3}\,\allowbreak{}x_{3}^{4}\,\allowbreak{}x_{4}^{5}\,\allowbreak{}x_{5}^{9}\,\allowbreak{}x_{6}^{12}}
\appmon{289}{x_{1}^{3}\,\allowbreak{}x_{2}^{3}\,\allowbreak{}x_{3}^{4}\,\allowbreak{}x_{4}^{7}\,\allowbreak{}x_{5}^{9}\,\allowbreak{}x_{6}^{10}}
\appmon{290}{x_{1}^{3}\,\allowbreak{}x_{2}^{3}\,\allowbreak{}x_{3}^{4}\,\allowbreak{}x_{4}^{8}\,\allowbreak{}x_{5}^{13}\,\allowbreak{}x_{6}^{5}}
\appmon{291}{x_{1}^{3}\,\allowbreak{}x_{2}^{3}\,\allowbreak{}x_{3}^{4}\,\allowbreak{}x_{4}^{8}\,\allowbreak{}x_{5}^{5}\,\allowbreak{}x_{6}^{13}}
\appmon{292}{x_{1}^{3}\,\allowbreak{}x_{2}^{3}\,\allowbreak{}x_{3}^{4}\,\allowbreak{}x_{4}^{8}\,\allowbreak{}x_{5}^{9}\,\allowbreak{}x_{6}^{9}}
\appmon{293}{x_{1}^{3}\,\allowbreak{}x_{2}^{3}\,\allowbreak{}x_{3}^{5}\,\allowbreak{}x_{4}\,\allowbreak{}x_{5}^{12}\,\allowbreak{}x_{6}^{12}}
\appmon{294}{x_{1}^{3}\,\allowbreak{}x_{2}^{3}\,\allowbreak{}x_{3}^{5}\,\allowbreak{}x_{4}\,\allowbreak{}x_{5}^{8}\,\allowbreak{}x_{6}^{16}}
\appmon{295}{x_{1}^{3}\,\allowbreak{}x_{2}^{3}\,\allowbreak{}x_{3}^{5}\,\allowbreak{}x_{4}^{10}\,\allowbreak{}x_{5}^{3}\,\allowbreak{}x_{6}^{12}}
\appmon{296}{x_{1}^{3}\,\allowbreak{}x_{2}^{3}\,\allowbreak{}x_{3}^{5}\,\allowbreak{}x_{4}^{10}\,\allowbreak{}x_{5}^{6}\,\allowbreak{}x_{6}^{9}}
\appmon{297}{x_{1}^{3}\,\allowbreak{}x_{2}^{3}\,\allowbreak{}x_{3}^{5}\,\allowbreak{}x_{4}^{11}\,\allowbreak{}x_{5}^{4}\,\allowbreak{}x_{6}^{10}}
\appmon{298}{x_{1}^{3}\,\allowbreak{}x_{2}^{3}\,\allowbreak{}x_{3}^{5}\,\allowbreak{}x_{4}^{11}\,\allowbreak{}x_{5}^{6}\,\allowbreak{}x_{6}^{8}}
\appmon{299}{x_{1}^{3}\,\allowbreak{}x_{2}^{3}\,\allowbreak{}x_{3}^{5}\,\allowbreak{}x_{4}^{13}\,\allowbreak{}x_{5}^{4}\,\allowbreak{}x_{6}^{8}}
\appmon{300}{x_{1}^{3}\,\allowbreak{}x_{2}^{3}\,\allowbreak{}x_{3}^{5}\,\allowbreak{}x_{4}^{14}\,\allowbreak{}x_{5}^{5}\,\allowbreak{}x_{6}^{6}}
\appmon{301}{x_{1}^{3}\,\allowbreak{}x_{2}^{3}\,\allowbreak{}x_{3}^{5}\,\allowbreak{}x_{4}^{14}\,\allowbreak{}x_{5}^{6}\,\allowbreak{}x_{6}^{5}}
\appmon{302}{x_{1}^{3}\,\allowbreak{}x_{2}^{3}\,\allowbreak{}x_{3}^{5}\,\allowbreak{}x_{4}^{3}\,\allowbreak{}x_{5}^{10}\,\allowbreak{}x_{6}^{12}}
\appmon{303}{x_{1}^{3}\,\allowbreak{}x_{2}^{3}\,\allowbreak{}x_{3}^{5}\,\allowbreak{}x_{4}^{3}\,\allowbreak{}x_{5}^{12}\,\allowbreak{}x_{6}^{10}}
\appmon{304}{x_{1}^{3}\,\allowbreak{}x_{2}^{3}\,\allowbreak{}x_{3}^{5}\,\allowbreak{}x_{4}^{4}\,\allowbreak{}x_{5}^{10}\,\allowbreak{}x_{6}^{11}}
\appmon{305}{x_{1}^{3}\,\allowbreak{}x_{2}^{3}\,\allowbreak{}x_{3}^{5}\,\allowbreak{}x_{4}^{4}\,\allowbreak{}x_{5}^{11}\,\allowbreak{}x_{6}^{10}}
\appmon{306}{x_{1}^{3}\,\allowbreak{}x_{2}^{3}\,\allowbreak{}x_{3}^{5}\,\allowbreak{}x_{4}^{4}\,\allowbreak{}x_{5}^{14}\,\allowbreak{}x_{6}^{7}}
\appmon{307}{x_{1}^{3}\,\allowbreak{}x_{2}^{3}\,\allowbreak{}x_{3}^{5}\,\allowbreak{}x_{4}^{4}\,\allowbreak{}x_{5}^{7}\,\allowbreak{}x_{6}^{14}}
\appmon{308}{x_{1}^{3}\,\allowbreak{}x_{2}^{3}\,\allowbreak{}x_{3}^{5}\,\allowbreak{}x_{4}^{4}\,\allowbreak{}x_{5}^{9}\,\allowbreak{}x_{6}^{12}}
\appmon{309}{x_{1}^{3}\,\allowbreak{}x_{2}^{3}\,\allowbreak{}x_{3}^{5}\,\allowbreak{}x_{4}^{5}\,\allowbreak{}x_{5}^{12}\,\allowbreak{}x_{6}^{8}}
\appmon{310}{x_{1}^{3}\,\allowbreak{}x_{2}^{3}\,\allowbreak{}x_{3}^{5}\,\allowbreak{}x_{4}^{5}\,\allowbreak{}x_{5}^{8}\,\allowbreak{}x_{6}^{12}}
\appmon{311}{x_{1}^{3}\,\allowbreak{}x_{2}^{3}\,\allowbreak{}x_{3}^{5}\,\allowbreak{}x_{4}^{6}\,\allowbreak{}x_{5}^{11}\,\allowbreak{}x_{6}^{8}}
\appmon{312}{x_{1}^{3}\,\allowbreak{}x_{2}^{3}\,\allowbreak{}x_{3}^{5}\,\allowbreak{}x_{4}^{6}\,\allowbreak{}x_{5}^{12}\,\allowbreak{}x_{6}^{7}}
\appmon{313}{x_{1}^{3}\,\allowbreak{}x_{2}^{3}\,\allowbreak{}x_{3}^{5}\,\allowbreak{}x_{4}^{6}\,\allowbreak{}x_{5}^{13}\,\allowbreak{}x_{6}^{6}}
\appmon{314}{x_{1}^{3}\,\allowbreak{}x_{2}^{3}\,\allowbreak{}x_{3}^{5}\,\allowbreak{}x_{4}^{6}\,\allowbreak{}x_{5}^{14}\,\allowbreak{}x_{6}^{5}}
\appmon{315}{x_{1}^{3}\,\allowbreak{}x_{2}^{3}\,\allowbreak{}x_{3}^{5}\,\allowbreak{}x_{4}^{6}\,\allowbreak{}x_{5}^{5}\,\allowbreak{}x_{6}^{14}}
\appmon{316}{x_{1}^{3}\,\allowbreak{}x_{2}^{3}\,\allowbreak{}x_{3}^{5}\,\allowbreak{}x_{4}^{6}\,\allowbreak{}x_{5}^{6}\,\allowbreak{}x_{6}^{13}}
\appmon{317}{x_{1}^{3}\,\allowbreak{}x_{2}^{3}\,\allowbreak{}x_{3}^{5}\,\allowbreak{}x_{4}^{6}\,\allowbreak{}x_{5}^{7}\,\allowbreak{}x_{6}^{12}}
\appmon{318}{x_{1}^{3}\,\allowbreak{}x_{2}^{3}\,\allowbreak{}x_{3}^{5}\,\allowbreak{}x_{4}^{6}\,\allowbreak{}x_{5}^{8}\,\allowbreak{}x_{6}^{11}}
\appmon{319}{x_{1}^{3}\,\allowbreak{}x_{2}^{3}\,\allowbreak{}x_{3}^{5}\,\allowbreak{}x_{4}^{7}\,\allowbreak{}x_{5}^{8}\,\allowbreak{}x_{6}^{10}}
\appmon{320}{x_{1}^{3}\,\allowbreak{}x_{2}^{3}\,\allowbreak{}x_{3}^{5}\,\allowbreak{}x_{4}^{8}\,\allowbreak{}x_{5}\,\allowbreak{}x_{6}^{16}}
\appmon{321}{x_{1}^{3}\,\allowbreak{}x_{2}^{3}\,\allowbreak{}x_{3}^{5}\,\allowbreak{}x_{4}^{8}\,\allowbreak{}x_{5}^{8}\,\allowbreak{}x_{6}^{9}}
\appmon{322}{x_{1}^{3}\,\allowbreak{}x_{2}^{3}\,\allowbreak{}x_{3}^{5}\,\allowbreak{}x_{4}^{9}\,\allowbreak{}x_{5}^{4}\,\allowbreak{}x_{6}^{12}}
\appmon{323}{x_{1}^{3}\,\allowbreak{}x_{2}^{3}\,\allowbreak{}x_{3}^{7}\,\allowbreak{}x_{4}^{3}\,\allowbreak{}x_{5}^{12}\,\allowbreak{}x_{6}^{8}}
\appmon{324}{x_{1}^{3}\,\allowbreak{}x_{2}^{3}\,\allowbreak{}x_{3}^{7}\,\allowbreak{}x_{4}^{3}\,\allowbreak{}x_{5}^{8}\,\allowbreak{}x_{6}^{12}}
\appmon{325}{x_{1}^{3}\,\allowbreak{}x_{2}^{3}\,\allowbreak{}x_{3}^{7}\,\allowbreak{}x_{4}^{4}\,\allowbreak{}x_{5}^{9}\,\allowbreak{}x_{6}^{10}}
\appmon{326}{x_{1}^{3}\,\allowbreak{}x_{2}^{3}\,\allowbreak{}x_{3}^{7}\,\allowbreak{}x_{4}^{5}\,\allowbreak{}x_{5}^{10}\,\allowbreak{}x_{6}^{8}}
\appmon{327}{x_{1}^{3}\,\allowbreak{}x_{2}^{3}\,\allowbreak{}x_{3}^{7}\,\allowbreak{}x_{4}^{5}\,\allowbreak{}x_{5}^{12}\,\allowbreak{}x_{6}^{6}}
\appmon{328}{x_{1}^{3}\,\allowbreak{}x_{2}^{3}\,\allowbreak{}x_{3}^{7}\,\allowbreak{}x_{4}^{8}\,\allowbreak{}x_{5}^{3}\,\allowbreak{}x_{6}^{12}}
\appmon{329}{x_{1}^{3}\,\allowbreak{}x_{2}^{3}\,\allowbreak{}x_{3}^{7}\,\allowbreak{}x_{4}^{8}\,\allowbreak{}x_{5}^{5}\,\allowbreak{}x_{6}^{10}}
\appmon{330}{x_{1}^{3}\,\allowbreak{}x_{2}^{4}\,\allowbreak{}x_{3}\,\allowbreak{}x_{4}^{11}\,\allowbreak{}x_{5}^{5}\,\allowbreak{}x_{6}^{12}}
\appmon{331}{x_{1}^{3}\,\allowbreak{}x_{2}^{4}\,\allowbreak{}x_{3}\,\allowbreak{}x_{4}^{2}\,\allowbreak{}x_{5}^{9}\,\allowbreak{}x_{6}^{17}}
\appmon{332}{x_{1}^{3}\,\allowbreak{}x_{2}^{4}\,\allowbreak{}x_{3}\,\allowbreak{}x_{4}^{3}\,\allowbreak{}x_{5}^{12}\,\allowbreak{}x_{6}^{13}}
\appmon{333}{x_{1}^{3}\,\allowbreak{}x_{2}^{4}\,\allowbreak{}x_{3}\,\allowbreak{}x_{4}^{3}\,\allowbreak{}x_{5}^{13}\,\allowbreak{}x_{6}^{12}}
\appmon{334}{x_{1}^{3}\,\allowbreak{}x_{2}^{4}\,\allowbreak{}x_{3}\,\allowbreak{}x_{4}^{3}\,\allowbreak{}x_{5}^{9}\,\allowbreak{}x_{6}^{16}}
\appmon{335}{x_{1}^{3}\,\allowbreak{}x_{2}^{4}\,\allowbreak{}x_{3}^{11}\,\allowbreak{}x_{4}^{5}\,\allowbreak{}x_{5}\,\allowbreak{}x_{6}^{12}}
\appmon{336}{x_{1}^{3}\,\allowbreak{}x_{2}^{4}\,\allowbreak{}x_{3}^{11}\,\allowbreak{}x_{4}^{5}\,\allowbreak{}x_{5}^{3}\,\allowbreak{}x_{6}^{10}}
\appmon{337}{x_{1}^{3}\,\allowbreak{}x_{2}^{4}\,\allowbreak{}x_{3}^{3}\,\allowbreak{}x_{4}^{13}\,\allowbreak{}x_{5}\,\allowbreak{}x_{6}^{12}}
\appmon{338}{x_{1}^{3}\,\allowbreak{}x_{2}^{4}\,\allowbreak{}x_{3}^{3}\,\allowbreak{}x_{4}^{5}\,\allowbreak{}x_{5}^{10}\,\allowbreak{}x_{6}^{11}}
\appmon{339}{x_{1}^{3}\,\allowbreak{}x_{2}^{4}\,\allowbreak{}x_{3}^{3}\,\allowbreak{}x_{4}^{5}\,\allowbreak{}x_{5}^{11}\,\allowbreak{}x_{6}^{10}}
\appmon{340}{x_{1}^{3}\,\allowbreak{}x_{2}^{4}\,\allowbreak{}x_{3}^{3}\,\allowbreak{}x_{4}^{5}\,\allowbreak{}x_{5}^{14}\,\allowbreak{}x_{6}^{7}}
\appmon{341}{x_{1}^{3}\,\allowbreak{}x_{2}^{4}\,\allowbreak{}x_{3}^{3}\,\allowbreak{}x_{4}^{5}\,\allowbreak{}x_{5}^{7}\,\allowbreak{}x_{6}^{14}}
\appmon{342}{x_{1}^{3}\,\allowbreak{}x_{2}^{4}\,\allowbreak{}x_{3}^{3}\,\allowbreak{}x_{4}^{5}\,\allowbreak{}x_{5}^{9}\,\allowbreak{}x_{6}^{12}}
\appmon{343}{x_{1}^{3}\,\allowbreak{}x_{2}^{4}\,\allowbreak{}x_{3}^{3}\,\allowbreak{}x_{4}^{7}\,\allowbreak{}x_{5}^{9}\,\allowbreak{}x_{6}^{10}}
\appmon{344}{x_{1}^{3}\,\allowbreak{}x_{2}^{4}\,\allowbreak{}x_{3}^{3}\,\allowbreak{}x_{4}^{9}\,\allowbreak{}x_{5}^{5}\,\allowbreak{}x_{6}^{12}}
\appmon{345}{x_{1}^{3}\,\allowbreak{}x_{2}^{4}\,\allowbreak{}x_{3}^{7}\,\allowbreak{}x_{4}^{3}\,\allowbreak{}x_{5}^{9}\,\allowbreak{}x_{6}^{10}}
\appmon{346}{x_{1}^{3}\,\allowbreak{}x_{2}^{4}\,\allowbreak{}x_{3}^{7}\,\allowbreak{}x_{4}^{9}\,\allowbreak{}x_{5}^{3}\,\allowbreak{}x_{6}^{10}}
\appmon{347}{x_{1}^{3}\,\allowbreak{}x_{2}^{5}\,\allowbreak{}x_{3}\,\allowbreak{}x_{4}\,\allowbreak{}x_{5}^{2}\,\allowbreak{}x_{6}^{24}}
\appmon{348}{x_{1}^{3}\,\allowbreak{}x_{2}^{5}\,\allowbreak{}x_{3}\,\allowbreak{}x_{4}\,\allowbreak{}x_{5}^{8}\,\allowbreak{}x_{6}^{18}}
\appmon{349}{x_{1}^{3}\,\allowbreak{}x_{2}^{5}\,\allowbreak{}x_{3}\,\allowbreak{}x_{4}^{10}\,\allowbreak{}x_{5}^{5}\,\allowbreak{}x_{6}^{12}}
\appmon{350}{x_{1}^{3}\,\allowbreak{}x_{2}^{5}\,\allowbreak{}x_{3}\,\allowbreak{}x_{4}^{11}\,\allowbreak{}x_{5}^{4}\,\allowbreak{}x_{6}^{12}}
\appmon{351}{x_{1}^{3}\,\allowbreak{}x_{2}^{5}\,\allowbreak{}x_{3}\,\allowbreak{}x_{4}^{2}\,\allowbreak{}x_{5}^{11}\,\allowbreak{}x_{6}^{14}}
\appmon{352}{x_{1}^{3}\,\allowbreak{}x_{2}^{5}\,\allowbreak{}x_{3}\,\allowbreak{}x_{4}^{2}\,\allowbreak{}x_{5}^{14}\,\allowbreak{}x_{6}^{11}}
\appmon{353}{x_{1}^{3}\,\allowbreak{}x_{2}^{5}\,\allowbreak{}x_{3}\,\allowbreak{}x_{4}^{3}\,\allowbreak{}x_{5}^{12}\,\allowbreak{}x_{6}^{12}}
\appmon{354}{x_{1}^{3}\,\allowbreak{}x_{2}^{5}\,\allowbreak{}x_{3}\,\allowbreak{}x_{4}^{6}\,\allowbreak{}x_{5}^{10}\,\allowbreak{}x_{6}^{11}}
\appmon{355}{x_{1}^{3}\,\allowbreak{}x_{2}^{5}\,\allowbreak{}x_{3}\,\allowbreak{}x_{4}^{6}\,\allowbreak{}x_{5}^{11}\,\allowbreak{}x_{6}^{10}}
\appmon{356}{x_{1}^{3}\,\allowbreak{}x_{2}^{5}\,\allowbreak{}x_{3}\,\allowbreak{}x_{4}^{6}\,\allowbreak{}x_{5}^{14}\,\allowbreak{}x_{6}^{7}}
\appmon{357}{x_{1}^{3}\,\allowbreak{}x_{2}^{5}\,\allowbreak{}x_{3}\,\allowbreak{}x_{4}^{6}\,\allowbreak{}x_{5}^{7}\,\allowbreak{}x_{6}^{14}}
\appmon{358}{x_{1}^{3}\,\allowbreak{}x_{2}^{5}\,\allowbreak{}x_{3}\,\allowbreak{}x_{4}^{8}\,\allowbreak{}x_{5}^{11}\,\allowbreak{}x_{6}^{8}}
\appmon{359}{x_{1}^{3}\,\allowbreak{}x_{2}^{5}\,\allowbreak{}x_{3}\,\allowbreak{}x_{4}^{8}\,\allowbreak{}x_{5}^{3}\,\allowbreak{}x_{6}^{16}}
\appmon{360}{x_{1}^{3}\,\allowbreak{}x_{2}^{5}\,\allowbreak{}x_{3}\,\allowbreak{}x_{4}^{8}\,\allowbreak{}x_{5}^{8}\,\allowbreak{}x_{6}^{11}}
\appmon{361}{x_{1}^{3}\,\allowbreak{}x_{2}^{5}\,\allowbreak{}x_{3}\,\allowbreak{}x_{4}^{9}\,\allowbreak{}x_{5}^{2}\,\allowbreak{}x_{6}^{16}}
\appmon{362}{x_{1}^{3}\,\allowbreak{}x_{2}^{5}\,\allowbreak{}x_{3}\,\allowbreak{}x_{4}^{9}\,\allowbreak{}x_{5}^{8}\,\allowbreak{}x_{6}^{10}}
\appmon{363}{x_{1}^{3}\,\allowbreak{}x_{2}^{5}\,\allowbreak{}x_{3}^{11}\,\allowbreak{}x_{4}\,\allowbreak{}x_{5}^{8}\,\allowbreak{}x_{6}^{8}}
\appmon{364}{x_{1}^{3}\,\allowbreak{}x_{2}^{5}\,\allowbreak{}x_{3}^{2}\,\allowbreak{}x_{4}\,\allowbreak{}x_{5}^{11}\,\allowbreak{}x_{6}^{14}}
\appmon{365}{x_{1}^{3}\,\allowbreak{}x_{2}^{5}\,\allowbreak{}x_{3}^{2}\,\allowbreak{}x_{4}\,\allowbreak{}x_{5}^{14}\,\allowbreak{}x_{6}^{11}}
\appmon{366}{x_{1}^{3}\,\allowbreak{}x_{2}^{5}\,\allowbreak{}x_{3}^{2}\,\allowbreak{}x_{4}\,\allowbreak{}x_{5}^{8}\,\allowbreak{}x_{6}^{17}}
\appmon{367}{x_{1}^{3}\,\allowbreak{}x_{2}^{5}\,\allowbreak{}x_{3}^{2}\,\allowbreak{}x_{4}\,\allowbreak{}x_{5}^{9}\,\allowbreak{}x_{6}^{16}}
\appmon{368}{x_{1}^{3}\,\allowbreak{}x_{2}^{5}\,\allowbreak{}x_{3}^{2}\,\allowbreak{}x_{4}^{11}\,\allowbreak{}x_{5}^{3}\,\allowbreak{}x_{6}^{12}}
\appmon{369}{x_{1}^{3}\,\allowbreak{}x_{2}^{5}\,\allowbreak{}x_{3}^{2}\,\allowbreak{}x_{4}^{11}\,\allowbreak{}x_{5}^{5}\,\allowbreak{}x_{6}^{10}}
\appmon{370}{x_{1}^{3}\,\allowbreak{}x_{2}^{5}\,\allowbreak{}x_{3}^{2}\,\allowbreak{}x_{4}^{13}\,\allowbreak{}x_{5}\,\allowbreak{}x_{6}^{12}}
\appmon{371}{x_{1}^{3}\,\allowbreak{}x_{2}^{5}\,\allowbreak{}x_{3}^{2}\,\allowbreak{}x_{4}^{13}\,\allowbreak{}x_{5}^{3}\,\allowbreak{}x_{6}^{10}}
\appmon{372}{x_{1}^{3}\,\allowbreak{}x_{2}^{5}\,\allowbreak{}x_{3}^{2}\,\allowbreak{}x_{4}^{3}\,\allowbreak{}x_{5}^{10}\,\allowbreak{}x_{6}^{13}}
\appmon{373}{x_{1}^{3}\,\allowbreak{}x_{2}^{5}\,\allowbreak{}x_{3}^{2}\,\allowbreak{}x_{4}^{3}\,\allowbreak{}x_{5}^{11}\,\allowbreak{}x_{6}^{12}}
\appmon{374}{x_{1}^{3}\,\allowbreak{}x_{2}^{5}\,\allowbreak{}x_{3}^{2}\,\allowbreak{}x_{4}^{3}\,\allowbreak{}x_{5}^{12}\,\allowbreak{}x_{6}^{11}}
\appmon{375}{x_{1}^{3}\,\allowbreak{}x_{2}^{5}\,\allowbreak{}x_{3}^{2}\,\allowbreak{}x_{4}^{3}\,\allowbreak{}x_{5}^{13}\,\allowbreak{}x_{6}^{10}}
\appmon{376}{x_{1}^{3}\,\allowbreak{}x_{2}^{5}\,\allowbreak{}x_{3}^{2}\,\allowbreak{}x_{4}^{3}\,\allowbreak{}x_{5}^{14}\,\allowbreak{}x_{6}^{9}}
\appmon{377}{x_{1}^{3}\,\allowbreak{}x_{2}^{5}\,\allowbreak{}x_{3}^{2}\,\allowbreak{}x_{4}^{3}\,\allowbreak{}x_{5}^{9}\,\allowbreak{}x_{6}^{14}}
\appmon{378}{x_{1}^{3}\,\allowbreak{}x_{2}^{5}\,\allowbreak{}x_{3}^{2}\,\allowbreak{}x_{4}^{5}\,\allowbreak{}x_{5}^{13}\,\allowbreak{}x_{6}^{8}}
\appmon{379}{x_{1}^{3}\,\allowbreak{}x_{2}^{5}\,\allowbreak{}x_{3}^{2}\,\allowbreak{}x_{4}^{5}\,\allowbreak{}x_{5}^{8}\,\allowbreak{}x_{6}^{13}}
\appmon{380}{x_{1}^{3}\,\allowbreak{}x_{2}^{5}\,\allowbreak{}x_{3}^{2}\,\allowbreak{}x_{4}^{8}\,\allowbreak{}x_{5}^{9}\,\allowbreak{}x_{6}^{9}}
\appmon{381}{x_{1}^{3}\,\allowbreak{}x_{2}^{5}\,\allowbreak{}x_{3}^{2}\,\allowbreak{}x_{4}^{9}\,\allowbreak{}x_{5}\,\allowbreak{}x_{6}^{16}}
\appmon{382}{x_{1}^{3}\,\allowbreak{}x_{2}^{5}\,\allowbreak{}x_{3}^{2}\,\allowbreak{}x_{4}^{9}\,\allowbreak{}x_{5}^{8}\,\allowbreak{}x_{6}^{9}}
\appmon{383}{x_{1}^{3}\,\allowbreak{}x_{2}^{5}\,\allowbreak{}x_{3}^{3}\,\allowbreak{}x_{4}\,\allowbreak{}x_{5}^{12}\,\allowbreak{}x_{6}^{12}}
\appmon{384}{x_{1}^{3}\,\allowbreak{}x_{2}^{5}\,\allowbreak{}x_{3}^{3}\,\allowbreak{}x_{4}^{10}\,\allowbreak{}x_{5}^{5}\,\allowbreak{}x_{6}^{10}}
\appmon{385}{x_{1}^{3}\,\allowbreak{}x_{2}^{5}\,\allowbreak{}x_{3}^{3}\,\allowbreak{}x_{4}^{10}\,\allowbreak{}x_{5}^{6}\,\allowbreak{}x_{6}^{9}}
\appmon{386}{x_{1}^{3}\,\allowbreak{}x_{2}^{5}\,\allowbreak{}x_{3}^{3}\,\allowbreak{}x_{4}^{12}\,\allowbreak{}x_{5}^{8}\,\allowbreak{}x_{6}^{5}}
\appmon{387}{x_{1}^{3}\,\allowbreak{}x_{2}^{5}\,\allowbreak{}x_{3}^{3}\,\allowbreak{}x_{4}^{13}\,\allowbreak{}x_{5}^{4}\,\allowbreak{}x_{6}^{8}}
\appmon{388}{x_{1}^{3}\,\allowbreak{}x_{2}^{5}\,\allowbreak{}x_{3}^{3}\,\allowbreak{}x_{4}^{13}\,\allowbreak{}x_{5}^{6}\,\allowbreak{}x_{6}^{6}}
\appmon{389}{x_{1}^{3}\,\allowbreak{}x_{2}^{5}\,\allowbreak{}x_{3}^{3}\,\allowbreak{}x_{4}^{14}\,\allowbreak{}x_{5}^{3}\,\allowbreak{}x_{6}^{8}}
\appmon{390}{x_{1}^{3}\,\allowbreak{}x_{2}^{5}\,\allowbreak{}x_{3}^{3}\,\allowbreak{}x_{4}^{14}\,\allowbreak{}x_{5}^{6}\,\allowbreak{}x_{6}^{5}}
\appmon{391}{x_{1}^{3}\,\allowbreak{}x_{2}^{5}\,\allowbreak{}x_{3}^{3}\,\allowbreak{}x_{4}^{2}\,\allowbreak{}x_{5}^{10}\,\allowbreak{}x_{6}^{13}}
\appmon{392}{x_{1}^{3}\,\allowbreak{}x_{2}^{5}\,\allowbreak{}x_{3}^{3}\,\allowbreak{}x_{4}^{2}\,\allowbreak{}x_{5}^{11}\,\allowbreak{}x_{6}^{12}}
\appmon{393}{x_{1}^{3}\,\allowbreak{}x_{2}^{5}\,\allowbreak{}x_{3}^{3}\,\allowbreak{}x_{4}^{2}\,\allowbreak{}x_{5}^{12}\,\allowbreak{}x_{6}^{11}}
\appmon{394}{x_{1}^{3}\,\allowbreak{}x_{2}^{5}\,\allowbreak{}x_{3}^{3}\,\allowbreak{}x_{4}^{2}\,\allowbreak{}x_{5}^{13}\,\allowbreak{}x_{6}^{10}}
\appmon{395}{x_{1}^{3}\,\allowbreak{}x_{2}^{5}\,\allowbreak{}x_{3}^{3}\,\allowbreak{}x_{4}^{2}\,\allowbreak{}x_{5}^{14}\,\allowbreak{}x_{6}^{9}}
\appmon{396}{x_{1}^{3}\,\allowbreak{}x_{2}^{5}\,\allowbreak{}x_{3}^{3}\,\allowbreak{}x_{4}^{2}\,\allowbreak{}x_{5}^{9}\,\allowbreak{}x_{6}^{14}}
\appmon{397}{x_{1}^{3}\,\allowbreak{}x_{2}^{5}\,\allowbreak{}x_{3}^{3}\,\allowbreak{}x_{4}^{4}\,\allowbreak{}x_{5}^{13}\,\allowbreak{}x_{6}^{8}}
\appmon{398}{x_{1}^{3}\,\allowbreak{}x_{2}^{5}\,\allowbreak{}x_{3}^{3}\,\allowbreak{}x_{4}^{4}\,\allowbreak{}x_{5}^{8}\,\allowbreak{}x_{6}^{13}}
\appmon{399}{x_{1}^{3}\,\allowbreak{}x_{2}^{5}\,\allowbreak{}x_{3}^{3}\,\allowbreak{}x_{4}^{5}\,\allowbreak{}x_{5}^{10}\,\allowbreak{}x_{6}^{10}}
\appmon{400}{x_{1}^{3}\,\allowbreak{}x_{2}^{5}\,\allowbreak{}x_{3}^{3}\,\allowbreak{}x_{4}^{5}\,\allowbreak{}x_{5}^{12}\,\allowbreak{}x_{6}^{8}}
\appmon{401}{x_{1}^{3}\,\allowbreak{}x_{2}^{5}\,\allowbreak{}x_{3}^{3}\,\allowbreak{}x_{4}^{5}\,\allowbreak{}x_{5}^{8}\,\allowbreak{}x_{6}^{12}}
\appmon{402}{x_{1}^{3}\,\allowbreak{}x_{2}^{5}\,\allowbreak{}x_{3}^{3}\,\allowbreak{}x_{4}^{6}\,\allowbreak{}x_{5}^{11}\,\allowbreak{}x_{6}^{8}}
\appmon{403}{x_{1}^{3}\,\allowbreak{}x_{2}^{5}\,\allowbreak{}x_{3}^{3}\,\allowbreak{}x_{4}^{6}\,\allowbreak{}x_{5}^{12}\,\allowbreak{}x_{6}^{7}}
\appmon{404}{x_{1}^{3}\,\allowbreak{}x_{2}^{5}\,\allowbreak{}x_{3}^{3}\,\allowbreak{}x_{4}^{6}\,\allowbreak{}x_{5}^{13}\,\allowbreak{}x_{6}^{6}}
\appmon{405}{x_{1}^{3}\,\allowbreak{}x_{2}^{5}\,\allowbreak{}x_{3}^{3}\,\allowbreak{}x_{4}^{6}\,\allowbreak{}x_{5}^{14}\,\allowbreak{}x_{6}^{5}}
\appmon{406}{x_{1}^{3}\,\allowbreak{}x_{2}^{5}\,\allowbreak{}x_{3}^{3}\,\allowbreak{}x_{4}^{6}\,\allowbreak{}x_{5}^{5}\,\allowbreak{}x_{6}^{14}}
\appmon{407}{x_{1}^{3}\,\allowbreak{}x_{2}^{5}\,\allowbreak{}x_{3}^{3}\,\allowbreak{}x_{4}^{6}\,\allowbreak{}x_{5}^{6}\,\allowbreak{}x_{6}^{13}}
\appmon{408}{x_{1}^{3}\,\allowbreak{}x_{2}^{5}\,\allowbreak{}x_{3}^{3}\,\allowbreak{}x_{4}^{6}\,\allowbreak{}x_{5}^{7}\,\allowbreak{}x_{6}^{12}}
\appmon{409}{x_{1}^{3}\,\allowbreak{}x_{2}^{5}\,\allowbreak{}x_{3}^{3}\,\allowbreak{}x_{4}^{6}\,\allowbreak{}x_{5}^{8}\,\allowbreak{}x_{6}^{11}}
\appmon{410}{x_{1}^{3}\,\allowbreak{}x_{2}^{5}\,\allowbreak{}x_{3}^{3}\,\allowbreak{}x_{4}^{6}\,\allowbreak{}x_{5}^{9}\,\allowbreak{}x_{6}^{10}}
\appmon{411}{x_{1}^{3}\,\allowbreak{}x_{2}^{5}\,\allowbreak{}x_{3}^{3}\,\allowbreak{}x_{4}^{7}\,\allowbreak{}x_{5}^{6}\,\allowbreak{}x_{6}^{12}}
\appmon{412}{x_{1}^{3}\,\allowbreak{}x_{2}^{5}\,\allowbreak{}x_{3}^{3}\,\allowbreak{}x_{4}^{7}\,\allowbreak{}x_{5}^{8}\,\allowbreak{}x_{6}^{10}}
\appmon{413}{x_{1}^{3}\,\allowbreak{}x_{2}^{5}\,\allowbreak{}x_{3}^{3}\,\allowbreak{}x_{4}^{8}\,\allowbreak{}x_{5}^{12}\,\allowbreak{}x_{6}^{5}}
\appmon{414}{x_{1}^{3}\,\allowbreak{}x_{2}^{5}\,\allowbreak{}x_{3}^{3}\,\allowbreak{}x_{4}^{8}\,\allowbreak{}x_{5}^{8}\,\allowbreak{}x_{6}^{9}}
\appmon{415}{x_{1}^{3}\,\allowbreak{}x_{2}^{5}\,\allowbreak{}x_{3}^{3}\,\allowbreak{}x_{4}^{8}\,\allowbreak{}x_{5}^{9}\,\allowbreak{}x_{6}^{8}}
\appmon{416}{x_{1}^{3}\,\allowbreak{}x_{2}^{5}\,\allowbreak{}x_{3}^{3}\,\allowbreak{}x_{4}^{9}\,\allowbreak{}x_{5}^{4}\,\allowbreak{}x_{6}^{12}}
\appmon{417}{x_{1}^{3}\,\allowbreak{}x_{2}^{5}\,\allowbreak{}x_{3}^{3}\,\allowbreak{}x_{4}^{9}\,\allowbreak{}x_{5}^{6}\,\allowbreak{}x_{6}^{10}}
\appmon{418}{x_{1}^{3}\,\allowbreak{}x_{2}^{5}\,\allowbreak{}x_{3}^{3}\,\allowbreak{}x_{4}^{9}\,\allowbreak{}x_{5}^{8}\,\allowbreak{}x_{6}^{8}}
\appmon{419}{x_{1}^{3}\,\allowbreak{}x_{2}^{5}\,\allowbreak{}x_{3}^{6}\,\allowbreak{}x_{4}\,\allowbreak{}x_{5}^{10}\,\allowbreak{}x_{6}^{11}}
\appmon{420}{x_{1}^{3}\,\allowbreak{}x_{2}^{5}\,\allowbreak{}x_{3}^{6}\,\allowbreak{}x_{4}\,\allowbreak{}x_{5}^{11}\,\allowbreak{}x_{6}^{10}}
\appmon{421}{x_{1}^{3}\,\allowbreak{}x_{2}^{5}\,\allowbreak{}x_{3}^{6}\,\allowbreak{}x_{4}\,\allowbreak{}x_{5}^{14}\,\allowbreak{}x_{6}^{7}}
\appmon{422}{x_{1}^{3}\,\allowbreak{}x_{2}^{5}\,\allowbreak{}x_{3}^{6}\,\allowbreak{}x_{4}\,\allowbreak{}x_{5}^{7}\,\allowbreak{}x_{6}^{14}}
\appmon{423}{x_{1}^{3}\,\allowbreak{}x_{2}^{5}\,\allowbreak{}x_{3}^{6}\,\allowbreak{}x_{4}^{11}\,\allowbreak{}x_{5}^{3}\,\allowbreak{}x_{6}^{8}}
\appmon{424}{x_{1}^{3}\,\allowbreak{}x_{2}^{5}\,\allowbreak{}x_{3}^{6}\,\allowbreak{}x_{4}^{11}\,\allowbreak{}x_{5}^{5}\,\allowbreak{}x_{6}^{6}}
\appmon{425}{x_{1}^{3}\,\allowbreak{}x_{2}^{5}\,\allowbreak{}x_{3}^{6}\,\allowbreak{}x_{4}^{14}\,\allowbreak{}x_{5}^{3}\,\allowbreak{}x_{6}^{5}}
\appmon{426}{x_{1}^{3}\,\allowbreak{}x_{2}^{5}\,\allowbreak{}x_{3}^{6}\,\allowbreak{}x_{4}^{3}\,\allowbreak{}x_{5}^{11}\,\allowbreak{}x_{6}^{8}}
\appmon{427}{x_{1}^{3}\,\allowbreak{}x_{2}^{5}\,\allowbreak{}x_{3}^{6}\,\allowbreak{}x_{4}^{3}\,\allowbreak{}x_{5}^{12}\,\allowbreak{}x_{6}^{7}}
\appmon{428}{x_{1}^{3}\,\allowbreak{}x_{2}^{5}\,\allowbreak{}x_{3}^{6}\,\allowbreak{}x_{4}^{3}\,\allowbreak{}x_{5}^{13}\,\allowbreak{}x_{6}^{6}}
\appmon{429}{x_{1}^{3}\,\allowbreak{}x_{2}^{5}\,\allowbreak{}x_{3}^{6}\,\allowbreak{}x_{4}^{3}\,\allowbreak{}x_{5}^{14}\,\allowbreak{}x_{6}^{5}}
\appmon{430}{x_{1}^{3}\,\allowbreak{}x_{2}^{5}\,\allowbreak{}x_{3}^{6}\,\allowbreak{}x_{4}^{3}\,\allowbreak{}x_{5}^{5}\,\allowbreak{}x_{6}^{14}}
\appmon{431}{x_{1}^{3}\,\allowbreak{}x_{2}^{5}\,\allowbreak{}x_{3}^{6}\,\allowbreak{}x_{4}^{3}\,\allowbreak{}x_{5}^{6}\,\allowbreak{}x_{6}^{13}}
\appmon{432}{x_{1}^{3}\,\allowbreak{}x_{2}^{5}\,\allowbreak{}x_{3}^{6}\,\allowbreak{}x_{4}^{3}\,\allowbreak{}x_{5}^{7}\,\allowbreak{}x_{6}^{12}}
\appmon{433}{x_{1}^{3}\,\allowbreak{}x_{2}^{5}\,\allowbreak{}x_{3}^{6}\,\allowbreak{}x_{4}^{3}\,\allowbreak{}x_{5}^{8}\,\allowbreak{}x_{6}^{11}}
\appmon{434}{x_{1}^{3}\,\allowbreak{}x_{2}^{5}\,\allowbreak{}x_{3}^{6}\,\allowbreak{}x_{4}^{7}\,\allowbreak{}x_{5}^{3}\,\allowbreak{}x_{6}^{12}}
\appmon{435}{x_{1}^{3}\,\allowbreak{}x_{2}^{5}\,\allowbreak{}x_{3}^{6}\,\allowbreak{}x_{4}^{7}\,\allowbreak{}x_{5}^{6}\,\allowbreak{}x_{6}^{9}}
\appmon{436}{x_{1}^{3}\,\allowbreak{}x_{2}^{5}\,\allowbreak{}x_{3}^{7}\,\allowbreak{}x_{4}^{10}\,\allowbreak{}x_{5}^{5}\,\allowbreak{}x_{6}^{6}}
\appmon{437}{x_{1}^{3}\,\allowbreak{}x_{2}^{5}\,\allowbreak{}x_{3}^{7}\,\allowbreak{}x_{4}^{10}\,\allowbreak{}x_{5}^{6}\,\allowbreak{}x_{6}^{5}}
\appmon{438}{x_{1}^{3}\,\allowbreak{}x_{2}^{5}\,\allowbreak{}x_{3}^{7}\,\allowbreak{}x_{4}^{3}\,\allowbreak{}x_{5}^{10}\,\allowbreak{}x_{6}^{8}}
\appmon{439}{x_{1}^{3}\,\allowbreak{}x_{2}^{5}\,\allowbreak{}x_{3}^{7}\,\allowbreak{}x_{4}^{3}\,\allowbreak{}x_{5}^{12}\,\allowbreak{}x_{6}^{6}}
\appmon{440}{x_{1}^{3}\,\allowbreak{}x_{2}^{5}\,\allowbreak{}x_{3}^{7}\,\allowbreak{}x_{4}^{3}\,\allowbreak{}x_{5}^{6}\,\allowbreak{}x_{6}^{12}}
\appmon{441}{x_{1}^{3}\,\allowbreak{}x_{2}^{5}\,\allowbreak{}x_{3}^{7}\,\allowbreak{}x_{4}^{3}\,\allowbreak{}x_{5}^{8}\,\allowbreak{}x_{6}^{10}}
\appmon{442}{x_{1}^{3}\,\allowbreak{}x_{2}^{5}\,\allowbreak{}x_{3}^{7}\,\allowbreak{}x_{4}^{6}\,\allowbreak{}x_{5}^{10}\,\allowbreak{}x_{6}^{5}}
\appmon{443}{x_{1}^{3}\,\allowbreak{}x_{2}^{5}\,\allowbreak{}x_{3}^{7}\,\allowbreak{}x_{4}^{6}\,\allowbreak{}x_{5}^{3}\,\allowbreak{}x_{6}^{12}}
\appmon{444}{x_{1}^{3}\,\allowbreak{}x_{2}^{5}\,\allowbreak{}x_{3}^{7}\,\allowbreak{}x_{4}^{6}\,\allowbreak{}x_{5}^{6}\,\allowbreak{}x_{6}^{9}}
\appmon{445}{x_{1}^{3}\,\allowbreak{}x_{2}^{5}\,\allowbreak{}x_{3}^{7}\,\allowbreak{}x_{4}^{6}\,\allowbreak{}x_{5}^{9}\,\allowbreak{}x_{6}^{6}}
\appmon{446}{x_{1}^{3}\,\allowbreak{}x_{2}^{5}\,\allowbreak{}x_{3}^{7}\,\allowbreak{}x_{4}^{8}\,\allowbreak{}x_{5}^{3}\,\allowbreak{}x_{6}^{10}}
\appmon{447}{x_{1}^{3}\,\allowbreak{}x_{2}^{5}\,\allowbreak{}x_{3}^{7}\,\allowbreak{}x_{4}^{9}\,\allowbreak{}x_{5}^{6}\,\allowbreak{}x_{6}^{6}}
\appmon{448}{x_{1}^{3}\,\allowbreak{}x_{2}^{5}\,\allowbreak{}x_{3}^{9}\,\allowbreak{}x_{4}\,\allowbreak{}x_{5}^{10}\,\allowbreak{}x_{6}^{8}}
\appmon{449}{x_{1}^{3}\,\allowbreak{}x_{2}^{5}\,\allowbreak{}x_{3}^{9}\,\allowbreak{}x_{4}\,\allowbreak{}x_{5}^{2}\,\allowbreak{}x_{6}^{16}}
\appmon{450}{x_{1}^{3}\,\allowbreak{}x_{2}^{5}\,\allowbreak{}x_{3}^{9}\,\allowbreak{}x_{4}^{9}\,\allowbreak{}x_{5}^{2}\,\allowbreak{}x_{6}^{8}}
\appmon{451}{x_{1}^{3}\,\allowbreak{}x_{2}^{7}\,\allowbreak{}x_{3}\,\allowbreak{}x_{4}\,\allowbreak{}x_{5}^{8}\,\allowbreak{}x_{6}^{16}}
\appmon{452}{x_{1}^{3}\,\allowbreak{}x_{2}^{7}\,\allowbreak{}x_{3}\,\allowbreak{}x_{4}^{8}\,\allowbreak{}x_{5}^{8}\,\allowbreak{}x_{6}^{9}}
\appmon{453}{x_{1}^{3}\,\allowbreak{}x_{2}^{7}\,\allowbreak{}x_{3}\,\allowbreak{}x_{4}^{8}\,\allowbreak{}x_{5}^{9}\,\allowbreak{}x_{6}^{8}}
\appmon{454}{x_{1}^{3}\,\allowbreak{}x_{2}^{7}\,\allowbreak{}x_{3}\,\allowbreak{}x_{4}^{9}\,\allowbreak{}x_{5}^{8}\,\allowbreak{}x_{6}^{8}}
\appmon{455}{x_{1}^{3}\,\allowbreak{}x_{2}^{7}\,\allowbreak{}x_{3}^{3}\,\allowbreak{}x_{4}^{3}\,\allowbreak{}x_{5}^{12}\,\allowbreak{}x_{6}^{8}}
\appmon{456}{x_{1}^{3}\,\allowbreak{}x_{2}^{7}\,\allowbreak{}x_{3}^{3}\,\allowbreak{}x_{4}^{3}\,\allowbreak{}x_{5}^{8}\,\allowbreak{}x_{6}^{12}}
\appmon{457}{x_{1}^{3}\,\allowbreak{}x_{2}^{7}\,\allowbreak{}x_{3}^{3}\,\allowbreak{}x_{4}^{4}\,\allowbreak{}x_{5}^{9}\,\allowbreak{}x_{6}^{10}}
\appmon{458}{x_{1}^{3}\,\allowbreak{}x_{2}^{7}\,\allowbreak{}x_{3}^{3}\,\allowbreak{}x_{4}^{5}\,\allowbreak{}x_{5}^{10}\,\allowbreak{}x_{6}^{8}}
\appmon{459}{x_{1}^{3}\,\allowbreak{}x_{2}^{7}\,\allowbreak{}x_{3}^{3}\,\allowbreak{}x_{4}^{5}\,\allowbreak{}x_{5}^{12}\,\allowbreak{}x_{6}^{6}}
\appmon{460}{x_{1}^{3}\,\allowbreak{}x_{2}^{7}\,\allowbreak{}x_{3}^{3}\,\allowbreak{}x_{4}^{8}\,\allowbreak{}x_{5}^{3}\,\allowbreak{}x_{6}^{12}}
\appmon{461}{x_{1}^{3}\,\allowbreak{}x_{2}^{7}\,\allowbreak{}x_{3}^{3}\,\allowbreak{}x_{4}^{8}\,\allowbreak{}x_{5}^{5}\,\allowbreak{}x_{6}^{10}}
\appmon{462}{x_{1}^{3}\,\allowbreak{}x_{2}^{7}\,\allowbreak{}x_{3}^{4}\,\allowbreak{}x_{4}^{3}\,\allowbreak{}x_{5}^{9}\,\allowbreak{}x_{6}^{10}}
\appmon{463}{x_{1}^{3}\,\allowbreak{}x_{2}^{7}\,\allowbreak{}x_{3}^{5}\,\allowbreak{}x_{4}^{10}\,\allowbreak{}x_{5}^{5}\,\allowbreak{}x_{6}^{6}}
\appmon{464}{x_{1}^{3}\,\allowbreak{}x_{2}^{7}\,\allowbreak{}x_{3}^{5}\,\allowbreak{}x_{4}^{10}\,\allowbreak{}x_{5}^{6}\,\allowbreak{}x_{6}^{5}}
\appmon{465}{x_{1}^{3}\,\allowbreak{}x_{2}^{7}\,\allowbreak{}x_{3}^{5}\,\allowbreak{}x_{4}^{3}\,\allowbreak{}x_{5}^{10}\,\allowbreak{}x_{6}^{8}}
\appmon{466}{x_{1}^{3}\,\allowbreak{}x_{2}^{7}\,\allowbreak{}x_{3}^{5}\,\allowbreak{}x_{4}^{3}\,\allowbreak{}x_{5}^{12}\,\allowbreak{}x_{6}^{6}}
\appmon{467}{x_{1}^{3}\,\allowbreak{}x_{2}^{7}\,\allowbreak{}x_{3}^{5}\,\allowbreak{}x_{4}^{3}\,\allowbreak{}x_{5}^{6}\,\allowbreak{}x_{6}^{12}}
\appmon{468}{x_{1}^{3}\,\allowbreak{}x_{2}^{7}\,\allowbreak{}x_{3}^{5}\,\allowbreak{}x_{4}^{3}\,\allowbreak{}x_{5}^{8}\,\allowbreak{}x_{6}^{10}}
\appmon{469}{x_{1}^{3}\,\allowbreak{}x_{2}^{7}\,\allowbreak{}x_{3}^{5}\,\allowbreak{}x_{4}^{5}\,\allowbreak{}x_{5}^{8}\,\allowbreak{}x_{6}^{8}}
\appmon{470}{x_{1}^{3}\,\allowbreak{}x_{2}^{7}\,\allowbreak{}x_{3}^{5}\,\allowbreak{}x_{4}^{6}\,\allowbreak{}x_{5}^{10}\,\allowbreak{}x_{6}^{5}}
\appmon{471}{x_{1}^{3}\,\allowbreak{}x_{2}^{7}\,\allowbreak{}x_{3}^{5}\,\allowbreak{}x_{4}^{6}\,\allowbreak{}x_{5}^{3}\,\allowbreak{}x_{6}^{12}}
\appmon{472}{x_{1}^{3}\,\allowbreak{}x_{2}^{7}\,\allowbreak{}x_{3}^{5}\,\allowbreak{}x_{4}^{6}\,\allowbreak{}x_{5}^{6}\,\allowbreak{}x_{6}^{9}}
\appmon{473}{x_{1}^{3}\,\allowbreak{}x_{2}^{7}\,\allowbreak{}x_{3}^{5}\,\allowbreak{}x_{4}^{6}\,\allowbreak{}x_{5}^{9}\,\allowbreak{}x_{6}^{6}}
\appmon{474}{x_{1}^{3}\,\allowbreak{}x_{2}^{7}\,\allowbreak{}x_{3}^{5}\,\allowbreak{}x_{4}^{8}\,\allowbreak{}x_{5}\,\allowbreak{}x_{6}^{12}}
\appmon{475}{x_{1}^{3}\,\allowbreak{}x_{2}^{7}\,\allowbreak{}x_{3}^{5}\,\allowbreak{}x_{4}^{8}\,\allowbreak{}x_{5}^{3}\,\allowbreak{}x_{6}^{10}}
\appmon{476}{x_{1}^{3}\,\allowbreak{}x_{2}^{7}\,\allowbreak{}x_{3}^{5}\,\allowbreak{}x_{4}^{9}\,\allowbreak{}x_{5}^{6}\,\allowbreak{}x_{6}^{6}}
\appmon{477}{x_{1}^{3}\,\allowbreak{}x_{2}^{7}\,\allowbreak{}x_{3}^{8}\,\allowbreak{}x_{4}^{5}\,\allowbreak{}x_{5}\,\allowbreak{}x_{6}^{12}}
\appmon{478}{x_{1}^{3}\,\allowbreak{}x_{2}^{7}\,\allowbreak{}x_{3}^{8}\,\allowbreak{}x_{4}^{5}\,\allowbreak{}x_{5}^{3}\,\allowbreak{}x_{6}^{10}}
\appmon{479}{x_{1}^{3}\,\allowbreak{}x_{2}^{7}\,\allowbreak{}x_{3}^{9}\,\allowbreak{}x_{4}\,\allowbreak{}x_{5}^{8}\,\allowbreak{}x_{6}^{8}}
\appmon{480}{x_{1}^{7}\,\allowbreak{}x_{2}\,\allowbreak{}x_{3}\,\allowbreak{}x_{4}\,\allowbreak{}x_{5}^{2}\,\allowbreak{}x_{6}^{24}}
\appmon{481}{x_{1}^{7}\,\allowbreak{}x_{2}\,\allowbreak{}x_{3}\,\allowbreak{}x_{4}^{10}\,\allowbreak{}x_{5}^{12}\,\allowbreak{}x_{6}^{5}}
\appmon{482}{x_{1}^{7}\,\allowbreak{}x_{2}\,\allowbreak{}x_{3}\,\allowbreak{}x_{4}^{10}\,\allowbreak{}x_{5}^{5}\,\allowbreak{}x_{6}^{12}}
\appmon{483}{x_{1}^{7}\,\allowbreak{}x_{2}\,\allowbreak{}x_{3}\,\allowbreak{}x_{4}^{3}\,\allowbreak{}x_{5}^{12}\,\allowbreak{}x_{6}^{12}}
\appmon{484}{x_{1}^{7}\,\allowbreak{}x_{2}\,\allowbreak{}x_{3}^{10}\,\allowbreak{}x_{4}^{5}\,\allowbreak{}x_{5}\,\allowbreak{}x_{6}^{12}}
\appmon{485}{x_{1}^{7}\,\allowbreak{}x_{2}\,\allowbreak{}x_{3}^{3}\,\allowbreak{}x_{4}\,\allowbreak{}x_{5}^{8}\,\allowbreak{}x_{6}^{16}}
\appmon{486}{x_{1}^{7}\,\allowbreak{}x_{2}\,\allowbreak{}x_{3}^{3}\,\allowbreak{}x_{4}^{12}\,\allowbreak{}x_{5}\,\allowbreak{}x_{6}^{12}}
\appmon{487}{x_{1}^{7}\,\allowbreak{}x_{2}\,\allowbreak{}x_{3}^{3}\,\allowbreak{}x_{4}^{12}\,\allowbreak{}x_{5}^{5}\,\allowbreak{}x_{6}^{8}}
\appmon{488}{x_{1}^{7}\,\allowbreak{}x_{2}\,\allowbreak{}x_{3}^{3}\,\allowbreak{}x_{4}^{12}\,\allowbreak{}x_{5}^{8}\,\allowbreak{}x_{6}^{5}}
\appmon{489}{x_{1}^{7}\,\allowbreak{}x_{2}\,\allowbreak{}x_{3}^{3}\,\allowbreak{}x_{4}^{3}\,\allowbreak{}x_{5}^{12}\,\allowbreak{}x_{6}^{10}}
\appmon{490}{x_{1}^{7}\,\allowbreak{}x_{2}\,\allowbreak{}x_{3}^{3}\,\allowbreak{}x_{4}^{5}\,\allowbreak{}x_{5}^{12}\,\allowbreak{}x_{6}^{8}}
\appmon{491}{x_{1}^{7}\,\allowbreak{}x_{2}\,\allowbreak{}x_{3}^{3}\,\allowbreak{}x_{4}^{6}\,\allowbreak{}x_{5}^{9}\,\allowbreak{}x_{6}^{10}}
\appmon{492}{x_{1}^{7}\,\allowbreak{}x_{2}\,\allowbreak{}x_{3}^{3}\,\allowbreak{}x_{4}^{8}\,\allowbreak{}x_{5}^{12}\,\allowbreak{}x_{6}^{5}}
\appmon{493}{x_{1}^{7}\,\allowbreak{}x_{2}\,\allowbreak{}x_{3}^{3}\,\allowbreak{}x_{4}^{8}\,\allowbreak{}x_{5}^{5}\,\allowbreak{}x_{6}^{12}}
\appmon{494}{x_{1}^{7}\,\allowbreak{}x_{2}\,\allowbreak{}x_{3}^{3}\,\allowbreak{}x_{4}^{8}\,\allowbreak{}x_{5}^{8}\,\allowbreak{}x_{6}^{9}}
\appmon{495}{x_{1}^{7}\,\allowbreak{}x_{2}\,\allowbreak{}x_{3}^{3}\,\allowbreak{}x_{4}^{8}\,\allowbreak{}x_{5}^{9}\,\allowbreak{}x_{6}^{8}}
\appmon{496}{x_{1}^{7}\,\allowbreak{}x_{2}\,\allowbreak{}x_{3}^{3}\,\allowbreak{}x_{4}^{9}\,\allowbreak{}x_{5}^{8}\,\allowbreak{}x_{6}^{8}}
\appmon{497}{x_{1}^{7}\,\allowbreak{}x_{2}\,\allowbreak{}x_{3}^{6}\,\allowbreak{}x_{4}^{3}\,\allowbreak{}x_{5}^{9}\,\allowbreak{}x_{6}^{10}}
\appmon{498}{x_{1}^{7}\,\allowbreak{}x_{2}^{3}\,\allowbreak{}x_{3}\,\allowbreak{}x_{4}\,\allowbreak{}x_{5}^{8}\,\allowbreak{}x_{6}^{16}}
\appmon{499}{x_{1}^{7}\,\allowbreak{}x_{2}^{3}\,\allowbreak{}x_{3}\,\allowbreak{}x_{4}^{12}\,\allowbreak{}x_{5}^{5}\,\allowbreak{}x_{6}^{8}}
\appmon{500}{x_{1}^{7}\,\allowbreak{}x_{2}^{3}\,\allowbreak{}x_{3}\,\allowbreak{}x_{4}^{12}\,\allowbreak{}x_{5}^{8}\,\allowbreak{}x_{6}^{5}}
\appmon{501}{x_{1}^{7}\,\allowbreak{}x_{2}^{3}\,\allowbreak{}x_{3}\,\allowbreak{}x_{4}^{5}\,\allowbreak{}x_{5}^{12}\,\allowbreak{}x_{6}^{8}}
\appmon{502}{x_{1}^{7}\,\allowbreak{}x_{2}^{3}\,\allowbreak{}x_{3}\,\allowbreak{}x_{4}^{5}\,\allowbreak{}x_{5}^{8}\,\allowbreak{}x_{6}^{12}}
\appmon{503}{x_{1}^{7}\,\allowbreak{}x_{2}^{3}\,\allowbreak{}x_{3}\,\allowbreak{}x_{4}^{8}\,\allowbreak{}x_{5}^{12}\,\allowbreak{}x_{6}^{5}}
\appmon{504}{x_{1}^{7}\,\allowbreak{}x_{2}^{3}\,\allowbreak{}x_{3}\,\allowbreak{}x_{4}^{8}\,\allowbreak{}x_{5}^{5}\,\allowbreak{}x_{6}^{12}}
\appmon{505}{x_{1}^{7}\,\allowbreak{}x_{2}^{3}\,\allowbreak{}x_{3}\,\allowbreak{}x_{4}^{8}\,\allowbreak{}x_{5}^{8}\,\allowbreak{}x_{6}^{9}}
\appmon{506}{x_{1}^{7}\,\allowbreak{}x_{2}^{3}\,\allowbreak{}x_{3}\,\allowbreak{}x_{4}^{8}\,\allowbreak{}x_{5}^{9}\,\allowbreak{}x_{6}^{8}}
\appmon{507}{x_{1}^{7}\,\allowbreak{}x_{2}^{3}\,\allowbreak{}x_{3}\,\allowbreak{}x_{4}^{9}\,\allowbreak{}x_{5}^{8}\,\allowbreak{}x_{6}^{8}}
\appmon{508}{x_{1}^{7}\,\allowbreak{}x_{2}^{3}\,\allowbreak{}x_{3}^{3}\,\allowbreak{}x_{4}^{3}\,\allowbreak{}x_{5}^{12}\,\allowbreak{}x_{6}^{8}}
\appmon{509}{x_{1}^{7}\,\allowbreak{}x_{2}^{3}\,\allowbreak{}x_{3}^{3}\,\allowbreak{}x_{4}^{3}\,\allowbreak{}x_{5}^{8}\,\allowbreak{}x_{6}^{12}}
\appmon{510}{x_{1}^{7}\,\allowbreak{}x_{2}^{3}\,\allowbreak{}x_{3}^{3}\,\allowbreak{}x_{4}^{4}\,\allowbreak{}x_{5}^{9}\,\allowbreak{}x_{6}^{10}}
\appmon{511}{x_{1}^{7}\,\allowbreak{}x_{2}^{3}\,\allowbreak{}x_{3}^{3}\,\allowbreak{}x_{4}^{5}\,\allowbreak{}x_{5}^{10}\,\allowbreak{}x_{6}^{8}}
\appmon{512}{x_{1}^{7}\,\allowbreak{}x_{2}^{3}\,\allowbreak{}x_{3}^{3}\,\allowbreak{}x_{4}^{5}\,\allowbreak{}x_{5}^{12}\,\allowbreak{}x_{6}^{6}}
\appmon{513}{x_{1}^{7}\,\allowbreak{}x_{2}^{3}\,\allowbreak{}x_{3}^{3}\,\allowbreak{}x_{4}^{8}\,\allowbreak{}x_{5}^{3}\,\allowbreak{}x_{6}^{12}}
\appmon{514}{x_{1}^{7}\,\allowbreak{}x_{2}^{3}\,\allowbreak{}x_{3}^{3}\,\allowbreak{}x_{4}^{8}\,\allowbreak{}x_{5}^{5}\,\allowbreak{}x_{6}^{10}}
\appmon{515}{x_{1}^{7}\,\allowbreak{}x_{2}^{3}\,\allowbreak{}x_{3}^{4}\,\allowbreak{}x_{4}^{3}\,\allowbreak{}x_{5}^{9}\,\allowbreak{}x_{6}^{10}}
\appmon{516}{x_{1}^{7}\,\allowbreak{}x_{2}^{3}\,\allowbreak{}x_{3}^{5}\,\allowbreak{}x_{4}^{10}\,\allowbreak{}x_{5}^{5}\,\allowbreak{}x_{6}^{6}}
\appmon{517}{x_{1}^{7}\,\allowbreak{}x_{2}^{3}\,\allowbreak{}x_{3}^{5}\,\allowbreak{}x_{4}^{10}\,\allowbreak{}x_{5}^{6}\,\allowbreak{}x_{6}^{5}}
\appmon{518}{x_{1}^{7}\,\allowbreak{}x_{2}^{3}\,\allowbreak{}x_{3}^{5}\,\allowbreak{}x_{4}^{3}\,\allowbreak{}x_{5}^{10}\,\allowbreak{}x_{6}^{8}}
\appmon{519}{x_{1}^{7}\,\allowbreak{}x_{2}^{3}\,\allowbreak{}x_{3}^{5}\,\allowbreak{}x_{4}^{3}\,\allowbreak{}x_{5}^{12}\,\allowbreak{}x_{6}^{6}}
\appmon{520}{x_{1}^{7}\,\allowbreak{}x_{2}^{3}\,\allowbreak{}x_{3}^{5}\,\allowbreak{}x_{4}^{3}\,\allowbreak{}x_{5}^{6}\,\allowbreak{}x_{6}^{12}}
\appmon{521}{x_{1}^{7}\,\allowbreak{}x_{2}^{3}\,\allowbreak{}x_{3}^{5}\,\allowbreak{}x_{4}^{3}\,\allowbreak{}x_{5}^{8}\,\allowbreak{}x_{6}^{10}}
\appmon{522}{x_{1}^{7}\,\allowbreak{}x_{2}^{3}\,\allowbreak{}x_{3}^{5}\,\allowbreak{}x_{4}^{5}\,\allowbreak{}x_{5}^{8}\,\allowbreak{}x_{6}^{8}}
\appmon{523}{x_{1}^{7}\,\allowbreak{}x_{2}^{3}\,\allowbreak{}x_{3}^{5}\,\allowbreak{}x_{4}^{6}\,\allowbreak{}x_{5}^{10}\,\allowbreak{}x_{6}^{5}}
\appmon{524}{x_{1}^{7}\,\allowbreak{}x_{2}^{3}\,\allowbreak{}x_{3}^{5}\,\allowbreak{}x_{4}^{6}\,\allowbreak{}x_{5}^{3}\,\allowbreak{}x_{6}^{12}}
\appmon{525}{x_{1}^{7}\,\allowbreak{}x_{2}^{3}\,\allowbreak{}x_{3}^{5}\,\allowbreak{}x_{4}^{6}\,\allowbreak{}x_{5}^{6}\,\allowbreak{}x_{6}^{9}}
\appmon{526}{x_{1}^{7}\,\allowbreak{}x_{2}^{3}\,\allowbreak{}x_{3}^{5}\,\allowbreak{}x_{4}^{6}\,\allowbreak{}x_{5}^{9}\,\allowbreak{}x_{6}^{6}}
\appmon{527}{x_{1}^{7}\,\allowbreak{}x_{2}^{3}\,\allowbreak{}x_{3}^{5}\,\allowbreak{}x_{4}^{8}\,\allowbreak{}x_{5}\,\allowbreak{}x_{6}^{12}}
\appmon{528}{x_{1}^{7}\,\allowbreak{}x_{2}^{3}\,\allowbreak{}x_{3}^{5}\,\allowbreak{}x_{4}^{8}\,\allowbreak{}x_{5}^{3}\,\allowbreak{}x_{6}^{10}}
\appmon{529}{x_{1}^{7}\,\allowbreak{}x_{2}^{3}\,\allowbreak{}x_{3}^{5}\,\allowbreak{}x_{4}^{9}\,\allowbreak{}x_{5}^{6}\,\allowbreak{}x_{6}^{6}}
\appmon{530}{x_{1}^{7}\,\allowbreak{}x_{2}^{3}\,\allowbreak{}x_{3}^{8}\,\allowbreak{}x_{4}^{5}\,\allowbreak{}x_{5}\,\allowbreak{}x_{6}^{12}}
\appmon{531}{x_{1}^{7}\,\allowbreak{}x_{2}^{3}\,\allowbreak{}x_{3}^{8}\,\allowbreak{}x_{4}^{5}\,\allowbreak{}x_{5}^{3}\,\allowbreak{}x_{6}^{10}}
\appmon{532}{x_{1}^{7}\,\allowbreak{}x_{2}^{3}\,\allowbreak{}x_{3}^{9}\,\allowbreak{}x_{4}\,\allowbreak{}x_{5}^{8}\,\allowbreak{}x_{6}^{8}}
\appmon{533}{x_{1}^{7}\,\allowbreak{}x_{2}^{3}\,\allowbreak{}x_{3}^{9}\,\allowbreak{}x_{4}^{5}\,\allowbreak{}x_{5}^{6}\,\allowbreak{}x_{6}^{6}}
\appmon{534}{x_{1}^{7}\,\allowbreak{}x_{2}^{3}\,\allowbreak{}x_{3}^{9}\,\allowbreak{}x_{4}^{6}\,\allowbreak{}x_{5}^{5}\,\allowbreak{}x_{6}^{6}}
\appmon{535}{x_{1}^{7}\,\allowbreak{}x_{2}^{3}\,\allowbreak{}x_{3}^{9}\,\allowbreak{}x_{4}^{6}\,\allowbreak{}x_{5}^{6}\,\allowbreak{}x_{6}^{5}}
\appmon{536}{x_{1}^{7}\,\allowbreak{}x_{2}^{8}\,\allowbreak{}x_{3}^{3}\,\allowbreak{}x_{4}^{5}\,\allowbreak{}x_{5}^{3}\,\allowbreak{}x_{6}^{10}}
\appmon{537}{x_{1}^{7}\,\allowbreak{}x_{2}^{9}\,\allowbreak{}x_{3}^{3}\,\allowbreak{}x_{4}^{5}\,\allowbreak{}x_{5}^{6}\,\allowbreak{}x_{6}^{6}}
\appmon{538}{x_{1}^{7}\,\allowbreak{}x_{2}^{9}\,\allowbreak{}x_{3}^{3}\,\allowbreak{}x_{4}^{6}\,\allowbreak{}x_{5}^{5}\,\allowbreak{}x_{6}^{6}}
\appmon{539}{x_{1}^{7}\,\allowbreak{}x_{2}^{9}\,\allowbreak{}x_{3}^{3}\,\allowbreak{}x_{4}^{6}\,\allowbreak{}x_{5}^{6}\,\allowbreak{}x_{6}^{5}}
\end{multicols}
\endgroup

\section{Software archive and computational source files}\label{app:software}

The research code associated with this paper is organized as a downloadable Zenodo archive:
\begin{center}
\SoftwareZenodo.
\end{center}
The archive contains the following source files and drivers.

\begin{center}
\small
\begin{tabular}{>{\raggedright\arraybackslash}p{0.42\linewidth}p{0.48\linewidth}}
\toprule
File & Purpose \\
\midrule
\path{src/AlgebraicTransfer.jl} & A new \texttt{Julia} package built on \texttt{Nemo} and \texttt{AbstractAlgebra}, providing exponent-vector arithmetic, streamed hit reduction, Kameko matrices, bit-packed nullspaces over $\mathbb{F}_2$, invariant solvers, and orbit-type post-processing.\\
\path{src/AlgebraicTransfer_OSCAR_Algorithm1_full_reference.jl} & Direct \texttt{OSCAR} implementation of Algorithm~\ref{alg:main}, including the production routines for the degree-$15$ target, the degree-$36$ Kameko kernel, lower-weight correction, target-lift correction, and orbit-type tables. \\
\path{drivers/run_singer_6_15.jl} & Recomputes the degree-$15$ invariant space $[(QP_6)_{15}]^{GL(6)}$ and the representative $\xi$. \\
\path{drivers/run_singer_6_36.jl} & Recomputes the degree-$36$ Kameko kernel, the correction systems, and the two representatives $\zeta_1,\zeta_2$. \\
\path{drivers/run_dickson_validation.jl} & \texttt{Julia}/\texttt{OSCAR} driver producing the Dickson validation data in Tables~\ref{tab:dickson-q3-validation} and~\ref{tab:dickson-q4-validation}. \\
\path{drivers/run_permutation_module_validation.jl} & \texttt{Julia}/\texttt{OSCAR} driver producing the finite permutation-module validation data in Table~\ref{tab:permutation-module-validation}. \\
\path{drivers/postprocess_orbit_tables.jl} & Parses reduced representatives and produces the support-slice and stabilizer tables. \\
\path{python/run_dickson_validation.py} and \path{python/run_permutation_module_validation.py} & Lightweight independent validation scripts for the Dickson and finite-group examples. \\
\bottomrule
\end{tabular}
\end{center}



\begin{thebibliography}{99}

\bibitem{Boardman}
J. M. Boardman,
\textit{Modular representations on the homology of powers of real projective space}, in \textit{Algebraic Topology: Oaxtepec 1991}, Contemporary Mathematics, vol. 146, American Mathematical Society, Providence, RI, 1993, pp. 49--70, DOI: \url{https://doi.org/10.1090/conm/146/01215}.



\bibitem{Chen}
T. W. Chen,
\textit{Determination of $\operatorname{Ext}^{5,*}_{\mathscr A}(\mathbb Z/2,\mathbb Z/2)$}, Topology and its Applications \textbf{158} (2011), no. 5, 660--689, DOI: \url{https://doi.org/10.1016/j.topol.2011.01.002}.

\bibitem{Chen2}
T. W. Chen,
\textit{The structure of decomposable elements in $\operatorname{Ext}^{6,*}_{\mathscr A}(\mathbb Z/2,\mathbb Z/2)$}, preprint, 2012, 35 pp.

\bibitem{Chen3}
T.W. Chen, 
\textit{Indecomposable elements in ${\rm Ext}^{6,*}_{\mathscr{A}}(\mathbb{Z}/2, \mathbb{Z}/2),$} preprint, 2013, 3 pp.

\bibitem{CrossleyTurgay}
M. Crossley and N.D. Turgay,
\textit{Conjugation invariants in the Leibniz--Hopf algebra}, J. Pure Appl. Algebra \textbf{217} (2013), 2247--2254.

\bibitem{OSCAR2025}
W. Decker, C. Eder, C. Fieker, M. Horn, M. Joswig (eds.), \textit{The Computer Algebra System OSCAR: Algorithms and Examples}, Algorithms and Computation in Mathematics, Vol.~32, Springer, Cham, 2025. \url{https://doi.org/10.1007/978-3-031-62127-7}

\bibitem{Dickson1911}
L. E. Dickson,
\textit{A fundamental system of invariants of the general modular linear group with a solution of the form problem}, Transactions of the American Mathematical Society \textbf{12} (1911), no. 1, 75--98, DOI: \url{https://doi.org/10.1090/S0002-9947-1911-1500882-4}.

\bibitem{Dold}
A. Dold,
\textit{Erzeugende der Thomschen Algebra $\mathfrak N$},
Mathematische Zeitschrift \textbf{65} (1956), no. 1, 25--35, DOI: \url{https://doi.org/10.1007/BF01473868}.

\bibitem{Ha}
L. M. H\`a, \textit{Sub-Hopf algebras of the Steenrod algebra and the Singer transfer}, Geom. Topol. Monogr. \textbf{11} (2007), 101-124. DOI: \url{https://doi.org/10.2140/gtm.2007.11.81}.

\bibitem{Hung}
N. H. V. Hung,
\textit{The cohomology of the Steenrod algebra and representations of the general linear groups}, Transactions of the American Mathematical Society \textbf{357} (2005), no. 10, 4065--4089, DOI: \url{https://doi.org/10.1090/S0002-9947-05-03889-4}.

\bibitem{Hung2}
N. H. V. Hung and V. T. N. Quynh,
\textit{The image of Singer's fourth transfer}, Comptes Rendus Math\'ematique \textbf{347} (2009), no. 23--24, 1415--1418, DOI: \url{https://doi.org/10.1016/j.crma.2009.10.018}.

\bibitem{HungSquares}
N. H. V. Hung,
\textit{The action of the Steenrod squares on the modular invariants of linear groups}, Proceedings of the American Mathematical Society \textbf{113} (1991), no. 4, 1097--1104, DOI: \url{https://doi.org/10.2307/2048789}.

\bibitem{HungMinh1995}
N. H. V. Hung and P. A. Minh,
\textit{The action of the mod $p$ Steenrod operations on the modular invariants of linear groups}, Vietnam Journal of Mathematics \textbf{23} (1995), no. 1, 39--56.

\bibitem{HungPetersonDickson}
N. H. V. Hung and F. P. Peterson,
\textit{$\mathscr A$-generators for the Dickson algebra}, Transactions of the American Mathematical Society \textbf{347} (1995), no. 12, 4687--4728, DOI: \url{https://doi.org/10.1090/S0002-9947-1995-1316852-X}.

\bibitem{HungPetersonSpherical}
N. H. V. Hung and F. P. Peterson,
\textit{Spherical classes and the Dickson algebra}, Mathematical Proceedings of the Cambridge Philosophical Society \textbf{124} (1998), no. 2, 253--264, DOI: \url{https://doi.org/10.1017/S0305004198002667}.

\bibitem{HungNamModular}
N. H. V. Hung and T. N. Nam,
\textit{The hit problem for modular invariants of linear groups}, Journal of Algebra \textbf{246} (2001), no. 1, 367--384, DOI: \url{https://doi.org/10.1006/jabr.2001.8974}.

\bibitem{HungNamDickson}
N. H. V. Hung and T. N. Nam,
\textit{The hit problem for the Dickson algebra}, Transactions of the American Mathematical Society \textbf{353} (2001), no. 12, 5029--5040, DOI: \url{https://doi.org/10.1090/S0002-9947-01-02705-2}.

\bibitem{HungTriviality}
N. H. V. Hung,
\textit{On triviality of Dickson invariants in the homology of the Steenrod algebra}, Mathematical Proceedings of the Cambridge Philosophical Society \textbf{134} (2003), no. 1, 103--113, DOI: \url{https://doi.org/10.1017/S0305004102006187}.

\bibitem{HungMargolis2}
N. H. V. Hung,
\textit{The mod 2 Margolis homology of the Dickson algebra}, Comptes Rendus Math\'ematique \textbf{358} (2020), no. 4, 505--510, DOI: \url{https://doi.org/10.5802/crmath.68}.

\bibitem{HungMargolisP}
N. H. V. Hung,
\textit{The mod $p$ Margolis homology of the Dickson--M\`ui algebra}, Comptes Rendus Math\'ematique \textbf{359} (2021), no. 3, 229--236, DOI: \url{https://doi.org/10.5802/crmath.151}.

\bibitem{Hung3}
N.H.V. H\uhorn ng,
\textit{Images of the Singer transfers and their possibility to be injective}, Journal of Mathematics and Mathematical Sciences \textbf{4} (2025), no. 1, 95-103, URL: \url{https://science.thanglong.edu.vn/index.php/volc/article/view/235}.

\bibitem{Janfada}
A. S. Janfada,
\textit{A criterion for a monomial in $P(3)$ to be hit}, Mathematical Proceedings of the Cambridge Philosophical Society \textbf{145} (2008), no. 3, 587--599, DOI: \url{https://doi.org/10.1017/S0305004108001515}.

\bibitem{Kameko}
M. Kameko,
\textit{Products of projective spaces as Steenrod modules}, Ph.D. thesis, The Johns Hopkins University, ProQuest LLC, Ann Arbor, MI, 1990, 29 pp.

\bibitem{Karaca}
I. Karaca,
\textit{Unstable modules over the Steenrod algebra}, Journal of Pure and Applied Algebra \textbf{171} (2002), 249-255.


\bibitem{Lin}
W. H. Lin,
\textit{$\operatorname{Ext}^{4,*}_{\mathscr A}(\mathbb Z/2,\mathbb Z/2)$ and $\operatorname{Ext}^{5,*}_{\mathscr A}(\mathbb Z/2,\mathbb Z/2)$}, Topology and its Applications \textbf{155} (2008), no. 5, 459--496, DOI: \url{https://doi.org/10.1016/j.topol.2007.11.003}.

\bibitem{Lin2}
W. Lin,
\textit{Charts of the cohomology of the mod 2 Steenrod algebra}, Zenodo, 2023, DOI: \url{https://doi.org/10.5281/zenodo.7786290}.

\bibitem{Minami}
N. Minami,
\textit{The iterated transfer analogue of the new doomsday conjecture}, Transactions of the American Mathematical Society \textbf{351} (1999), no. 6, 2325--2351, DOI: \url{https://doi.org/10.1090/S0002-9947-99-02037-1}.

\bibitem{Mui1975}
H. M\`ui,
\textit{Modular invariant theory and cohomology algebras of symmetric groups}, Journal of the Faculty of Science, University of Tokyo, Section IA, Mathematics \textbf{22} (1975), 319--369.

\bibitem{Nam}
T. N. Nam,
\textit{Transfert alg\'ebrique et action du groupe lin\'eaire sur les puissances divis\'ees modulo 2}, Annales de l'Institut Fourier (Grenoble) \textbf{58} (2008), no. 5, 1785--1837, DOI: \url{https://doi.org/10.5802/aif.2399}.

\bibitem{Oscar}
The OSCAR Development Team,
\textit{OSCAR -- Open Source Computer Algebra System}, \url{https://www.oscar-system.org/}.


\bibitem{Peterson}
F. P. Peterson,
\textit{Generators of $H^*(\mathbb RP^{\infty}\times \mathbb RP^{\infty})$ as a module over the Steenrod algebra}, Abstracts of Papers Presented to the American Mathematical Society \textbf{833} (1987), 55--89.

\bibitem{Phuc}
\DJ.V. Ph\'uc,
\textit{The affirmative answer to Singer's conjecture on the algebraic transfer of rank four}, Proceedings of the Royal Society of Edinburgh, Section A: Mathematics \textbf{153} (2023), no. 5, 1529--1542, DOI: \url{https://doi.org/10.1017/prm.2022.57}. Full online version (24 pages) available at OSF:\url{https://osf.io/preprints/osf/4ckf8_v9}.

\bibitem{Phuc3}
\DJ.V. Ph\'uc,
\textit{A note on the hit problem for the polynomial algebra of six variables and the sixth algebraic transfer}, Journal of Algebra \textbf{613} (2023), 1--31, DOI: \url{https://doi.org/10.1016/j.jalgebra.2022.08.028}.

\bibitem{Phuc7}
\DJ.V. Ph\'uc,
\textit{Computational Approaches to the Singer Transfer: Preimages in the Lambda Algebra and $G_k$-Invariant Theory}, Preprint, 2025, 100 pages, arXiv:2507.10108v4, DOI: \url{https://doi.org/10.48550/arXiv.2507.10108}.

\bibitem{P8}
\DJ.V. Ph\'uc, \textit{On Singer's conjecture for the fourth algebraic transfer in certain generic degrees}, Preprint (2025), 34 pages, arXiv:2506.10232, \url{https://doi.org/10.48550/arXiv.2506.10232}.

\bibitem{P9}
\DJ.V. Ph\'uc, \textit{On the algebraic transfers of ranks 4 and 6 at generic degrees}, Preprint (2026), 93 pages, Available online at OSF: \url{https://doi.org/10.31219/osf.io/6tkxh_v10}.

\bibitem{ZenodoDeg15}
\DJ.V. Ph\'uc,
\textit{Detailed computational output for the case $(q,n)=(6,15)$}, Zenodo, 2026, DOI: \url{https://doi.org/10.5281/zenodo.17620661}.

\bibitem{ZenodoDeg36}
\DJ.V. Ph\'uc,
\textit{Detailed computational output for the case $(q,n)=(6,36)$}, Zenodo, 2026, DOI: \url{https://doi.org/10.5281/zenodo.17620680}.

\bibitem{Phuc2026}
\DJ.V. Ph\'uc, \textit{On a counterexample to Singer's conjecture for the fifth algebraic transfer in degree 108}, Preprint (2026), 18 pages. Available at Zenodo \url{https://doi.org/10.5281/zenodo.22750602}.

\bibitem{PS}
\DJ.V. Ph\'uc and N. Sum,
\textit{On the generators of the polynomial algebra as a module over the Steenrod algebra}, Comptes Rendus Math\'ematique \textbf{353} (2015), no. 11, 1035--1040, DOI: \url{https://doi.org/10.1016/j.crma.2015.09.003}.

\bibitem{PS2}
\DJ.V. Ph\'uc and N. Sum,
\textit{On a minimal set of generators for the polynomial algebra of five variables as a module over the Steenrod algebra}, Acta Mathematica Vietnamica \textbf{42} (2017), no. 1, 149--162, DOI: \url{https://doi.org/10.1007/s40306-016-0194-6}.

\bibitem{Serre}
J.-P. Serre,
\textit{Linear Representations of Finite Groups}, Graduate Texts in Mathematics \textbf{42}, Springer, New York, 1977, DOI: \url{https://doi.org/10.1007/978-1-4684-9458-7}.

\bibitem{Singer}
W. M. Singer,
\textit{The transfer in homological algebra}, Mathematische Zeitschrift \textbf{202} (1989), no. 4, 493--523, DOI: \url{https://doi.org/10.1007/BF01221587}.

\bibitem{Sum1}
N. Sum,
\textit{The squaring operation and the Singer algebraic transfer}, Vietnam Journal of Mathematics \textbf{49} (2021), no. 4, 1079--1096, DOI: \url{https://doi.org/10.1007/s10013-020-00423-1}.

\bibitem{Sum2}
N. Sum, \textit{A counter-example to Singer's conjecture for the algebraic transfer}, Preprint (2025), arXiv:2408.06669.

\bibitem{Thom}
R. Thom,
\textit{Quelques propri\'et\'es globales des vari\'et\'es diff\'erentiables},
Commentarii Mathematici Helvetici \textbf{28} (1954), no. 1, 17--86, DOI: \url{https://doi.org/10.1007/BF02566923}.

\bibitem{Turgay1}
N. D. Turgay,  
\textit{An alternative approach to the Adem relations in the mod 2 Steenrod algebra}, Turkish Journal of Mathematics \textbf{38} (2014), 924-934.

\bibitem{Turgay2}
N. D. Turgay,  
\textit{A remark on the conjugation in the Steenrod algebra}, Communications of the Korean Mathematical Society \textbf{30} (2015), 269-276.

\bibitem{Turgay3}
N. D. Turgay and S. Kaji,
\textit{The mod 2 dual Steenrod algebra as a subalgebra of the mod 2 dual Leibniz--Hopf algebra},  Journal of Homotopy and Related Structures  \textbf{12} (2017), 727-739.

\bibitem{Turgay2020}
N. D. Turgay,
\textit{On the mod $p$ Steenrod algebra and the Leibniz--Hopf algebra}, Electronic Research Archive \textbf{28} (2020),  951--959.

\bibitem{TK}
N. D. Turgay and I. Karaca,  
\textit{The Arnon bases in the Steenrod algebra}, Georgian Mathematical Journal \textbf{27} (2020), 649-654.

\bibitem{WW}
G. Walker and R. M. W. Wood,
\textit{Polynomials and the mod 2 Steenrod Algebra. Volume 1: The Peterson hit problem}, London Mathematical Society Lecture Note Series, vol. 441, Cambridge University Press, Cambridge, 2018, DOI: \url{https://doi.org/10.1017/9781108333368}.

\bibitem{WW2}
G. Walker and R. M. W. Wood,
\textit{Polynomials and the mod 2 Steenrod Algebra. Volume 2: Representations of $GL(n,\mathbb F_2)$}, London Mathematical Society Lecture Note Series, vol. 442, Cambridge University Press, Cambridge, 2018, DOI: \url{https://doi.org/10.1017/9781108304092}.

\bibitem{Wood2}
R. M. W. Wood,
\textit{Steenrod squares of polynomials and the Peterson conjecture}, Mathematical Proceedings of the Cambridge Philosophical Society \textbf{105} (1989), no. 2, 307--309, DOI: \url{https://doi.org/10.1017/S0305004100067797}.

\end{thebibliography}
\end{document}